\documentclass[12pt]{amsart}

\usepackage{amsmath,amsthm,amssymb, amscd, amsfonts}

\numberwithin{equation}{section}

\newtheorem{theorem}{Theorem}[section]
\newtheorem{proposition}[theorem]{Proposition}
\newtheorem{conjecture}[theorem]{Conjecture}
\newtheorem{corollary}[theorem]{Corollary}
\newtheorem{lemma}[theorem]{Lemma}
\newtheorem{claim}[theorem]{Claim}
\newtheorem{maintheorem}[theorem]{Main Theorem}

\theoremstyle{definition}

\newtheorem{remark}[theorem]{Remark}
\newtheorem{example}[theorem]{Example}
\newtheorem{definition}[theorem]{Definition}

\renewcommand{\eqref}[1]{{\rm (\ref{#1})}}

\def\endproof{\hfill$\square$\medskip}

\newcommand {\osr}{\overset\sim\rightarrow}

\def\AA{\mathbb{A}}

\def\ZZ{\mathbb{Z}}
\def\CC{\mathbb{C}}

\def\RR{\mathbb{R}}
\def\QQ{\mathbb{Q}}
\def\GG{\mathbb{G}}

\def\BB{\mathcal{B}}
\def\XX{\mathcal{X}}

\def\Supp{\operatorname{Supp}}
\def\Norm{\operatorname{Norm}}
\def\Trop{\operatorname{Trop}}
\def\Hom{\operatorname{Hom}}
\def\dom{\operatorname{dom}}
\def\ran{\operatorname{ran}}
\def\Vert{\operatorname{Vert}}
\def\reg{\operatorname{reg}}
\def\dual{\operatorname{dual}}
\def\id{\operatorname{id}}
\def\Fun{\operatorname{Fun}}

\def\bb{\mathfrak{b}}
\def\ll{\mathfrak{l}}
\def\gg{\mathfrak{g}}
\def\hh{\mathfrak{h}}
\def\nn{\mathfrak{n}}

\def\uu{\mathfrak{u}}

\def\ii{\mathbf{i}}

\def\BBB{\mathfrak{B}}

\addtocounter{section}{-1}

\begin{document}

\title[Geometric and Unipotent Crystals II]{Geometric and Unipotent Crystals II:\\
 From Unipotent Bicrystals to Crystal Bases}

\author[A.~BERENSTEIN]{Arkady Berenstein}
\address{Department of Mathematics, University of Oregon,
Eugene, OR 97403, USA} \email{arkadiy@math.uoregon.edu}

\author[D.~KAZHDAN]{David Kazhdan}
\address{\noindent Department of Mathematics, Hebrew University, Jerusalem, Israel}
\email{kazhdan@math.huji.ac.il}


\dedicatory{To the memory of Joseph Donin}

\date{January 15, 2006}

\thanks{The first author was supported in part
by NSF grants  and the second author was supported in part by ISF
and NSF grants. }



\begin{abstract}

For each reductive algebraic group $G$, we introduce and study
{\it unipotent bicrystals} which serve as a regular version of
birational {\it geometric and unipotent crystals} introduced
earlier by the authors. The framework of unipotent bicrystals
allows, on the one hand, to study systematically such varieties as
Bruhat cells in $G$ and their convolution products and, on the
other hand, to give a new construction of many normal Kashiwara
crystals including those for $G^\vee$-modules, where $G^\vee$ is
the Langlands dual groups. In fact, our  analogues of crystal
bases (which we refer to as crystals {\it associated} to
$G^\vee$-modules) are associated to $G^\vee$-modules directly,
i.e., without quantum deformations.

One of the main  results of the present paper is an explicit construction  of the  crystal  $\BB_0$ for the coordinate ring
of the dual flag variety $X^\vee_0=G^\vee/U^\vee$ based on the
{\it positive} unipotent bicrystal on the open Bruhat cell $X_0=Bw_0B$.  Our general {\it tropicalization} procedure assigns to each
{\it strongly positive} unipotent bicrystal a normal Kashiwara crystal $\BB$ equipped with the {\it multiplicity erasing} homomorphism
$\BB\to \BB_0$ and the {\it combinatorial central charge}  $\BB\to \ZZ$ which
is invariant under all crystal operators. Applying the construction to $\BB_0\times \BB_0$ gives a crystal multiplication
$\BB_0\times \BB_0\to \BB_0$ and an invariant  grading $\BB_0\times \BB_0\to \ZZ$.

\end{abstract}

\maketitle


\section{Introduction}
\label{sec:introduction}

The present paper is a continuation  of the study of geometric and unipotent crystals  initiated in \cite{bk}.
However, all necessary definitions and constructions are included, so that the paper can be read independently.

The aim of this paper is two-fold: first, to introduce {\it unipotent bicrystals} as regular versions of geometric and unipotent crystals
from \cite{bk} and, second, to construct a large class of Kashiwara's crystal bases (or,  rather, the combinatorial crystals directly {\it associated} to appropriate modules) via the {\it tropicalization} of {\it positive} unipotent bicrystals.

More precisely, let $G$ be a split reductive algebraic group and
$B\subset G$ be a Borel subgroup. Denote by  $G^\vee$ the
Langlands dual group of $G$, and by $\gg^\vee$  the Lie algebra of
$G^\vee$. For each dominant integral weight $\lambda$ of
$\gg^\vee$, denote by $V_\lambda$ the irreducible
$\gg^\vee$-module and by $\BB(V_\lambda)$ the corresponding normal
Kashiwara crystal (see e.g., \cite{k93}). We will  explicitly
construct all $\BB(V_\lambda)$ in terms of the positive unipotent
bicrystal structure on the open Bruhat cell $X=Bw_0B$ or, which is
the same, we construct the {\it associated crystal}
$\BB_0=\bigsqcup_\lambda \BB(V_\lambda)$ for the coordinate
algebra of the basic affine space $X^\vee_0=G^\vee/U^\vee$
(Theorem \ref{th:Blambda}).

More generally, for each positive unipotent bicrystal on a
$U\times U$-variety $X$ (where $U$ is the unipotent radical of a
Borel subgroup $B$ of $G$), we construct an infinite normal
Kashiwara crystal $\BB$ which is a union of finite normal crystals
$\BB^\lambda$ (Section \ref{subsect:from positive unipotent
bicrystals to combinatorial crystals}). If  the unipotent
bicrystal were {\it strongly positive}, then the resulting normal
crystal gets equipped with {\it multiplicity erasing map} $\tilde
{\bf m}:\BB\to \BB_0$ such that $\tilde {\bf
m}(\BB^\lambda)=\BB(V_\lambda)$ if $\BB^\lambda$ is not empty
(Corollary \ref{cor:projection tilde m}) and with the {\it
combinatorial central charge} $\tilde \Delta:\BB\to \ZZ$ invariant
under all crystal operators $\tilde e_i^{\,n}$ (Claim
\ref{cl:tilde delta}). In particular, if $\BB=\BB_0\times \BB_0$
is  the crystal associated to the algebra $\CC[X^\vee_0\times
X^\vee_0]=\CC[X^\vee_0]\otimes \CC[X^\vee_0]$,  we obtain:

\noindent $\bullet$ the associative {\it crystal multiplication} $\tilde {\bf m}_0:\BB_0\times \BB_0\to \BB_0$. This multiplication turns
$\BB_0$ into a monoid in the category of normal Kashiwara crystals.

\noindent $\bullet$ the {\it combinatorial central charge} $\tilde
\Delta_0: \BB_0\times \BB_0\to \ZZ$ invariant under all crystal
operators $\tilde e_i^{\,n}$ acting on  $\BB_0\times
\BB_0\cong\oplus_{\lambda,\lambda'} \BB(V_\lambda\otimes
V_{\lambda'})$. It is a combinatorial analogue of a certain
operator $\Delta_0:\CC[X^\vee_0]\otimes \CC[X^\vee_0]\to
\CC[X^\vee_0]\otimes \CC[X^\vee_0]$ commuting with the
$\gg^\vee$-action. In  turn, $\tilde \Delta_0$ provides a
$q$-analog of the tensor product multiplicities:
$$[V_{\lambda''}:V_\lambda\otimes V_{\lambda'}]_q=\sum_{\tilde b} q^{\tilde \Delta_0(\tilde b)}\ ,$$
where the summation is over all  $\tilde b\in \BB(V_\lambda \otimes V_{\lambda'})^{\lambda''}=(\BB(V_\lambda) \times \BB(V_{\lambda'}))^{\lambda''}$. Similarly, we  define the $q$-tensor multiplicities for products of
several irreducible $\gg^\vee$-modules (Section \ref{subsect:Crystal multiplication and the central charge}).

This shows that strongly positive unipotent bicrystals are the
closest geometric ``relatives" of Kashiwara crystal bases of
$U_q(\gg)$-modules. Moreover, we expect (Conjecture
\ref{con:strongly positive crystal bases}) that each $\BB^\lambda$
obtained from a strongly positive unipotent bicrystal is always
isomorphic to a union of copies of $\BB(V_\lambda)$. To make the
analogy between geometric and combinatorial objects more precise,
we lift  the Kashiwara construction of the product of Kashiwara
crystals to the geometric level --  unipotent bicrystals
(similarly to the unipotent crystals introduced in \cite{bk}) form
a nice monoidal category under the convolution product (Claims
\ref{cl:monoidal} and \ref{cl:unipotent chi-linear bicrystals}),
and it turns out that strong positivity is preserved under this
product (Theorem \ref{th:positive unipotent monoidal}).

As an application of this method, for each parabolic subgroup $P$
of $G$, we construct  the unipotent bicrystal  on the {\it reduced
Bruhat cell} $U\overline {w_P}U$
which, in turn,  produces the crystal associated to $\CC[U_P]$,
the coordinate algebra  of the unipotent radical $U_P^\vee$ of the
dual parabolic $P^\vee$ (Theorem~\ref{th:asociated crystal
Schubert}).

One of the main tools of this paper is the transition from the
unipotent bicrystals and geometric crystals to Kashiwara crystals
which we refer to as the {\it tropicalization}. It is based on
Theorem \ref{th:Trop} (originally proved in our first paper
\cite{bk}) which establishes the functoriality of the transition
from the category of {\it positive algebraic tori} to the category
of (marked) sets. In Section \ref{sect:tropicalization}  we
discuss in detail the semi-field of polytopes and related
structures as a combinatorial foundation of the total positivity
and the carrier of the tropicalization functor. To extend the
total positivity to other rational varieties and unipotent
bicrystals, we introduce in Section \ref{subsect:toric charts and
positive structures} a notion of a {\it positive variety} and
discuss the properties of these new algebro-geometric objects.

In order to relate the combinatorial crystals with the actual
bases for $\gg^\vee$-modules, we develop in Section
\ref{sect:Kashiwara crystals, perfect bases, and crystal bases}
the theory of {\it perfect bases} of modules and their {\it
associated crystals}. The novelty of this approach is that our
definition of the Kashiwara crystal {\it associated} to a
$\gg^\vee$-module does not require any quantum deformation of
$\gg^\vee$ or its modules. Our main results in this direction
(Theorems \ref{th:the upper explicit crystal bijection} and
\ref{th:the crystal}) guarantee that the crystal associated to a
module is well-defined and independent of the choice of an
underlying perfect basis. Note, however, that our approach allows
to prove uniqueness, rather than existence, of the appropriate
Kashiwara crystals. The existence follows, for instance, from the
fact that the canonical, global crystal, and semi-canonical bases
are all perfect (see Remark \ref{rem:perfect canonical}).

\smallskip
A number of statements in the paper, which we refer to as {\it claims} and {\it corollaries},  are almost immediate. Their proofs  are left to the reader.

\smallskip
\noindent {\bf Acknowledgments.} We are grateful to Dr.~O.~Schramm
of  Microsoft Research and Professor G.~Kalai of Hebrew University
for very useful discussions about piecewise-linear continuous
maps. The first author would like to express his gratitude to
Professor M.~Okado of Osaka University, Professor Kuniba of
University of Tokyo, and Dr. Anatol Kirillov of RIMS for the
opportunity to present the results of this paper in RIMS in the
Summer of 2004.

\section{ $U\times U$-varieties }

\subsection{Definitions and notation}


\label{subsect:notation}

Throughout the paper, $\GG_m$ and $\GG_a$ are respectively the multiplicative and the additive groups, and $G$ is a split reductive algebraic
group defined over $\QQ$ (or over any other field of characteristic $0$).

We fix  a maximal torus $T\subset G$ and  a  Borel subgroup   $B$  of $G$ containing $T$.
Let $B^-\subset G$ be the Borel subgroup opposite to $B$, i.e., $B\cap B^-=T$.  Denote by $U$ and $U^-$
respectively the unipotent radicals of $B$ and $B^-$.

Let $X^\star(T)=\Hom(T,\GG_m)$ and $X_\star(T)=\Hom(\GG_m,T)$ be
respectively the lattices of {\it characters} and {\it
co-characters} of $T$. These lattices are also known as weight and
the co-weight lattices of $G$. By definition, the lattices are
dual to each other via the canonical pairing $\left<\cdot,\cdot
\right>:X^\star(T)\times X_\star(T)\to \ZZ$.  Denote by  $I$  the
set of vertices of the Dynkin diagram of $G$ and for any $i\in I$
denote by $\alpha_i\in X^\star(T)$ the simple root, and by
$\alpha_i^\vee\in X_\star(T)$ the corresponding simple coroot.

For each $i\in I$, we fix a group homomorphism $\phi_i:SL_2\to G$
such that
$$\phi_i\begin{pmatrix} 1 & 0 \\
                 \GG_a & 1
\end{pmatrix}\subset B^- ,~ \phi_i\begin{pmatrix} 1 & \GG_a \\
                 0 & 1
\end{pmatrix}\subset B,~\phi_i\begin{pmatrix} c & 0 \\
                 0 & c^{-1}
\end{pmatrix}=\alpha_i^\vee(c)$$
for $c\in \GG_m$. Such a simultaneous choice of homomorphisms $\phi_i$ is also called a {\it splitting} of $G$. Note that any two splittings
are conjugate by an element of $T\cap [G,G]$, where $[G,G]$ is the adjoint group of $G$.

 Using a splitting $\{\phi_i\}$, $i\in I$, we define the co-characters
$y_i:\GG_a\to U^-$ and $x_i:\GG_a\to U$ by
$$y_i(a):=\phi_i\begin{pmatrix} 1 & 0 \\
                 a & 1
\end{pmatrix}\in B^- ,~ x_i(a):=~\phi_i\begin{pmatrix} 1 & a \\
                 0 & 1
\end{pmatrix}\in B \ .$$

We denote by $U_i\subset U$  the image $x_i(\GG_a)$, by
$U^-_i\subset U^-$ the image $y_i(\GG_a)$. Clearly, $U$
(resp.~$U^-$) is generated by  $U_i$ (resp.~$U_i^-$), $i\in I$.

Denote by $\widehat U$ the set of all characters of $U$, (i.e.,
the set of group homomorphisms  $\chi:U\to \GG_a$). And for each
$i\in I$, define the elementary character $\chi_i\in \widehat U$
by
$$\chi_i(x_j(a))=\delta_{ij}\cdot a$$
for $a\in \GG_a$. The family
$\chi_i$, $i\in I$, is a basis in the vector space $\widehat U$.

Following \cite[Section 4.2]{bz-invent}, we define  the ``positive
inverse" anti-automorphism $\iota:G\to G$ by:
\begin{equation}
\label{eq:iota}
\iota(y_i(a))=y_i(a),~\iota(t)=t^{-1},\iota(x_i(a))=x_i(a)
\end{equation}
for $a\in \GG_a$, $t\in T$.

\begin{claim}
\label{cl:chi_i iota} For each character $\chi:U\to \GG_a$ and
$u\in U$, one has $$\chi(\iota(u))=\chi(u)=-\chi(u^{-1}) \ .$$

\end{claim}

Denote by $W$ the Weyl group of $G$. By definition, $W$ is
generated by the simple reflections $s_i$, $i\in I$. Let $l:W\to
\ZZ_{\ge 0}$ ($w\mapsto l(w)$) be the length function. For any
sequence $\ii=(i_1,\ldots,i_\ell)\in I^\ell$, we write
$w(\ii)=s_{i_1}\cdots s_{i_\ell}$. A sequence $\ii\in I^l$ is
called a {\it reduced decomposition} of
 $w\in W$ if $w=w(\ii)$
and the length $l(w)$ of $w$ is equal to $\ell$.
And let $R(w)$ be the set of all reduced decompositions of $w$. We denote by  $w_0\in W$ the element
of the maximal length in $W$ and refer to it as the {\it longest
element} of $W$ . For each $w\in W$, denote by
\begin{equation}
\label{eq:|w|}
|w|
\end{equation}
the smallest subset $J$ of $I$ such that $w$ belongs to the sub-group of
$W$ generated by $s_j$, $j\in J$; and refer to $|w|$ as to the {\it support}
of $w$. In other words, $|w|$ is the set of all $j\in I$ such that $s_j$
occurs in each reduced decomposition of $w$. By definition $|w_0|=I$.

We say that a parabolic subgroup $P$ is {\it standard} if
$P\supset B$. Clearly, for each standard parabolic subgroup $P$,
there exists a unique subset $J=J(P)$ of $I$ such that $P$ is
generated by $B$ and all $U^-_j$, $j\in J(P)$.  Denote by $L_P$
the Levi factor of $P$. It is clear that $L_P$ is generated by $T$
and all $U_j^-,U_j$, $j\in J$. Let $W_P:=\Norm_{L_P}(T)/T$ be the
Weyl group of $L_P$.

\begin{claim} For each standard parabolic subgroup $P$ of $G$ the Weyl group,
$W_P$ is a Coxeter subgroup of $W$ generated by all $s_j$, $j\in J(P)$.
\end{claim}

Denote by $w_0^P\in W_P$ the longest element of $W_P$ and define
the {\it the parabolic element} $w_P\in W$ by
\begin{equation}
\label{eq:parabolic element}
w_P=w_0^P\cdot w_0
\end{equation}

By definition,  $w_P$ is the minimal representative of $w_0$ in
the set of right cosets $W_P\backslash W$. Also $|w_0^P|=J(P)$.

For $i\in I,$ define $\overline {s_i}\in G$ by
\begin{equation}
\label{eq:split si}
\overline s_i=x_i(-1)y_i(1)x_i(-1) =\phi_i
\begin{pmatrix}
0 & -1 \\
1 & 0
\end{pmatrix} \ .
\end{equation}

Each $\overline s_i$ belongs to $\Norm_G(T)$ and is a
representative of $s_i\in W$. It is  well-known (\cite {bou}) that
the elements $\overline {s_i}$, $i\in I$, satisfy the braid
relations. Therefore, we can associate to each $w\in W$ its {\it
standard representative} $\overline  w\in \Norm_G(T)$ in such a
way that for any $(i_1,\ldots,i_\ell)\in R(w),$ we have
\begin{equation}
\label{eq:split w}
\overline w=\overline {s_{i_1}}\cdots \overline {s_{i_\ell}} \ .
\end{equation}

\begin{claim}
\label{cl:w bar iota} For each $w\in W$, one has $\iota(\overline
w)=\overline{w^{-1}}$.

\end{claim}

Define the monoid of dominant co-weights  $X^\star(T)^+$ and the monoid of
dominant weights $X_\star(T)^+$ by
$$X_\star(T)^+=\{\lambda\in X_\star(T):\left<\alpha_i,\lambda\right>\ge 0~\forall i\in I\} \ ,$$
$$X^\star(T)^+=\{\mu\in X^\star(T):\left<\mu,\alpha_i^\vee\right>\ge 0~\forall i\in I\}\ .$$

Following \cite[Section 4.3]{bz-invent}, for each $\mu\in X^\star
(T)^+$,  define the {\it principal minor} $\Delta_\mu$ to be the
regular function $G\to \AA^1$ uniquely determined by the property
$$\Delta_\mu(u_-tu_+)=\mu(t)$$
for all $u_-\in U^-$, $t\in T$, $u_+\in U^+$. For any extremal
weights $\gamma,\delta\in X^\star(T)$ of the form $\gamma=u\mu$,
$\delta=v\mu$ for some $u,v\in W$  and a dominant weight $\mu\in
X^\star (T)$, define the {\it generalized minor}
$\Delta_{\gamma,\delta}$ to be the regular function on $G$ given
by
\begin{equation}
\label{eq:generalized minor}
\Delta_{\gamma,\delta}(g)=\Delta_\mu(\overline u^{\,-1}g\overline v)
\end{equation}
for all $g\in G$. Clearly, ${\Delta_{\mu, \mu}(u)}=1$ and
\begin{equation}
\label{eq:chi Delta}
\chi_i(u)=\Delta_{\mu, s_i\mu}(u)
\end{equation}
for any $u\in U$, $i\in I$ and any $\mu\in X^\star(T)$ such that $\left<\mu,\alpha_i^\vee\right>=1$.
If $G$ is simply-connected, we will
use the formula
\begin{equation}
\label{eq:chi Delta omega}
\chi_i(u)=\Delta_{\omega_i, s_i\omega_i}(u)
\end{equation}
where $\omega_i$ is a {\it fundamental weight}, i.e, $\left<\omega_i,\alpha_j^\vee\right>=\delta_{i,j}$ for all $i,j\in I$.

\subsection{Basic facts on $U\times U$-varieties}
\label{subsect:basic facts on U-U varieties}

A $U\times U$-{\it variety} ${\bf X}$
is a pair
$(X,\alpha)$, where $X$ is an irreducible affine variety over $\QQ$ and
$\alpha:U\times X\times U\to X$ is a $U\times U$-action
on $X$, where the first $U$-action is {\it left} and second is {\it right},  such that each group $e\times U$ and $U\times e$ acts freely on $X$.
We will write the action as $\alpha(u,x,u')= uxu'$.

The $U\times U$-varieties form a category which morphisms are morphisms of underlying varieties commuting with the $U\times U$-action.

Define the {\it convolution product} $*$ of $U\times U$-varieties ${\bf X}=(X,\alpha)$ and ${\bf Y}=(Y,\alpha')$:
$${\bf X}*{\bf Y}:=(X * Y,\beta)\ ,$$ where
the variety $X * Y$  is the quotient of $X\times Y$  by the following  left action of $U$ on $X \times Y$:
$$u(x,y)=(xu^{-1},u y)\ .$$

For each $x\in X$, $y\in Y$, we denote by $x*y=\{(xu^{-1},uy)|u\in
U\}=U(x,y)$  the corresponding  point of $X*Y$. Then the action
$\beta:U\times X*Y\times U\to X$ is given by
$u(x*y)u'=(ux)*(yu')$,

Clearly, both actions $U\times e$ and $e\times U$ on $X*Y$ are free. So ${\bf X}*{\bf Y}$ is, indeed, a $U\times U$-variety.

\begin{claim}
\label{cl:strict monoidal Uvarieties}
The category of $U\times U$-varieties is naturally monoidal with respect to $*$ (and the unit object $U$).
\end{claim}

\begin{remark} Even though the category of $U$-varieties is not strict, we can ignore it for all practical purposes (see e.g., \cite{sch}).
The same applies to other monoidal categories  (unipotent
bicrystals, geometric crystals, Kashiwara crystals, etc.) which we
will deal with in this paper.
\end{remark}

Given a $U\times U$-variety ${\bf X}$, the right quotient $X/U$
and the left quotient $U\backslash X$ are well-defined
$U$-varieties. For each $x\in X$ define $U_{x\bullet}=\{u\in
U:ux\in xU\}$, the stabilizer in $U$ of the point $xU $ in $X/U$,
and $U_{\bullet x}=\{u\in U:xu\in Ux\}$ -- the  stabilizer in $U$
of the point $Ux$ in $U\backslash X$.

\begin{claim}
\label{cl:Uorbits}
Let ${\bf X}$ be a $U\times U$-variety, and $x\in X$ be a point. There is a unique group isomorphism $\varphi_x:U_{\bullet x}\to U_{x\bullet}$ such that
$xu=\varphi_x(u)x$ for all $u\in U_{\bullet x}$. In particular, the orbit $UxU$ is
isomorphic to the quotient of
$U\times U$ by the following action $\diamond$ of $U_{\bullet x}$: $$u_{\bullet x}\diamond(u,u')= (u\varphi_x(u_{\bullet x})^{-1},u_{\bullet x}u')\ .$$

\end{claim}

\begin{example}

\label{ex:Uw} For any  $w\in W$ and a representative $\tilde w\in
\Norm_G(T)$  of $w$, the {\it reduced Bruhat cell} \ $U\tilde wU$
in $G$ is a $U\times U$ variety. Denote $U(w):=U\cap \tilde
wU\tilde w^{-1}$ and $V(w)=U\cap \tilde  w U^-\tilde w^{\,-1}$
(clearly, $U(w)$ and $V(w)$ depend only on $w$). Then
$U_{\bullet\tilde w}=U(w^{-1}), U_{\tilde w\bullet}=U(w)$, the
isomorphism $\varphi_{\tilde w}:U(w^{-1})\to U(w)$ is given by
$\varphi(u)=\tilde w u\tilde w^{-1}$, and $U\tilde wU$ has two
unique factorizations
$$U\tilde wU=V(w)\tilde w U= U\tilde w V(w^{-1}) \ .$$
In particular, for $w=w_P$ as in (\ref{eq:parabolic element}), one
has $U_{\bullet \tilde {w_P}}=U\cap L_P$ and $V(w_P)=U_P$, where
$L_P$ is the Levi factor of $P$ and $U_P$ is the unipotent radical
of $P$,
 and one has two unique factorizations:
\begin{equation}
\label{eq:UPop}
U\tilde w_P U=U_P\tilde w_P U= U\tilde w U_{P^{op}} \ ,
\end{equation}
where $P^{op}=(w_0L_Pw_0)B$ is the standard parabolic subgroup opposite of $P$.

\end{example}

\begin{claim}
\label{cl:stabilizers of convolution product} For each point $x*y$
of $X*Y$, one has
$$U_{x*y\bullet}=\varphi_x(U_{\bullet x}\cap U_{y \bullet}),~U_{\bullet x*y}=\varphi_y(U_{\bullet x}\cap U_{y \bullet }) \ .$$
\end{claim}

\begin{claim}
\label{cl:xU,yU} For any $U$-varieties ${\bf X}$ and ${\bf Y},$
the morphism
$$ X*Y/U\to X/U\times Y/U$$
given by $x*yU\mapsto (xU,yU)$ is surjective; and the fiber over
any point $(xU,yU)$ is isomorphic to $U/U_{y\bullet}$ via
$xu*yU\mapsto u\cdot U_{y \bullet}$ for any $u\in U$.

\end{claim}

We say that a $U\times U$-variety ${\bf X}$ is {\it regular} if the action of $U\times U$ on $X$ is free or, equivalently,
each $U\times U$-orbit in $X$ is isomorphic to $U\times U$.

\begin{claim}
\label{cl:product regular} The convolution product of regular
$U\times U$-varieties is also regular.

\end{claim}

For any $U\times U$-variety ${\bf X},$ define  ${\bf
X}^{op}=(X,\alpha^{op})$, where $\alpha^{op}$ is the twisted
$U\times U$-action given by
$$(u,x,u')\mapsto \iota(u')\cdot x \cdot \iota(u) \ ,$$
where $\iota:G\to G$ is the positive inverse anti-automorphism (defined in \eqref{eq:iota}. Clearly, ${\bf X}^{op}$ is a well-defined
$U\times U$-variety. We will refer to ${\bf X}^{op}$ as the {\it opposite} $U\times U$-variety of ${\bf X}$.

\begin{claim}
\label{cl:op}
The correspondence ${\bf X}\to {\bf X}^{op}$ defines an involutive covariant functor from the category of $U\times U$-varieties into itself.
This functor reverses the convolution product, i.e., the $U\times U$-variety $({\bf X}*{\bf Y})^{op}$ is naturally isomorphic to
${\bf Y}^{op}*{\bf X}^{op}$.

\end{claim}

\subsection{Standard $U\times U$-orbits and their convolution products} We will consider here the $U\times U$-orbits in $G$ of the form
$X=U\overline wU$, where $w\in W$, and $\overline w$ is the
standard representative of $w$ in $\Norm_G(T)$ as defined in
\eqref{eq:split w}. We will refer to such a $U\times U$-orbit as
{\it standard}.

In the this section we investigate convolution products of standard $U\times U$-varieties. The following facts are well-known.

\begin{claim}

\label{cl:star}

{}\quad \begin{enumerate}\item[(a)] For any $w,w'\in W$, there is
a unique element $w''=w\star w'\in W$ such that $Bw''B$ is a dense
open subset of the product $BwB\cdot Bw'B=BwBw'B$. \item[(b)] The
map $(w,w')\to w\star w'\in W$ defines a monoid structure on $W$
such that and $s_i\star s_i=s_i$ for all $i\in I$ and
$$w\star w'=ww'$$
whenever $l(w'w')=l(w)+l(w')$. \item[(c)] For any $w,w',$ one has
$$(w\star w')^{-1}=w'^{-1}\star w^{-1} \ .$$
\item[(d)] For any $w,w'\in W$, there exists $w''\in W$ such that
$l(ww'')=l(w)+l(w'')$ and
$$w\star w'=ww''\ .$$\end{enumerate}
\end{claim}

The following result describes the convolution product of standard $U\times U$-orbits.

\begin{proposition}
\label{pr:generic} For each $w,w'\in W,$  all $U\times U$-orbits
in the convolution product $(U\overline wU)*(U\overline {w'}U)$
are isomorphic to $U\overline{w\star w'}U$. In particular,  if
$w\star w'=ww'$ (i.e., if  $l(ww')=l(w)+l(w')$), then
$$(U\overline wU)*(U\overline {w'}U)\cong U\overline {ww'}U \ .$$

\end{proposition}

\begin{proof} We will prove both parts by induction on $l(w)$. If $l(w)=0$, i.e., $w$ is
  the identity
element of $W$, we have nothing to prove.

Now let $w=s_i$ for some $i\in I$. We have to show  that each $U\times U$-orbit in $(U\overline {s_i}U)*(U\overline {w'}U)$ is isomorphic to
$U\overline {s_i\star w'}U$ for any $w'\in W$. In the proof we will use the obvious equality $U\overline {s_i}U=U\overline s_iU_i$

We first consider the case when $l(s_iw')>l(w')$, i.e., $s_iw'=s_i\star w'$ (we will  implicitly use here the obvious fact that
$(U\overline {s_i}U)(U\overline {s_iw'}U)=U\overline {s_iw'}U$ in $G$).
Since $U_i\overline {w'}U=\overline {w'} U$, we obtain
$$(U\overline {s_i}U)*(U\overline {w'}U)=(U\overline {s_i}U_i)*(\overline {w'}U)=(U\overline {s_i})*(U_i\overline {w'}U)=
(U\overline {s_i})*(\overline {w'}U) =U\overline {s_i}*\overline {w'}U\ .$$

Then, taking into account that $U_{\overline
{s_i}\bullet}=U(s_i)=U\cap \overline {s_i}U\overline
{s_i}^{\,-1}$, in the notation of Example \ref{ex:Uw} and Claim
\ref{cl:stabilizers of convolution product}, we have
$$U_{\overline {s_i}*\overline {w'} \bullet}=\overline s_i(U(s_i)\cap U(w'))\overline s_i^{\,-1}=U(s_i)\cap U(s_iw')=U(s_iw')=
U_{\overline {s_iw'}\bullet} \ .$$
Therefore, $(U\overline {s_i}U)*(U\overline {w'}U)=U\overline {s_i}*\overline {w'}U$ is isomorphic to
$U\overline {s_iw'}U=U\overline {s_i\star w'}U$ as a $U\times U$-orbit.

Second, consider the case when $l(s_iw')<l(w')$, that is, $w'=s_i\star w'$. Note that
$U\overline {s_i}U_i=U\overline {s_i}U_i=U\cdot U_i^-\cdot \overline {s_i}$.
Then we obtain the following decomposition into $U\times U$-orbits:
$$(U\overline {s_i}U)*(U\overline {w'}U)=\bigsqcup_{x\in U_i^-\cdot \overline {s_i}} \ Ux* \overline {w'}U \ .$$
Since $U(s_i)$ is the unipotent radical of that $i$-th minimal parabolic $P_i\supset B$ which has the Levi factor $L_i=T\cdot \phi_i(SL_2)$ and
$U_i^-\cdot \overline {s_i}\subset L_i$, we obtain for each $x\in U_i^-\cdot \overline {s_i}$:
$$U(s_i)\cdot x=x\cdot U(s_i) \ . $$
This implies that for any $x\in U_i^-\cdot \overline {s_i}$ we have $U_{x\bullet}=U(s_i)=U\cap \overline {s_i}U\overline {s_i}^{\,-1}$.

Therefore,  we obtain (again in the notation of Example \ref{ex:Uw} and Claim \ref{cl:stabilizers of convolution product}):
$$U_{x* \overline {w'} \bullet }=\overline s_i(U(s_i)\cap U(w'))\overline s_i^{\,-1}=U(s_i)\cap U(s_iw')=U(w')=U_{\overline {w'}\bullet} \ ,$$
that is, each orbit $Ux* \overline {w'}U$ of $(U\overline {s_i}U)*(U\overline {w'}U)$ is isomorphic to $U\overline {w'}U$.

This proves the assertion for $w=s_i$.

Furthermore, let $l(w)>1$. Let us consider the  convolution
product
$$Z=(U\overline wU)*(U\overline {w'}U) \ .$$
Since $l(w)>1$, there exists $i\in I$ such that $l(s_iw)<l(w)$.
Based on the above, we have an isomorphism of $U\times
U$-varieties $U\overline {w}U\cong (U\overline {s_i}U)*(U\overline
{s_iw}U)$, which, in turn, using Claim \ref{cl:strict monoidal
Uvarieties}, implies the isomorphism:
$$Z \cong ((U\overline {s_i}U)*(U\overline {s_iw}U))*(U\overline {w'}U)= (U\overline {s_i}U)*((U\overline {s_iw}U)*(U\overline {w'}U)) \ .$$

On the other hand, by the inductive hypothesis, each orbit ${\mathcal O}$ of the $U\times U$-variety $(U\overline {s_iw}U)*(U\overline {w'}U)$
is isomorphic to $U\overline {(s_iw)\star w'}U$. Therefore, each $U\times U$-orbit in $Z$ is isomorphic to an orbit in
$$(U\overline {s_i}U)*(U\overline {(s_iw)\star w'}U)\ ,$$
and, by the already proved assertion, each $U\times U$-orbit in the latter $U\times U$-variety is isomorphic to
$$U\overline {s_i\star ((s_iw)\star w')}U=U\overline {(s_i\star (s_iw))\star w'}U=U\overline {w\star w'}U \ .$$

This finishes the proof of the proposition.
\end{proof}

\begin{remark} Unlike for the product of Bruhat cells $BwB$ and $Bw'B$ in $G$
 (as in Claim \ref{cl:star}),
Proposition \ref{pr:generic} guarantees that the convolution product $(BwB)*(Bw'B)$ contains only isomorphic $U\times U$-orbits.

\end{remark}

\begin{corollary} For any sequence $\ii=(i_1,\ldots,i_\ell)\in I^\ell$, each orbit in the $\ell$-fold convolution product
$(U\overline {s_{i_1}}U)*\cdots *(U\overline {s_{i_\ell}}U)$ is isomorphic to $U\overline {s_{i_1}\star\cdots \star s_{i_\ell}}U$.
In particular, if $\ii$ is a reduced decomposition of $w\in W$, then
$$(U\overline {s_{i_1}}U)*\cdots *(U\overline {s_{i_\ell}}U)\cong U\overline wU \ .$$
\end{corollary}



%
%
%


\subsection{Linear functions on $U\times U$-varieties}
\label{subsect:linear functions on UxU-varieties}
Let ${\bf X}$ be a $U\times U$-variety and let $\chi:U\to \AA^1$ be a character.
We say that a  function $f:X\to \AA^1$ is
$(U\times U,\chi)$-{\it linear} if
\begin{equation}
f(u\cdot x\cdot u')=\chi(u)+f(x)+\chi(u')
\end{equation}
for any $x\in X, u,u'\in  U$.

Denote by $\chi^{st}:U\to \AA^1$ the {\it standard} regular character
\begin{equation}
\label{eq:chi standard}
\chi^{st}:=\sum_{i\in I} \chi_i \ .
\end{equation}

\begin{lemma}
\label{le:op on unipotent bicrystals} A $(U\times
U,\chi^{st})$-linear function on $U\times U$-variety ${\bf X}$ is
also \linebreak $(U\times U,\chi^{st})$-linear on the $U\times
U$-variety ${\bf X}^{op}$ (see Claim \ref{cl:op}).

\end{lemma}

\begin{proof} Indeed, under the twisted action of $U\times U$ on $X$, we obtain for $u,u'\in U$, $x\in X$:
$$f(\iota(u')\cdot x\cdot\iota(u))=\chi(\iota(u'))+f(x)+\chi(\iota(u))=\chi(u')+f(x)+\chi(u)$$
because $\chi_i(\iota(u))=\chi_i(u)$ for all $i\in I$ by by Claim \ref{cl:chi_i iota}.
The lemma is proved.
\end{proof}

\begin{claim}
\label{cl:criterion of linearity}
 Let ${\bf X}$ be a $U\times U$-variety, $x\in X$ be a point, and $\chi:U\to \AA^1$ be a non-zero character of $U$.
 Then the orbit $UxU$ admits  a $(U\times U,\chi)$-linear function if and only if
\begin{equation}
\label{eq:chi x}
\chi(u)=\chi(\varphi_x(u))
\end{equation}
for each $u\in U_{\bullet x}$, where $\varphi_x:U_{\bullet x}\to
U_{x\bullet}$ is as
in Claim \ref{cl:Uorbits}.

\end{claim}

Clearly, any $(U\times U,\chi)$-linear function on a single orbit $UxU$ is determined by its value at $x$.


\begin{corollary} Let $X=U\tilde w U$ be a $U\times U$-orbit in $G$,
where $\tilde w\in \Norm(T)$ is a representative of $w\in W$. Let
$\chi\ne 0$ be a character of $U$. Then $X$ admits a $(U\times
U,\chi)$-linear function $f:X\to \AA^1$ if and only if
\begin{equation}
\label{eq:chiw}
\chi(u)=\chi(\tilde w^{-1}u\tilde w)
\end{equation}
for any $u\in U(w)=\tilde wU\tilde w^{-1}\cap U$ and  $f$ is determined by $f(\tilde w)$.
\end{corollary}


We say that a character $\chi:U\to \AA^1$ is {\it regular} if  $\chi(U_i)\ne 0$ for each $i\in I$.

The following result gives a surprisingly simple classification of standard $U\times U$-orbits that admit a $(U\times U,\chi)$-linear function.

\begin{proposition}
\label{pr:bruhat linear} Let $\chi$ be a regular character of $U$
and let $\widetilde w\in \Norm_G(T)$ be a representative of $w\in
W$. Assume that the $U\times U$-orbit  $U\widetilde wU$ admits a
$(U\times U,\chi)$-linear function. Then  $w=w_P$ for some standard
parabolic  subgroup $P$ of $G$ (see \eqref{eq:parabolic element}).

\end{proposition}

\begin{proof} Let $R^+\subset X^\star(T)$ be the set of all positive roots of $G$. For each $\alpha\in R^+$, let $U_\alpha$ be the
corresponding $1$-dimensional subgroup of $U$. Note that $\tilde w
U_\alpha \tilde w^{-1}=U_{w\alpha}$ for any $w\in W$; hence
$U_\alpha\subset U(w)=U\cap \tilde wU\tilde w^{-1}$ if and only if
$\alpha\in R^+\cap w(R^+)$.

Therefore, if a regular character $\chi:U\to \AA^1$ and $w\in W$ satisfy \eqref{eq:chiw}, then for each $i\in I$ one has:
either $w^{-1}(\alpha_i)\notin R^+$ or $w^{-1}(\alpha_i)=\alpha_j$ for some $j\in
I$. Equivalently, if we denote $\sigma=w^{-1}w_0$, the latter condition
reads: for each $i\in I$ either $\sigma(\alpha_i)\in R^+$ or
$\sigma(\alpha_i)=-\alpha_j$ for some $j\in I$. Let us show that
the latter condition holds if and only if $\sigma$ is equal to the
longest element $w_0^P$ of some standard parabolic subgroup $P$ of
$G$. Indeed, denote by $\Pi=\{\alpha_i|i\in I\}$ the set of simple
roots of $G$ and define
$$\Pi_1:=\Pi\cap \sigma^{-1}(-\Pi),~\Pi_2:=\Pi\cap
\sigma(-\Pi)\ .$$
Obviously, $\sigma(\Pi_1)=-\Pi_2$,
and $\sigma(\Pi\setminus \Pi_1)\subset R^+$.

Let $P=P(\Pi_1)$ be the corresponding standard parabolic and $\tau:=\sigma w_0^P$.
All we have to prove is that $\tau=e$. In order to do so it suffices to
show that $\tau(\Pi)\subset R^+$.
By definition of $\tau$,
$$\tau(\Pi_1)=\sigma w_0^P(\Pi_1)=\sigma
(-\Pi_1)\subset R^+ \ .$$
So, it suffices to show that
$\tau(\Pi\setminus \Pi_1)\subset R^+$ as well.

Assume, by contradiction, that this is false and
there exists  $\alpha_j\in \Pi\setminus \Pi_1$ such that
$\tau(\alpha_j)\in -R^+$.
First of all, since $\alpha_j\notin  \Pi_1$, we obtain
$w_0^P(\alpha_j)=\alpha_j+\beta$
where $\beta$ is a non-negative combination of $\Pi_1$. This implies that
$$\tau(\alpha_j)=\sigma w_0^P(\alpha_j)=\sigma(\alpha_j+\beta)=\sigma(\alpha_j)+\sigma(\beta)=\sigma(\alpha_j)-\beta' \ ,$$
where $\beta'$ is a non-negative combination of $\Pi_2$.

Therefore, $\sigma(\alpha_j)=\tau(\alpha_j)+\beta'$. Since $\tau(\alpha_j)\in -R^+$ and
$\sigma(\alpha_j)\in R^+$, this immediately implies that both $-\tau(\alpha_j)$ and $\sigma(\alpha_j)$ are
 positive combinations of $\Pi_2$. In  turn, this
implies that $\alpha_j$ is a positive combination of
$\sigma^{-1}(\Pi_2)=-\Pi_1$, i.e.,
$\alpha_j$ is a negative combination of $\Pi_1$. This
contradiction proves the proposition.  \end{proof}

Now let us construct such a $(U\times U,\chi^{st})$-linear function $f_P$ (where $\chi^{st}$ is defined in \eqref{eq:chi standard} on a
$U\times U$-sub-variety of $Bw_PB$.

\begin{claim}
\label{cl:chi U(w)} Let $w\in W$ and let $i,j\in I$ be such that $w(\alpha_i)=\alpha_j$. Then
$$\overline w ^{\,-1}x_i(a)\overline w=x_j(a) $$
for all $a\in \GG_a$,  and, therefore,
$$\chi^{st}(\overline w^{\,-1}u\overline {w})=\chi^{st}(u)$$
for any $u\in U_i$. In particular, for each standard parabolic
$P$, one has
$$\chi^{st}(\overline {w_P}^{\,-1}u\overline {w_P})=\chi^{st}(u)$$
for any $u\in L_P\cap U$, where $L_P$ is the Levi factor of $P$.
\end{claim}


Let $P$ be a  standard parabolic subgroup of $G$ with the Levi factor $L$. Denote by $Z(L_P)\subset T$  the center of $L_P$.
Proposition \ref{pr:bruhat linear} and Claim \ref{cl:chi U(w)} imply the following corollary.

\begin{corollary} A $U\times U$-orbit $U\tilde wU$ admits a $\chi^{st}$-linear function, if and only if
$\tilde w\in Z(L_P)\cdot \overline {w_P}$ for some  standard
parabolic subgroup $P\subset G$.
\end{corollary}

Define a function $f_P:UZ(L_P)\overline {w_P} U\to \AA^1$ by
\begin{equation}
\label{eq:fLchi standard}
f_P(ut\overline {w_P}u')=\chi^{st}(u)+\chi^{st}(u')
\end{equation}
for $u,u'\in U$, $t\in Z(L_P)$.

In particular, taking  $P=B$ (and $L_P=T$, $Z(L_P)=T$), we obtain
a function  $f_B:UT\overline {w_0} U\to \AA^1$ given by
\begin{equation}
\label{eq:fBchi standard}
f_B(ut\overline {w_0}u')=\chi^{st}(u)+\chi^{st}(u')
\end{equation}
for $u,u'\in U$, $t\in T$. The latter function, unlike the former, can be extended to any character $\chi:U\to \AA^1$ via $f_{B,\chi}\to \AA^1$:
\begin{equation}
\label{eq:fBchi}
f_{B,\chi}(ut\overline {w_0}u')=\chi(u)+\chi(u')
\end{equation}

Using Claims \ref{cl:criterion of linearity} and \ref{cl:chi
U(w)}, we obtain the following  result.

\begin{claim} For any standard parabolic $P,$ the function $f_P$ is a $(U\times U,\chi^{st})$-linear function on $UZ(L_P)\overline {w_P} U$.

\end{claim}

\begin{remark} Since each regular character $\chi:U\to \AA^1$  belongs to the $T$-orbit of $\chi^{st}$ under the adjoint action, i.e.,
$\chi(u)=\chi^{st}(t_0ut_0^{-1})$ for $t_0\in T$, one can obtain the $(U\times U,\chi)$-linear function $f_{P,\chi}$ on
$t_0UZ(L_P)\overline {w_P}Ut_0^{-1}=UZ(L_P)t_0\overline {w_P}t_0^{-1}U$ by twisting \eqref{eq:fLchi standard} with this $t_0\in T$.
\end{remark}

Define a projection $\pi^+:B^-\cdot U\to U$ by
\begin{equation}
\label{eq:piplus}
\pi^+(bu)=u
\end{equation}
for $b\in B^-$, $u\in U$.

The following result provides a very useful formula for computing $f_P$.
\begin{lemma}
\label{le:formula fL} For any $g\in UZ(L_P)\overline {w_P} U$, one
has\begin{equation} \label{eq:formula fL}
f_P(g)=\chi^{st}(\pi^+(\overline {w_P}^{\,-1}g))+\sum_{i\in I
\setminus J(P)} \chi_i(\pi^+(\overline
{w_P^{-1}}^{\,-1}\iota(g)))\ ,
\end{equation}
where $\iota:G\to G$ is defined in \eqref{eq:iota} (and $J(P)$ is the set of all those $i\in I$ for which $U_i^- \in P$). In particular,
$$f_B(g)=\chi^{st}(\pi^+(\overline {w_0}^{\,-1}g))+ \chi^{st}(\pi^+(\overline {w_0}^{\,-1}\iota(g)))$$
for $g\in Bw_0B$.
\end{lemma}

\begin{proof} We proceed similarly to the proof of  \cite[Proposition 5.12]{bk}. Indeed,
\linebreak
 according to Example \ref{ex:Uw} for
$\tilde w=t\overline w_P$, $t\in Z(L)$, one has a unique factorization $g=u_Pt\overline
 {w_P} u$, where $t\in Z(L)$, $u_P\in U_P$, and $u\in U$.
Clearly, $\pi^+(\overline {w_P}^{\,-1}g)=u$. Furthermore, let us write a unique
 factorization $u=u_Lu'$, where $u_L\in U(w_P^{-1})=U\cap w_0Lw_0$
and $u'\in V(w_P^{-1})$. Let $u_L:= \overline {w_P}u'\overline {w_P}^{-1}$. Clearly,
 $u_P'\in U\cap L$, therefore $tu_p'=u_P't$ and
$g=u_Pu_Lt\overline {w_P} u'$. This implies that
$\iota(g)=\iota(u') \iota(t\overline {w_P})\iota(u_Pu_L)=\iota(u')
\overline {w_P^{-1}}^{\,-1}t^{-1}\iota(u_Pu_L)$. Therefore,
$\pi^+(\overline {w_P^{-1}}^{\,-1}\iota(g))=\iota(u_L)\iota(u_P)$.
Note that\linebreak $\chi_i(\iota(u_Pu_L))=\chi_i(u_Pu_L)$ for all
$i\in I$ by Claim \ref{cl:chi_i iota} and
$\chi_i(u_Pu_L)=\chi_i(u_P)$ for all $i\in I\setminus J(P)$.
Putting it together, the right hand side of
 \eqref{eq:formula fL} equals
$$\chi^{st}(u)+\sum_{i\in I \setminus J(P)} \chi_i(u_P)=f_P(g) \ .$$
The lemma is proved.
\end{proof}

Taking into   account that for a simply-connected group $G$ one
has
$$\chi_i(\pi^+(g))=\frac{\Delta_{\omega_i,s_i\omega_i}(g)}{\Delta_{\omega_i,\omega_i}(g)} \ ,$$
we obtain the following result.

\begin{corollary} For any $g\in UZ(L_P)\overline {w_P} U$, one has
\begin{equation}
\label{eq:formula fL minors}
f_P(g)=\sum_{i\in I} \frac{\Delta_{w_P\omega_i,s_i\omega_i}(g)}{\Delta_{w_P\omega_i,\omega_i}(g)}+\sum_{i\in I \setminus J(P^{op})}
\frac{\Delta_{w_0s_i\omega_i,\omega_i}(g)}{\Delta_{w_0\omega_i,\omega_i}(g)}\ ,
\end{equation}
where $P^{op}$ is  the standard parabolic opposite of $P$ (as in \eqref{eq:UPop}). In particular,
$$f_B(g)=\sum_{i\in I}\frac{\Delta_{w_0s_i\omega_i,\omega_i}(g)+\Delta_{w_0s_i\omega_i,\omega_i}(g)}{\Delta_{w_0\omega_i,\omega_i}(g)}$$
for $g\in Bw_0B$.
\end{corollary}

We finish the section with a technical result which will be used in Section \ref{sect:Unipotent bicrystals}.

For each $w\in W$, let $\pi^w$ and $ {}^w\pi$ be regular morphisms
$BwB\to T$ given by
\begin{equation}
\label{eq:pi0} \pi^w(ut\overline wu')=t, {}^w\pi(u\overline
wtu')=t \
\end{equation}
for $u,u'\in U$, $w\in W$, and $t\in T$.

Clearly, both  $\pi^w$ and ${}^w\pi$ are
 $U\times U$-invariant,  and
$$w({}^w\pi(g))=\pi^w(g)$$
for all $g\in BwB$, $w\in W$.

\begin{lemma} \label{le:piw iota} For any $w\in W$, we have
$$\pi^w(\iota(g)^{-1})=t_w\pi^w(g)$$ for any $g\in BwB$,
where $t_w=\overline{w^{-1}}^{\, -1}\overline w^{\,-1}\in T$.

\end{lemma}

\begin{proof}Indeed, let $g=ut\overline w u'$ for some $u,u'\in U$,
$t\in T$. By definition, $\pi^w(g)=t$. Then
$$\iota(g)^{-1}=u_1t\overline {w^{-1}}^{\,-1}u'_1 \ ,$$
where $u_1=\iota(u)^{-1}\in U,u_1=\iota(u)^{-1}\in U$. Therefore,
$\iota(g)^{-1}=u_1t\cdot t_w\overline {w}u'_1$ and $\pi^w(\iota(g)^{-1})= t\cdot t_w=t_w\cdot t$. \end{proof}

\section{Unipotent bicrystals and geometric crystals}

\label{sect:Unipotent bicrystals}

\subsection{Unipotent bicrystals}

\label{subsect:Basic definitions unipotent}

A {\it unipotent bicrystal}, or $U$-bicrystal is a pair $({\bf
X},{\bf p})$, where  ${\bf X}$ is a $U\times U$-variety, and ${\bf
p}:X\to G$ is a $U\times U$-equivariant morphism, where the  action
$U\times G\times U\to G$ is given by $(u,g,u')\mapsto ugu'$.

A morphism of $U$-bicrystals is any structure-preserving morphism of underlying $U\times U$-varieties.

For $U$-bicrystals $({\bf X},{\bf p})$, $({\bf Y},{\bf p}')$, let
us define the convolution product by
$$({\bf X},{\bf p})*({\bf Y},{\bf p}')=({\bf X}*{\bf Y},{\bf p}'') \ ,$$
where ${\bf X}*{\bf Y}$ is the convolution product of $U\times U$-varieties as in Section \ref{subsect:basic facts on U-U varieties} and
${\bf p}'':X*Y\to G$ is defined by
${\bf p}''(x*y)={\bf p}(x){\bf p}'(y)$ (clearly, the morphism ${\bf p}''$ is well-defined).

\begin{claim}
\label{cl:monoidal} The category of $U$-bicrystals with the product $*$ and structure-preserving morphisms as arrows is naturally monoidal
(where $U$ is the  unit object).
\end{claim}

For each unipotent bicrystal $({\bf X},{\bf p},)$, let us define
the sub-variety $X^-\subset X$ by
\begin{equation}
\label{eq:Xminus}
X^-={\bf p}^{-1}(B^-)
\end{equation}
(i.e., $X^-={\bf p}^{-1}({\bf p}(X)\cap B^-)$).

We refer to the variety $X^-$ as the {\it unipotent crystal} associated to $({\bf X},{\bf p})$. The variety $X^-$ is never empty because any
$U\times U$-orbit in $G$ has a non-trivial intersection with $B^-$.

\begin{remark}  The above definition of a unipotent crystal is equivalent to the original definition in \cite{bk}: a rational
$U$-action on $X$ is given below by \eqref{eq:new U action}
and ${\bf f}:X^-\to B^-$ is the restriction of ${\bf p}$ to $X^-$. And this ${\bf f}$ is $U$-equivariant.

\end{remark}

Next, we will list some obvious facts about unipotent crystals.

\begin{claim}
\label{cl:Xminus} For any unipotent bicrystal $({\bf X},{\bf p})$,
the restriction to $X^-$ of the quotient map $X\to X/U$ is an open
inclusion (and, therefore, a birational isomorphism)
$j_X:X^-\hookrightarrow X/U$. In particular, $X^-$ possesses a
rational $U$-action via
\begin{equation}
\label{eq:new U action}
u(x):=u\cdot x\cdot \pi^+({\bf p}(u\cdot x))^{-1}
\end{equation}
for $x\in X, u\in U$, where $\pi^+:B^-\cdot U\to U$ is the projection to the second factor (defined in \eqref{eq:piplus}).

\end{claim}

For any unipotent bicrystals  $({\bf X},{\bf p})$, $({\bf
Y},{\bf p}')$ and any  $x^-\in X^-$, $y^-\in Y^-$, we have
${\bf p}(x^-){\bf p}'(y^-)\in B^-$. Therefore, the correspondence
\begin{equation}
\label{eq:jXY}
(x,y)\mapsto j_{X,Y}(x^-,y^-):=x^-*y^-
\end{equation}
for $x^-\in X^-,y^-\in Y^-$ is a well-defined morphism
$j_{X,Y}:X^-\times Y^-\to (X*Y)^-$.

\begin{lemma}
\label{le:birational minus} For any unipotent bicrystals $({\bf
X},{\bf p})$ and $({\bf Y},{\bf p}')$, the morphism  $j_{X,Y}$ is
an open inclusion (hence a  birational isomorphism) $X^-\times
Y^-\osr (X*Y)^-$.
\end{lemma}

\begin{proof} First, we need the following obvious fact.

\begin{claim}
\label{cl:x-,y-} For any unipotent bicrystals $({\bf X},{\bf p})$
and $({\bf Y},{\bf p}')$, the restriction of the surjective
morphism from Claim \ref{cl:xU,yU} to $X^-*Y^-\subset
X^-*Y/U\subset X*Y/U$ defines a biregular isomorphism
$$ X^-*Y^-\osr X^-\times Y^-\ ,$$
where we identified $X^-$ (resp.~$Y^-$) with its image in $X/U$
(resp.~in $Y/U$).

\end{claim}

This guarantees that the morphism \eqref{eq:jXY} is injective. To show that the image is dense it suffices to use the irreducibility of
$(X*Y)^-$ and to count the dimension:
$$\dim (X*Y)^-=\dim X*Y-\dim U=\dim X+\dim Y-2\dim U=\dim X^-+\dim Y^- \ .$$

The lemma is proved.
\end{proof}

\begin{claim} For any unipotent bicrystal $({\bf X},{\bf p})$, there exists a unique element $w\in W$ such that the intersection
${\bf p}(X)\cap BwB$ is dense in ${\bf p}(X)$.

\end{claim}

We will refer to $w$ as the {\it type} of $({\bf X},{\bf p})$.

\begin{example}
\label{ex:standard bicrystal} For each $w\in W$, let
$X_w=U\overline w U$. Clearly, $X_w^-=B^-_w$ where we denote, following \cite{bk},
\begin{equation}
\label{eq:B-w}
B^-_w=B^-\cap U\overline w U \ .
\end{equation}
The pair $({\bf
X}_w,\id)$, where $\id$ stands for the natural inclusion
$\id:X_w\hookrightarrow G$ is a $U$-bicrystal.  Also the pair
$(BwB,\id)$ is a $U$-bicrystal such that $(BwB)^-=TB^-_w$. Each of these
$U$-bicrystals has type $w$.

\end{example}

\begin{claim} Let $({\bf X},{\bf p})$ and $({\bf Y},{\bf p'})$ be unipotent bicrystals of types $w$ and $w'$, respectively. Then the product
$({\bf X},{\bf p})*({\bf Y},{\bf p'})$ is of type $w\star w'$ (see Claim \ref{cl:star}).

\end{claim}

For any unipotent bicrystal $({\bf X},{\bf p})$ of type $w$,
define rational morphisms\linebreak  $hw_X, lw_X:X\to T$ by the
formula
\begin{equation}
\label{eq:hw_X}
hw_X(x)=\pi^w({\bf p}(x)), ~ lw_X(x)={}^w\pi({\bf p}(x)) \ .
\end{equation}
for $x\in X$, where $\pi^w$ and  ${}^w\pi$ are defined in \eqref{eq:pi0}.

We will refer to these morphisms as the {\it highest weight} of
$({\bf X},{\bf p})$ and the {\it lowest weight} of $({\bf X},{\bf
p})$, respectively.

\begin{claim} For any unipotent bicrystal $({\bf X},{\bf p})$, both the highest weight and the lowest weights $hw_X, lw_X:X\to T$ are
$U\times U$-invariant rational morphisms.

\end{claim}

In the notation of Claim \ref{cl:op}, for each $U$-bicrystal $({\bf X},{\bf p})$ denote $({\bf X},{\bf p})^{op}:=({\bf X}^{op},{\bf p}\circ \iota)$.
It follows immediately from Claim \ref{cl:op} that $({\bf X},{\bf p})^{op}$ is also a $U$-bicrystal. We will refer to it as the {\it opposite}
$U$-bicrystal of $({\bf X},{\bf p})$. Clearly, if $({\bf X},{\bf p})$ is of type $w$, then $({\bf X},{\bf p})^{op}$ is of type $w^{-1}$.

\begin{claim}
\label{cl:op bicrystal}
The correspondence $({\bf X},{\bf p})\mapsto ({\bf X},{\bf p})^{op}$ defines an involutive covariant functor from the category of
$U$-bicrystals into itself.

\end{claim}


\subsection{From unipotent bicrystals to  geometric crystals}

\label{subsect:unipotent to geometric}

We start with the recollection of some definitions and results of \cite{bk} on {\it geometric crystals} (in a slightly modified form).


\begin{definition}
\label{def:rational morphism} Given  varieties $X$ and $Y$ and a
rational morphism $f:X\to Y$,   denote by $\dom(f)\subset X$ the
maximal open subset of $X$ on which $f$ is defined; denote by
$f_{\reg}:\dom(f)\to Y$ the corresponding regular morphism. We
denote by  $\ran(f)\subset Y$ the closure of the constructible set
$f_{\reg}(\dom(f))$ in $Y$.
\end{definition}

It is easy to see that for any irreducible  algebraic varieties
$X,Y,Z$ and rational morphisms  $f:X\to Y,~g:Y\to Z$ such that
$\dom(g)$ intersects $\ran(f)$ non-trivially, the composition
$(f,g)\mapsto g\circ f$ is well-defined and is a  rational
morphism $X\to Z$.

For any algebraic group $H$, we call a rational action
$\alpha:H\times X\to X$ {\it unital} if $\dom(\alpha)\supset
\{e\}\times X$.


\begin{definition} A {\it geometric pre-crystal} is a $5$-tuple  ${\mathcal X}=(X,\gamma,
\varphi_i,\varepsilon_i,e_i^\cdot|$\linebreak $i\in I)$, where $X$
is an irreducible algebraic variety, $\gamma$ is rational morphism
$X\to T$, $\varphi_i,\varepsilon_i:X\to \AA^1$ are rational
functions, and each $e_i^\cdot:\GG_m\times X\to X$ is a unital
rational action of the multiplicative group $\GG_m$ (to be denoted
by $(c,x)\mapsto e_i^c(x)$) such that for each $i\in I$, one has
either:  $\varphi_i= \varepsilon_i=0$ and the action $e_i^\cdot$
is trivial, or: $\varphi_i\not = 0$, $\varepsilon_i\not = 0$, and
\begin{align*}\gamma(e_i^c(x))&=\alpha_i^\vee(c)\gamma(x),\varepsilon_i(x)
=\alpha_i(\gamma(x))\varphi_i(x), \\
\varepsilon_i(e_i^c(x))&=c\varepsilon_i(x),\varphi_i(e_i^c(x))=c^{-1}\varphi_i(x)\end{align*}
for $x\in X$, $c\in \GG_m$.

We will refer to the set of all $i\in I$ such that $\varphi_i\not
= 0$ as the {\it support} of ${\mathcal X}$ and denote it by
$\Supp {\mathcal X}$.

\end{definition}

\begin{example}
\label{ex:trivial geometric crystal} We consider  $T$ as a {\it trivial} geometric pre-crystal with $\varphi_i= \varepsilon_i= 0$ and the action
$e_i^\cdot$ is trivial for all $i\in I$. That is, the support of this geometric pre-crystal is the empty set $\emptyset$. Another example of a
trivial geometric pre-crystal is any sub-variety of $T$, in particular, a single point $e\in T$. Yet another example of a geometric pre-crystal on $T$ is a $5$-tuple ${\mathcal T}=(T,id_T,\varphi_i,\varepsilon_i,e_i^\cdot|i \in I)$,
where $e_i^c(t)=\alpha_i^\vee(c)\cdot t$ for all $c\in \GG_m$, $t\in T$, $i\in I$ and $\varepsilon_i,\varphi_i\in X^\star(T)$ are such that $\langle \varepsilon_i,\alpha_i^\vee\rangle=1$, $ \varphi_i=\varepsilon_i-\alpha_i^\vee$ for all $i\in I$ (e.g., $\varepsilon_i=\rho$, $\varphi_i=\rho-\alpha_i$ for all $i\in I$).
\end{example}

\begin{definition}
\label{def:morphism geometric pre-crystals} A {\it morphism} of
geometric pre-crystals ${\bf f}:\XX\to \XX'$ is a pair $(f,J)$, where $f$ is a  rational
morphism of underlying varieties $X\to X'$ and $J\subset \Supp \XX\cap \Supp \XX'$ such that $e_i^c \circ f=f\circ e_i^c$ for all $i\in \Supp \XX$, $c\in \GG_m$ and:
$$\varepsilon_j\circ f=\varepsilon'_j,~\varphi_j\circ f=\varphi'_j$$
for all $j\in J$
(we will refer to this $J$ as the {\it support} of $f$ and denote it by $\Supp {\bf f}$).

\end{definition}

Clearly, there is a category whose objects are geometric
pre-crystals and arrows are {\it dominant} morphisms of geometric
pre-crystals, where the composition of morphisms ${\bf f}=(f,J)$ and ${\bf f}'=(f',J')$ is defined by ${\bf f}'\circ {\bf f}:=(f'\circ f,J'\cap J)$
(i.e., $\Supp {\bf f'}\circ {\bf f}=\Supp {\bf f}'\cap \Supp {\bf f}$).

\medskip

Now we define the product of geometric pre-crystals.

\begin{definition}
\label{def:product geometric crystal}
Given geometric pre-crystals ${\mathcal X}=(X,\gamma,\varphi_i,\varepsilon_i,e_i^\cdot|i \in I)$ and
${\mathcal Y}=(Y,\gamma', \varphi'_i,\varepsilon'_i,e_i^\cdot|i \in I)$, define the $5$-tuple
$${\mathcal X}\times {\mathcal Y}=(X\times Y,\gamma'',e_i^\cdot, \varphi''_i,\varepsilon''_i|i \in I)\ ,$$
where

\noindent $\bullet$ the morphism $\gamma'':X\times Y\to T$ is given by
\begin{equation}
\label{eq:gamma of the product geometric}
\gamma''(x,y)=\gamma(x)\gamma'(y)
\end{equation}
for $x,y\in X\times Y$;

\noindent $\bullet$ the functions
$\varphi''_i,\varepsilon''_i:X\times Y\to \AA^1$ is given by
\begin{equation}
\label{eq:phi of the product geometric}
\varphi''_i(x,y)=
\varphi_i(x)+\frac{\varphi'_i(y)}{\alpha_i(\gamma(x))},~\varepsilon''_i(x,y)=\varepsilon'_i(y)+\varepsilon_i(x)\cdot \alpha_i(\gamma'(y))
\end{equation}
for $x,y\in X\times Y$;

\noindent $\bullet$ the rational morphisms $e_i^\cdot:\GG_m\times
X\times Y\to X\times Y$ are given by the formula (for $i\in \Supp
{\mathcal X}\cup \Supp {\mathcal Y}$)
\begin{equation}
\label{eq:tensor product of geometric crystals}
e_i^c(x,y)=(e_i^{c_1}(x),e_i^{c_2}(y)) \ ,
\end{equation}
for $x,y\in X\times Y$, where
\begin{equation}
\label{eq:tensor product of geometric crystals c_1 c_2}
c_1=\frac{c\varepsilon_i(x)+\varphi'_i(y)}
{\varepsilon_i(x)+\varphi'_i(y)},
~~c_2=\frac{\varepsilon_i(x)+\varphi'_i(y)}
{\varepsilon_i(x)+c^{-1}\varphi'_i(y)} \ .
\end{equation}


\end{definition}

\begin{claim} For any geometric pre-crystals ${\mathcal X}$ and ${\mathcal Y}$, the $5$-tuple
${\mathcal X}\times {\mathcal Y}$ is also a geometric pre-crystal
(with $\Supp {\mathcal X}\times {\mathcal Y}=\Supp {\mathcal
X}\cup \Supp {\mathcal Y}$). Moreover, the category of geometric
pre-crystals is  monoidal under the product $({\mathcal
X},{\mathcal Y})\mapsto {\mathcal X}\times {\mathcal Y}$ (where
the unit object is the single point $\{e\}$ as in Example
\ref{ex:trivial geometric crystal}).

\end{claim}

\begin{remark}  The formula
\eqref{eq:tensor product of geometric crystals} for the action of
$e_i^\cdot$ on $X\times Y$ provides a geometric analogue of the
tensor product of Kashiwara's crystals. See also
Remark~\ref{rem:kashiwara product inverted} below.

\end{remark}

\begin{example} For each geometric pre-crystal, $\XX$ the product $T\times \XX$ (where $T$ is considered the trivial geometric pre-crystal as in
Example \ref{ex:trivial geometric crystal}) is a geometric
pre-crystal  with $\Supp T\times \XX=\Supp \XX$. The projection to
the first factor is a morphism (with the empty support) of geometric pre-crystals $T\times
\XX\to T$.
\end{example}

Given a geometric pre-crystal ${\mathcal X}=(X,\gamma,\varphi_i,\varepsilon_i,e_i^\cdot|i \in I)$ denote by ${\mathcal X}^{op}$ the $5$-tuple
$(X,\gamma^{op},\varphi_i^{op},\varepsilon_i^{op},(e_i^\cdot)^{op}|i \in I)$, where $\gamma^{op}(x)=\gamma(x)^{-1}$, $\varphi_i^{op}=\varepsilon_i$,
$\varepsilon_i^{op}=\varphi_i$, and $(e_i^c)^{op}=e_i^{c^{-1}}$. Clearly, ${\mathcal X}^{op}$ is also a geometric pre-crystal which we will refer
to as the {\it opposite} of ${\mathcal X}$.

\begin{claim}
\label{cl:op geometric} The correspondence ${\mathcal X}\to
{\mathcal X}^{op}$ defines an involutive covariant functor from
the category of geometric pre-crystals into itself. This functor
reverses the product, i.e., $({\mathcal X}\times {\mathcal
Y})^{op}$ is naturally isomorphic to ${\mathcal Y}^{op}*{\mathcal
X}^{op}$. On the underlying varieties this isomorphism is the
permutation of factors $X\times Y\osr Y\times X$.

\end{claim}

For a geometric pre-crystal ${\mathcal X}$ and a sequence
$\ii=(i_1,\ldots,i_\ell)\in I^\ell$ we define a rational morphism
$e_{\ii}^\cdot :T\times X \to X$  by
\begin{equation}
\label{eq:general eii}
(t,x)\mapsto e_\ii^t(x)=e_{i_1}^{\alpha^{(1)}(t)}\circ
\cdots \circ e_{i_\ell}^{\alpha^{(\ell)}(t)}(x) \ ,
\end{equation}
where $\alpha^{(k)}=s_{i_\ell}s_{i_{\ell-1}}\cdots s_{i_{k+1}}(\alpha_{i_k})$, $k=1,\ldots,\ell$
are the associated roots.

\begin{definition}
A geometric pre-crystal
$(X,\gamma,\varphi_i,\varepsilon_i,e_i^\cdot|i \in I)$ is called a
{\it geometric crystal} if for any sequence
$\ii=(i_1,\ldots,i_\ell)\in I^\ell$ such that $s_{i_1}\cdots
s_{i_\ell}=1$  one has\begin{equation} \label{eq:verma2}
e_\ii=\id_X \ .
\end{equation}
\end{definition}

It is easy to see that the relations \eqref{eq:verma2} are
equivalent to the following relations between $e_i^\cdot,
e_j^\cdot$ for $i,j\in J$: \begin{align*} &e_i^{c_1}
e_j^{c_2}=e_j^{c_2} e_i^{c_1} \ \text{if}\
\langle\alpha_i,\alpha_j^\vee \rangle=0;\\
 &e_i^{c_1}e_j^{c_1c_2}e_i^{c_2}
=e_j^{c_2}e_i^{c_1c_2}e_j^{c_1} \  \text{if}\
 \langle\alpha_j,\alpha_i^\vee \rangle
= \langle\alpha_i,\alpha_j^\vee  \rangle=-1;\\
 &e_i^{c_1}e_j^{c_1^2c_2}e_i^{c_1c_2}e_j^{c_2}
=e_j^{c_2}e_i^{c_1c_2}e_j^{c_1^2c_2}e_i^{c_1}  \ \text{if}\
 \langle\alpha_j,\alpha_i^\vee \rangle=-2\langle\alpha_i,
 \alpha_j^\vee \rangle=-2;\\
  &e_i^{c_1} e_j^{c_1^3c_2} e_i^{c_1^2c_2}
e_j^{c_1^3c_2^2} e_i^{c_1c_2} e_j^{c_2} =e_j^{c_2} e_i^{c_1c_2}
e_j^{c_1^3c_2^2} e_i^{c_1^2c_2} e_j^{c_1^3c_2} e_i^{c_1}  
\text{if}\  \langle\alpha_j,\alpha_i^\vee \rangle=-3\langle\alpha_i,\alpha_j^\vee \rangle=-3. \end{align*}

\medskip

Since the above relations are invariant under taking the inverse of both hand sides,  the opposite ${\mathcal X}^{op}$ of a geometric crystal
${\mathcal X}$ is always a geometric crystal. Therefore, the correspondence ${\mathcal X}\mapsto {\mathcal X}^{op}$ is a covariant functor
from the category of geometric crystals into itself.

\begin{remark}
\label{rem:geometric crystals not monoidal}
Unlike geometric pre-crystals, the category of geometric crystals is not monoidal. This was a main reason for introducing unipotent crystals
in \cite{bk} and unipotent bicrystals in the present work.

\end{remark}

Next, we will construct geometric crystals out of unipotent bicrystals. For a unipotent bicrystal $({\bf X},{\bf p})$,  define a $5$-tuple
\begin{equation}
\label{eq:unipotent to geometric}
{\mathcal F}({\bf X},{\bf p})=(X^-,\gamma,\varphi_i,\varepsilon_i,e_i^\cdot|i\in I)
\end{equation}
as follows:

\noindent $\bullet$ The variety $X^-$ is given by  $X^-:={\bf p}^{-1}(B^-)$ as in \eqref{eq:Xminus}.

\noindent $\bullet$ A morphism $\gamma:X^-\to T$ is the composition of ${\bf p}:X^-\to B^-$ with the canonical projection $B^-\to B^-/U^-=T$.

\noindent $\bullet$ Regular functions
$\varphi_i,\varepsilon_i:X^-\to \AA^1$, $i\in I$ are as follows. Let $pr_i$ be the   natural projection $B^-\to B^-\cap \phi_i(SL_2)$
(where $\phi_i$ is a homomorphism $SL_2\to G$ defined in Section \ref{subsect:notation}).
Using the fact that $x\in X^-$ if and only if ${\bf p}(x)\in B^-$, we set:
\begin{equation}
\label{eq:phi on X}
\varphi_i(x):=\frac{b_{21}}{b_{11}},
~ \varepsilon_i(x):=\frac{b_{21}}{b_{22}}=\varphi_i(x)\alpha_i(x)
\end{equation}
for all $x\in X^-$, where $pr_i({\bf p}(x))=
\phi_i\begin{pmatrix}
b_{11} & 0 \\
b_{21} & b_{22}
\end{pmatrix}$.

\noindent $\bullet$ A rational morphism $e_i^\cdot:\GG_m \times X^- \to X$, $i\in I$ is given by ($x\in X^-$, $c\in \GG_m$):
\begin{equation}
\label{eq:simple multiplicative generator ei}
e_i^c(x)=x_i\left(\frac{c-1}{\varphi_i(x)}\right)\cdot x\cdot
x_i\left(\frac{c^{-1}-1}{\varepsilon_i(x)}\right)
\end{equation}
if $\varphi_i\not = 0$ and $e_i^c(x)=x$ if $\varphi_i= 0$.

\begin{remark} Even in the case when the homomorphism $\phi_i:SL_2\to G$ is not injective (i.e., when the fiber
$\phi_i^{-1}(b)\subset SL_2$ consists of two points which differ by a sign)  the functions
$\varphi_i$ and $\varepsilon_i$ are well-defined.

\end{remark}

\begin{remark}  It is easy to see that for any unipotent bicrystal $({\bf X},{\bf p}),$
we have: $\varphi_i\not = 0$ if and only if $i\in |w|$
(see \eqref{eq:|w|}), where $w$ is the type of $({\bf X},{\bf p})$.

\end{remark}

\begin{remark} For each $i\in I$
such that $\varphi_i\not = 0$, the action $e_i^\cdot$ is, indeed,
rational because it is undefined at the zero locus of $\varphi_i$
(or, which is the same, at the zero locus of $\varepsilon_i$).
\end{remark}

\newpage

\begin{proposition}
\label{pr:unipotent to geometric}

{}\quad

\begin{enumerate}
\item[(a)] For any unipotent bicrystal $({\bf X},{\bf p}),$ the
$5$-tuple ${\mathcal F}({\bf X},{\bf p})$ given by
\eqref{eq:unipotent to geometric} is a geometric crystal.
\item[(b)] For any unipotent bicrystals $({\bf X},{\bf p})$,
$({\bf Y},{\bf p}'),$ the formula \eqref{eq:jXY} defines an
isomorphism (of full support $I$) of geometric crystals
$${\mathcal F}({\bf X},{\bf p})\times {\mathcal F}({\bf Y},{\bf p}')\osr {\mathcal F}(({\bf X},{\bf p})*({\bf Y},{\bf p}')) \ .$$
\item[(c)] For any unipotent bicrystal $({\bf X},{\bf p})$, one
has
$${\mathcal F}(({\bf X},{\bf p})^{op})=({\mathcal F}({\bf X},{\bf p}))^{op}.$$
\end{enumerate}
\end{proposition}

\begin{proof} Prove (a). First, let us show that each $e_i^c$ takes $X^-$ to $X^-$. This immediately follows from the following fact.

\begin{lemma} For any $i\in I$ such that $\varphi_i\not = 0$, the rational action
\eqref{eq:new U action} is given on generators of $U$ by
\begin{equation}
\label{eq:unipotent crystal}
(x_i(a))(x) = x_i(a)\cdot x \cdot x_i\left(a'\right)
\end{equation}
for $x\in X^-$ and $a\in \GG_a$, where  $a'=-\frac{a}{\alpha_i(\gamma(x))(1+a\varphi_i(x))}$.

\end{lemma}

\begin{proof} Without loss of generality, it suffices to take $X=GL_2$, therefore, $X^-=\left\{\begin{pmatrix}
b_{11} & 0 \\
b_{21} & b_{22}
\end{pmatrix}: b_{11}b_{22}\ne 0\right \}$. We have
$$\begin{pmatrix}
1 & a \\
0 & 1
\end{pmatrix}\begin{pmatrix}
b_{11} & 0 \\
b_{21} & b_{22}
\end{pmatrix}\begin{pmatrix}
1 & a' \\
0 & 1
\end{pmatrix}=\begin{pmatrix}
b_{11}+ab_{21} & a'(b_{11}+ab_{21})+ab_{22} \\
b_{21} & b_{22}+a'b_{21}
\end{pmatrix} \ .
$$
Clearly, this product belongs to $X^-$ if and only if
$a'(b_{11}+ab_{21})+ab_{22}=0$. Since
$\alpha_i(\gamma(x))=\frac{b_{11}}{b_{22}}$ and
$\varphi_i(x)=\frac{b_{21}}{b_{11}}$, we obtain the desired
result.
\end{proof}

Furthermore, taking $a=\frac{c-1}{\varphi_i(x)}$ in \eqref{eq:unipotent crystal}, we obtain \eqref{eq:simple multiplicative generator ei}.
In turn,
each action $e_i^\cdot:X^-\to X^-$ given by \eqref{eq:simple multiplicative generator ei} coincides with that from \cite[Formula (3.7)]{bk}.
Since $X^-$ is a unipotent crystal, \cite[Theorem 3.8]{bk} implies that \eqref{eq:unipotent to geometric} defines a geometric crystal.
This proves part (a).

Part (b) now. For $x\in X^-$, $y\in Y^-$, we have
$\gamma(x*y)=pr({\bf p}(x){\bf p}(y))=pr({\bf p}(x))pr({\bf
p}(y))=\gamma(x)\gamma(y)$ because the natural projection
$pr:B^-\to T$ is a group homomorphism. This proves \eqref{eq:gamma
of the product geometric}. Prove \eqref{eq:phi of the product
geometric} now. Indeed, $pr_i({\bf p}(x))=b=\begin{pmatrix}
b_{11} & 0 \\
b_{21} & b_{22}
\end{pmatrix}$, $pr_i({\bf p}(y))=b'=\begin{pmatrix}
b'_{11} & 0 \\
b'_{21} & b'_{22}
\end{pmatrix}$
(where $pr_i:B^-\to B^-\cap \phi_i(SL_2)$ is defined above in \eqref{eq:unipotent to geometric}
) and
$$\varphi_i(x)=\frac{b_{21}}{b_{11}},\varepsilon_i(x)=\frac{b_{21}}{b_{22}},  \varphi_i(y)=\frac{b'_{21}}{b'_{11}},
\varepsilon_i(y)=\frac{b'_{21}}{b'_{22}},\alpha_i(\gamma(x))=\frac{b_{11}}{b_{22}},\alpha_i(\gamma(y))=\frac{b'_{11}}{b'_{22}} \ .$$

Therefore,
$$\varphi_i(x*y)=\frac{(bb')_{21}}{(bb')_{11}}=\frac{b_{21}b'_{11}+b_{22}b'_{21}}{b_{11}b'_{11}}=\varphi_i(x)+\frac{\varphi'_i(y)}{\alpha_i(\gamma(x))}$$
$$\varepsilon_i(x*y)=\frac{(bb')_{21}}{(bb')_{22}}=\frac{b_{21}b'_{11}+b_{22}b'_{21}}{b_{22}b'_{22}}
=\varepsilon_i(x)\cdot \alpha_i(\gamma'(y))+\varepsilon'_i(y) \ .$$
This proves \eqref{eq:phi of the product geometric}. It remains to prove \eqref{eq:tensor product of geometric crystals}.
Indeed, by definition \eqref{eq:simple multiplicative generator ei},
$$e_i^c(x*y)=x_i(a)\cdot x*y\cdot x_i(a')\ ,$$
where $a=\frac{c-1}{\varphi_i(x*y)}$, $a'=\frac{c^{-1}-1}{\varepsilon_i(x*y)}$. If we take $c_1$, $c_2$ as in
\eqref{eq:tensor product of geometric crystals c_1 c_2}, we easily see that $a=\frac{c_1-1}{\varphi_i(x)}$,
$b=\frac{c_2^{-1}-1}{\varepsilon_i(y)}$. This in conjunction with the identity
$\frac{c_1^{-1}-1}{\varepsilon_i(x)}=-\frac{c_2-1}{\varphi_i(y)}$ implies that
$$e_i^c(x*y)=x_i(a)\cdot x*y\cdot x_i(a')=x_i(a)\cdot x\cdot x_i(a'')*x_i(-a'')y\cdot x_i(a')=e_i^{c_1}(x)*e_i^{c_2}(y)\ ,$$
where $a''=\frac{c_1^{-1}-1}{\varepsilon_i(x)}=-\frac{c_2-1}{\varphi_i(y)}$.
This proves \eqref{eq:tensor product of geometric crystals}. Part (b) is proved.

Part (c) follows from the fact that $\varphi_i(\iota(b))=\varepsilon_i(b)$ and $\varepsilon_i(\iota(b))=\varphi_i(b)$ for all $b\in B^-$, $i\in I$.

The proposition is proved. \end{proof}

\begin{remark} As we argued above (Remark \ref{rem:geometric crystals not monoidal}),
the product of  geometric crystals is not always a  geometric crystal. However,
Proposition \ref{pr:unipotent to geometric}(b) implies that if the factors are geometric crystals
coming from unipotent bicrystals, then the product is always a geometric crystal.
\end{remark}

\begin{claim}
\label{cl:geometric highest weight morphism} For any unipotent
bicrystal $({\bf X},{\bf p})$ of type $w\in W$, one has:
\begin{enumerate}\item[(a)] The restriction  ${\bf p}|_{X^-}:X^-\to
TB^-_w$ defines a morphism of geometric crystals
\begin{equation}
\label{eq:homomorphism geometric crystals} {\bf f}_w:{\mathcal
F}({\bf X},{\bf p})\to {\mathcal F}(BwB,\id) \ ,
\end{equation}
where $(BwB,\id)$ is the $U$-bicrystal defined in Example \ref{ex:standard bicrystal}.  The support of ${\bf f}_w$ is
$\Supp {\mathcal F}({\bf X},{\bf p})=  \Supp {\mathcal
F}(BwB,\id)=|w|$ (see \eqref{eq:|w|}). \item[(b)] The restriction
of the highest weight morphism $hw_X$  (given by \eqref{eq:hw_X})
to $X^-$ defines a morphism of geometric crystals ${\mathcal
F}({\bf X},{\bf p})\to T$ (see Example~\ref{ex:trivial geometric
crystal}). The support of this morphism is $\emptyset$.
\end{enumerate}

\end{claim}

\subsection{Unipotent $\chi$-linear bicrystals}
\label{subsect:unipotent chi-linear bicrystals}
Let $\chi$ be a character of $U$.
A {\it unipotent $\chi$-linear bicrystal} or, $(U\times U,\chi)$-linear bicrystal is a triple
$({\bf X},{\bf p},f)$ where $({\bf X},{\bf p})$ is a unipotent bicrystal,
and  $f$ is a $(U\times U,\chi)$-linear function on $X$ (see Section \ref{subsect:linear functions on UxU-varieties}).

\medskip


A morphism of $(U\times U,\chi)$-linear bicrystals is any structure-preserving morphism of underlying $U$-bicrystals.

For any  $(U,\chi)$-linear bicrystals  $({\bf X},{\bf p},f)$ and
$({\bf Y},{\bf p}',f'),$ define the convolution product
\begin{equation}
\label{eq:bicrystal product}
({\bf X},{\bf p},f)*({\bf Y},{\bf p}',f'):=({\bf Z},{\bf p}'',f'') \ ,
\end{equation}
where $({\bf Z},{\bf p}'')=({\bf Z},{\bf p})*({\bf Y},{\bf p}')$ is the
product of $U$-bicrystals (see Section \ref{subsect:Basic definitions unipotent}) and
the function $f''$ on $Z=X*Y$ is defined by
$$f''(x*y)=f(x)+f'(y)\ .$$

For each $(U\times U,\chi)$-linear bicrystal $({\bf X},{\bf p},f)$
denote $({\bf X},{\bf p},f)^{op}:=(({\bf X},{\bf p})^{op},f)$.
Due to Lemma \ref{le:op on unipotent bicrystals}, this is also a
$(U\times U,\chi)$-linear bicrystal.

\begin{claim}
\label{cl:unipotent chi-linear bicrystals}

{}\quad

\begin{enumerate}\item[(a)] The category of unipotent $\chi$-linear
bicrystals is naturally monoidal with respect to the convolution
product $*$ (where the unit object is $(U,\id,\chi)$). \item[(b)]
The correspondence $({\bf X},{\bf p},f)\mapsto ({\bf X},{\bf
p},f)^{op}$ is an involutive covariant functor from the category
of unipotent $\chi$-linear bicrystals into itself. \end{enumerate}
\end{claim}

\medskip
We get a very surprising combinatorial corollary of the above arguments.

\begin{proposition}
In the notation of Section \ref{subsect:linear functions on UxU-varieties}, we have:
{}\quad

 \label{pr:Levi monoid} \begin{enumerate}
 \item[(a)]  For
any standard parabolic subgroups $P,P'$, there is a standard
parabolic subgroup $P''=P\star P'$ such that
$$w_P\star w_{P'}=w_{P''} \ .$$
\item[(b)] The operation $(P,P')\to P\star P'$  defines a monoid
structure on the set of standard parabolic subgroups. \item[(c)]
This monoid is commutative with the unit element $P=G$. \item[(d)]
$P\star P'\subset P\cap P'$ for any standard parabolic subgroups
$P$ and $P'$.
 \end{enumerate}
\end{proposition}

\begin{proof} Let  $f$ and $f'$ be  $(U\times U,\chi)$-linear
  functions
respectively on $Bw_PB$ and $Bw_{P'}B$ as in Proposition
\ref{pr:bruhat linear}. Then the triples $(U\overline {w_P}U,{\bf
  p},f)$
and $(U\overline {w_{P'}}U,{\bf p}',f')$, where ${\bf
p}:U\overline {w_P}U\hookrightarrow G$ and ${\bf p}':U\overline
{w_{P'}}U\hookrightarrow G$ are natural inclusions, are\linebreak
$(U\times U,\chi)$-linear bicrystals. Let $Z:=(U\overline {w_P}U)*(U\overline
{w_{P'}}U)$. Consider the convolution product $(U\overline
{w_P}U,{\bf p},f)*(U\overline {w_{P'}}U,{\bf p}',f')=(Z,{\bf
  p}'',f'')$. By definition,  $f''$ is
a $(U\times U,\chi^{st})$-linear function on $Z$.

By Proposition \ref{pr:generic} any $U\times U$ orbit in $Z$ is
isomorphic to $U\overline {w_P\star w_{P'}}U$. Hence the restriction of $f''$ to any such orbit becomes a $(U\times U,\chi)$-linear
function on this orbit, i.e., the standard $U\times U$-orbit
$U\overline {w_P\star w_{P'}}U$ admits a $(U\times U,\chi)$-linear
function. Therefore, according to Proposition
\ref{pr:bruhat linear},
the element $w_P\star w_{P'}$ must be of the form $w_{P''}$ for some standard parabolic subgroup $P''$. This proves (a)

It is clear that the operation $\star$ is associative, which proves (b). To prove (c), consider a map $\varphi:W\to W$ given by
$\varphi(w)=w_0w^{-1}w_0$. It is easy to see that $\varphi$ is
also an anti-automorphism of the monoid $(W,\star)$. On the other
hand, each  element of the form $w_P$ is fixed by $\varphi$
because
$$\varphi(w_P)=w_0(w_0^Pw_0)^{-1}w_0=w_0(w_0w_0^P)w_0=w_0^Pw_0=w_P \ .$$
Therefore,
$w_P\star w_{P'}=w_{P'}\star w_P$.
This proves the commutativity of the operation $P\star P'$. Part (c) is proved.

Prove (d) now. In order to prove the inclusion $P\star P'\subset
P\cap P'$, we use Claim~\ref{cl:star}. Recall that the {\it right}
Bruhat order on $W$ is defined by: $w\preceq w'$ if there exists
$w''$ such that $w''w=w'$ and $l(w')=l(w'')+l(w)$. One can easily
show that $w\preceq w'$ implies $w'w_0\preceq ww_0$. Take
$w=w_{P'}$, $w'=w_P\star w_{P'}=w_0^{P\star P'}w_0$. Then Claim
\ref{cl:star} guarantees that $w\preceq w'$, hence, by above,
$w'w_0\preceq ww_0$. Since $ww_0=w_0^{P'}$ and $w'w_0=w_0^{P\star
P'}$, we obtain the inequality
$$w_0^{P\star P'} \preceq w_0^{P'} \ ,$$
which is equivalent to ${P\star P'} \subseteq P'$. By   symmetry,
we also obtain ${P\star P'} \subseteq P$. Part~(d) is proved.

The proposition is proved. \end{proof}

\begin{remark} It would be interesting to find a combinatorial
or geometric description of the operation $(P,P')\mapsto P\star P'$. For instance, the product of $G$-orbits $G/P\times G/P'$
 contains a (unique) closed $G$-orbit isomorphic to $G/(P\cap P')$, which is a good approximation for  $G/(P\star P')$.
\end{remark}

\subsection{From unipotent $\chi$-linear bicrystals to decorated geometric  crystals}

\label{subsect:From unipotent chi-linear bicrystals to decorated geometric crystals}
We start the section with the useful notion of decorated geometric crystals.

\begin{definition}
\label{def:decorated geometric crystal}
A {\it decorated} geometric pre-crystal is a pair $({\mathcal X},f)$, where  ${\mathcal X}=(X,\gamma,\varphi_i,\varepsilon_i,e_i^\cdot|i\in I)$
is a geometric pre-crystal (as defined in  Section \ref{subsect:unipotent to geometric} above) with the support $I$, and $f$ is a
function on $X$ such that
\begin{equation}
\label{eq:f of e(x)}
f(e_i^c(x))=f(x)+
\frac{c-1}{\varphi_i(x)}+
\frac{c^{-1}-1}{\varepsilon_i(x)}
\end{equation}
for $x\in X$, $c\in \GG_m$, $i\in I$.

We say that a  decorated geometric pre-crystal $({\mathcal X},f)$ is a decorated geometric crystal if ${\mathcal X}$ is a geometric crystal.

\end{definition}

\begin{definition}
\label{def:product decorated geometric pre-crystals} For decorated
geometric pre-crystals $({\mathcal X},f)$ and $({\mathcal Y},f')$,
we define the {\it product}
\begin{equation}
\label{eq:product decorated geometric pre-crystals}
({\mathcal X},f)\times ({\mathcal Y},f'):=({\mathcal X}\times {\mathcal Y},f*f') \ ,
\end{equation}
where $f*f':X\times Y\to \AA^1$ is given by $(f*f')(x,y)=f(x)+f'(y)$.
\end{definition}

\begin{lemma}
\label{le:product decorated geometric crystal} The product \eqref{eq:product decorated geometric pre-crystals} of decorated geometric
pre-crystals is a well-defined decorated geometric pre-crystal.

\end{lemma}

\begin{proof} It suffices to show that the function $f*f'$ from \eqref{eq:product decorated geometric pre-crystals} satisfies
\eqref{eq:f of e(x)}. Indeed, in the notation of Definition
\ref{def:product geometric crystal}, we have
\begin{align*}(f*f')(e^c(x,y))&=(f*f')(e^{c_1}(x),e^{c_2}(y))=f(e^{c_1}(x))+f'(e^{c_2}(y))
\ .\\ &=f(x)+ \frac{c_1-1}{\varphi_i(x)}+
\frac{c_1^{-1}-1}{\varepsilon_i(x)}+f'(y)+
\frac{c_2-1}{\varphi'_i(y)}+
\frac{c_2^{-1}-1}{\varepsilon'_i(y)}.\end{align*} Taking into
account that
$$\varphi''_i(x,y)=\frac{
\varepsilon_i(x)+\varphi'_i(y)}{\alpha_i(\gamma(x))},~\varepsilon''_i(x,y)=(\varepsilon_i(x)+\varphi'_i(y))\cdot \alpha_i(\gamma'(y))\ ,$$
and
$$
\frac{c_1-1}{\varphi_i(x)}=\frac{c-1}{\varphi''_i(x,y)}, \frac{c_2^{-1}-1}{\varepsilon'_i(y)}=\frac{c^{-1}-1}{\varepsilon''_i(x,y)} \ ,
$$
$$\frac{c_1^{-1}-1}{\varepsilon_i(x)}=\frac{1-c}
{c\varepsilon_i(x)+\varphi'_i(y)},
~~\frac{c_2-1}{\varphi'_i(y)}=\frac{1-c^{-1}}
{\varepsilon_i(x)+c^{-1}\varphi'_i(y)}=-\frac{c_1^{-1}-1}{\varepsilon_i(x)}\ ,
$$
we obtain
$$(f*f')(e_i^c(x,y))=(f*f')(x,y)+
\frac{c-1}{\varphi''_i(x,y)}+
\frac{c^{-1}-1}{\varepsilon''_i(x,y)}  \ .$$
The lemma is proved.
\end{proof}

Clearly, the product of decorated geometric pre-crystals is associative. Therefore, the category of decorated geometric pre-crystals is monoidal.

It is easy to see that for a decorated geometric (pre)crystal
$({\mathcal X},f)$, the pair  $({\mathcal X}^{op},f)$ is also
decorated geometric (pre)crystal so that we will denote the latter
one by $({\mathcal X},f)^{op}$. Clearly, the association
$({\mathcal X},f)\mapsto ({\mathcal X},f)^{op}$ is a
multiplication-reversing functor from the category of decorated
geometric pre-crystals to itself.

For each $(U\times U,\chi)$-linear bicrystal $({\bf X},{\bf
p},f),$ we abbreviate
\begin{equation}
\label{eq:from unipotent linear decorated geometric}
{\mathcal F}({\bf X},{\bf p},f):=({\mathcal F}({\bf X},{\bf p}),f|_{X^-}) \ .
\end{equation}

The following fact is a ``decorated'' version of Proposition
\ref{pr:unipotent to geometric}.

\newpage

\begin{proposition}
\label{pr:uchi to decorated geometric}

{}\quad

 \begin{enumerate} \item[(a)]  In the notation of Proposition \ref{pr:unipotent
to geometric},   ${\mathcal F}({\bf X},{\bf p},f)$ is a decorated
geometric crystal for any $(U\times U,\chi^{st})$-linear bicrystal
$({\bf X},{\bf p},f)$. \item[(b)] ${\mathcal F}(({\bf X},{\bf
p},f)*({\bf Y},{\bf p}',f'))={\mathcal F}({\bf X},{\bf p},f)\times
{\mathcal F}({\bf Y},{\bf p}',f')$ for any $(U\times
U,\chi^{st})$-linear bicrystals $({\bf X},{\bf p},f)*({\bf Y},{\bf
p}',f')$. \item[(c)] In the notation of Claim \ref{cl:unipotent
chi-linear bicrystals}, we have $({\mathcal F}(({\bf X},{\bf
p},f)^{op}))=({\mathcal F}({\bf X},{\bf p},f)^{op}$.
\end{enumerate}
\end{proposition}

\begin{proof} Prove (a). It suffices to show that $f|_{X^-}$ satisfies \eqref{eq:f of e(x)}.
Indeed, by \eqref{eq:simple multiplicative generator ei}, we have
for each $x\in X^-$, $c\in \GG_m$, $i\in I$,
$$f(e_i^c(x))=f\left(x_i\left(\frac{c-1}{\varphi_i(x)}\right)\cdot x\cdot
x_i\left(\frac{c^{-1}-1}{\varepsilon_i(x)}\right)\right)$$
$$ =\chi^{st}\left(x_i\left(\frac{c-1}{\varphi_i(x)}\right)\right)
+f(x)+\chi^{st}\left(x_i\left(\frac{c^{-1}-1}{\varepsilon_i(x)}\right)\right)=f(x)+
\frac{c-1}{\varphi_i(x)}+
\frac{c^{-1}-1}{\varepsilon_i(x)} \ .$$

Parts (b) and (c) easily follow. The proposition is proved. \end{proof}


%
%

\subsection{Parabolic $(U\times U,\chi)$-linear bicrystals and the central charge}  We say
that a
$(U\times U,\chi)$-linear bicrystal $({\bf X},{\bf p},f)$ is
{\it parabolic} if a dense $U\times U$-invariant subset of the image ${\bf p}(X)$ in $G$
admits a $(U\times U,\chi)$-linear function.

Proposition \ref{pr:bruhat linear} guarantees that the type of each parabolic $(U\times U,\chi)$-linear bicrystal is $w_P$
for some standard parabolic subgroup $P\subset G$ (i.e., $P\supset B$).

Denote $X_P:=UZ(L_P)\overline {w_P}U$ for each standard parabolic subgroup $P$, where $L_P$ is the Levi factor of $P$.

\begin{claim}
\label{cl:standard parabolic and op} For each standard parabolic
subgroup $P$ of $G$, one has: \begin{enumerate} \item[(a)] $({\bf
X}_P,\id,f_P)$ is a parabolic $(U\times U,\chi^{st})$-linear
bicrystal of type $w_P$, where $f_P$ is the $(U\times
U,\chi^{st})$-linear function on $X_P$ defined in \eqref{eq:fLchi
standard} and\linebreak $\id:X_P\hookrightarrow G$ is the natural
(i.e., identical) inclusion. \item[(b)] The morphism $\iota:G\to
G$ (see \eqref{eq:iota}) defines an isomorphism of\linebreak
$(U\times U,\chi^{st})$-linear bicrystals (see Claim
\ref{cl:unipotent chi-linear bicrystals}):
$$({\bf X}_P,\id,f_P)^{op}\osr({\bf X}_{P^{op}},\id,f_{P^{op}})\ ,$$
where $P^{op}$ is the only standard parabolic which Levi factor
$L_{P^{op}}$ is related to the Levi factor $L_P$ of $P$ by
$L_{P^{op}}=w_0L_Pw_0$. In particular, $({\bf X}_B,\id,f_B)$ is
self-opposite. \item[(c)] If $P=G$, i.e.,  $({\bf X},{\bf p},f)$ is
a parabolic
  $(U\times U,\chi^{st}$)-linear bicrystal
of type $w_0$, then the pre-image $X_0={\bf p}^{-1}(Bw_0B)$ is a dense $U\times U$-invariant subset of
$X$ and the $U\times U$-action on $X_0$ is free and, therefore,
restriction of ${\bf p}$ to each $U\times U$-orbit ${\mathcal O}\subset
X_0$
is an inclusion ${\mathcal O}\hookrightarrow Bw_0B$.
\end{enumerate}
\end{claim}


According to Proposition \ref{pr:bruhat linear}, for any parabolic
 $(U\times U,\chi^{st})$-linear bicrystal $({\bf X},{\bf p},f)$ of type $w_P$
  the intersection of $X_P$ with ${\bf p}(X)$ is dense in ${\bf p}(X)$.

For any parabolic $(U\times U,\chi^{st})$-linear bicrystal $({\bf
X},{\bf p},f)$ of type $w_P$, define a rational function
$\Delta_X$ on $X$ by
\begin{equation}
\label{eq:Delta}\Delta_X(x):=f(x)-f_P({\bf p}(x))
\end{equation}
for any $x\in X$.

We will call this function  the {\it central charge}
of $({\bf X},{\bf p},f)$.

\begin{lemma}
\label{le:central charge} For any parabolic $(U\times
U,\chi^{st})$-linear bicrystal $({\bf X},{\bf p},f)$, the central
charge $\Delta_X$ is a $U\times U$-invariant rational function on
$X$.

\end{lemma}

\begin{proof} By definition,
\begin{align*}\Delta_X(uxu')=f(uxu')-f_P({\bf p}(uxu'))
&=f(x)+\chi^{st}(u)+\chi^{st}(u')-f_P(u{\bf p}(x)u')\\
&=f(x)-f_P({\bf p}(x)) =\Delta_X(x)\end{align*} for any $x\in X$,
$u,u'\in U$. This proves the lemma. \end{proof}

Similarly to \cite[Formula (4.4)]{bz-invent}, define an automorphism $\sigma:G\to G$ by
\begin{equation}
\label{eq:sigma}\sigma(g)=\overline {w_0}\cdot \iota\left(g\right)^{-1}\cdot\overline {w_0}^{\,-1} \ ,
\end{equation}
where the positive inverse $\iota:G\to G$ is given by \eqref{eq:iota}. By definition, $\sigma$ is an involution such that
$\sigma(B)=B^-$. According to \cite[Formula (4.5)]{bz-invent}, this involution satisfies
\begin{equation}
\label{eq:Delta sigma}
\Delta_{w_1\lambda,w_2\lambda}(\sigma(g))=\Delta_{w_0w_1\lambda,w_0w_2\lambda}(g) \
\end{equation}
for any $g\in G$, $w_1,w_2\in W$, and $\lambda\in X^\star(T)^+$.

\medskip
Next, define  a morphism $\pi:(Bw_0B)*(Bw_0B)\to U$ by
\begin{equation}
\label{eq:pi}
\pi(g*g')=uu' \ ,
\end{equation}
where $g=u_0t \overline {w_0}u,~g'=u't' \overline {w_0}u'_0 \ ,u,u',u_0,u_0'\in U,
t,t'\in T$.

\begin{claim} The morphism  $\pi$ is $U\times U$-invariant, i.e.,
$$\pi(ug*g'u')=\pi(g*g')$$
for any $g\in Bw_0B$, $g'\in Bw_0B$,$u,u'\in U$.
\end{claim}

Define  a morphism ${\bf m}:U\times T\times T\to B^-$ by
\begin{equation}
\label{eq:m}{\bf m}(u,t,t')= t\cdot \sigma\left(u\cdot t'\right) \ ,
\end{equation}
for $u\in U$, $t,t'\in T$.


\begin{proposition}
\label{pr:central charge of product} Let $({\bf X},{\bf p},f)$ and
$({\bf Y},{\bf p}',f')$ be parabolic $(U\times
U,\chi^{st})$-linear bicrystals of type $w_0$ each. Then the
central charge of the convolution product $({\bf X},{\bf p},f)* ({\bf Y},{\bf
p}',f')$ is given by the formula
\begin{equation}
\label{eq:central charge of product}
\Delta_{X*Y}(x*y)=\Delta_X(x)+\Delta_{Y}(y)+\chi^{st}(\pi(g,g'))+f_B({\bf
m}(\pi(g*g'),t,t')),
\end{equation}
 where $g={\bf p}(x)$, $g'={\bf p}'(y)$, $t=\pi^{w_0}(g),t'=\pi^{w_0}(g')$
  (and $\pi^{w_0}$ was defined in \eqref{eq:pi0}).

\end{proposition}


\noindent{\bf Proof of Proposition~\ref{pr:central charge of
product}}.\ We need the following  fact.

\begin{lemma}
\label{le:iota} For any standard parabolic subgroup $P$ with the
Levi factor $L_P$, the automorphism $g\mapsto (\iota(g))^{-1}$
preserves $X_P=UZ(L_P)\overline {w_P}U$,  and for any  $g\in X_P$,
we have  $f_P(\iota(g)^{-1})=-f_P(g)$.

\end{lemma}

\renewcommand{\qedsymbol}{\openbox}
\begin{proof} Note that
$\iota(U)^{-1}=U$ and $\iota(Z(L_P))^{-1}=Z(L_P)$.
Therefore, in order to show that $\iota(X_P)^{-1}=X_P$ it suffices to
prove that $\iota\left(\,\overline {w_P}\,\right)^{-1}\in Z(L_P)\overline {w_P}$.
Using Claim \ref{cl:w bar iota} and $w_0^{-1}=w_0$, $(w_0^P)^{-1}=w_0^P$, we obtain
$$\iota\left(\,\overline {w_P}\,\right)^{-1}=\left(\overline {(w_P)^{-1}}\right )^{-1}=
\overline {w_0w_0^P}^{\,-1}=\left(\overline {w_0}\overline {w_0^P}^{\,-1}\right)^{\,-1}
=\overline {w_0^P}\overline {w_0}^{\,-1} \ .$$
Therefore,
\begin{equation}
\label{eq:iotaL}
\iota\left(\,\overline {w_P}\,\right)^{-1}=\overline {w_0^P}^{\,2}\left(\overline
{w_0^P}^{\,-1}\overline {w_0}\right)\overline {w_0}^{\,-2}=\overline {w_0^P}^{\,\,2}\,\,\overline
{w_P}\,\,\overline {w_0}^{\,-2}
 \ .
\end{equation}
Finally, note that $\overline{w_0}^{\, -2}\in Z(G)$, $\overline
{w_0^P}^{\,\,2}\in Z(L_P)$.
Therefore,
$$\iota(\overline {w_P})^{\,-1}\in Z(L_P)\overline {w_P}Z(G)=Z(L_P)\overline {w_P} \ .$$

This proves that $\iota(X_P)^{-1}=X_P$. Next, prove that
$f_P(\iota(g)^{-1})=-f_P(g)$ for any $g\in X_P$. Indeed, for any
$g\in Z(L_P)\overline {w_P}$ and $u,u'\in U$ we obtain using
Claim~\ref{cl:chi_i iota}:
$$f_P(\iota(ugu')^{-1})=f_P(\iota(u)^{-1}\iota(g)^{-1}\iota(u')^{-1})=$$
$$\chi^{st}(\iota(u)^{-1})+\chi^{st}(\iota(u')^{-1})=-\chi^{st}(u)-\chi^{st}(u')=-f_P(ugu') \ .$$
This finishes the proof of the lemma.
\end{proof}

Furthermore, in order to prove \eqref{eq:central charge of product}, it
suffices to deal with the standard parabolic crystals of type $w_0$:
$({\bf X},{\bf p},f)=({\bf Y},{\bf p}',f')=({\bf X}_B,\id_{X_B},f_B)$.
Let $g=u_0t\overline {w_0}u$, $g'=u't'\overline {w_0}$ with $u,u'\in U$, $t,t'\in T$.

Thus (using the notation $\chi=\chi^{st}$ and the fact that
$\Delta=\Delta_{X_B*X_B}$ is $U\times U$-invariant), we obtain
\begin{align*}\Delta(g*g')&=f_B(g)+f_B(g')-f_B(gg')=\chi(u)+\chi(u')-f_{T}
(t\overline {w_0}uu't'\overline {w_0})\\
&=\chi(uu')-f_B(t\overline {w_0}uu't'\overline {w_0})=
\chi(u'')-f_{T}(t\overline {w_0}u''t'\overline {w_0})  \
,\end{align*} where $u''=uu'=\pi(g*g')$.

Furthermore, using Lemma \ref{le:iota} with $L_P=T$ and the fact
that $\iota(\overline {w})^{\, -1}=\overline {w^{-1}}^{\,-1}$ for
any $w\in W$, we obtain
\begin{align*}\Delta(g*g')&=\chi(u'')+f_B(\iota(t\overline
{w_0}u''t'\overline {w_0})^{-1})= \chi(u'')+f_B( t\cdot\overline
{w_0}^{\, -1}  \iota(u''t')^{-1}\cdot \overline {w_0}^{\,-1})
\\ &=\chi(u'')+f_B( t\cdot\overline {w_0}^{\, -2}  \overline
{w_0}\cdot \iota(u''t')^{-1}\cdot \overline {w_0}^{\,-1}) \\ &=
 \chi(u'')+f_B( \overline {w_0}^{\, -2}{\bf m}(u'',t,t'))=\chi(u'')+f_B( {\bf m}(u'',t,t'))
 \end{align*}
because $\overline {w_0}^{\, -2}\in Z(G)\subset T$.

This proves Proposition \ref{pr:central charge of product}. \endproof

\begin{remark} The formula \eqref{eq:central charge of product}
implies that $\Delta_{X*Y}$ is a rational (but not regular) function on $X*Y$ even when $\Delta_X$ and $\Delta_Y$ are regular
because the domain of $f_B$ is $Bw_0B$ while the range of the morphism
$g*g'\to {\bf m}(\pi(g*g'),\pi^{w_0}(g),\pi^{w_0}(g'))$ is the entire
$B^-$
(rather than the intersection $B^-\cap Bw_0B$).

\end{remark}

\subsection{The $(U\times U,\chi^{st})$-linear bicrystal $Z_w$}
\label{subsect:bicrystal Z_w} Now we explicitly describe the
unipotent bicrystal structure on $X_B*X_B=(Bw_0B)*(Bw_0B)$.
According to Claim~\ref{cl:product regular},  $X_B*X_B$ is a
regular $U\times U$-variety.

For each $w\in W$, denote by $(X_B*X_B)_w$ the set of all $g*g'\in
X_B*X_B$ such that $gg'\in BwB$. Clearly, $(X_B*X_B)_w$ is a
locally closed subset in $X_B*X_B$ and $X_B*X_B=\bigsqcup_{w\in W}
(X_B*X_B)_w$ and it is a $X_B*X_B$ is a regular $U\times
U$-variety.



Following \cite{bk}, define
\begin{equation}
\label{eq:Uw}
U^w:=U\cap B^-wB^-
\end{equation}
for each $w\in W$ (note a change of the notation). It is easy to
see that $\dim~U^w=\dim ~U^{w_0ww_0}=l(w)$. For each $w\in W$, let
\begin{equation}
\label{eq:Zw}
Z_w:=(U^{w_0ww_0}\times T\times T)\times_{T} BwB \ ,
\end{equation}
where the fibered product is taken with respect to  $\pi^w:BwB\to T$
(given by \eqref{eq:pi0}) and the following morphism $\mu_w:U^{w_0ww_0}\times T\times T\to T$:
\begin{equation}
\label{eq:mu w}
\mu_w(u,t,t')=\overline w_0^{\,2}\cdot (\overline w\, \overline {w^{-1}})\cdot\pi^w({\bf m}(u,t,t')) \ ,
\end{equation}
(where the morphism ${\bf m}:U\times T\times T\to B^-$ is defined by \eqref{eq:m}). By definition, $Z_w$ is a  $U\times U$-variety,
where the $U\times U$
acts on the second factor.

Now define a morphism $F_w:(X_B*X_B)_w\to U\times T\times T\times
BwB$ by the formula
$$F_w(g*g')=(\pi(g*g'),\pi^{w_0}(g),\pi^{w_0}(g'),gg')$$
for all $g*g'\in (X_B*X_B)_w$.

Clearly, each $F_w$ commutes with the $U\times U$-action
on
$U\times T\times T\times BwB$ (which is the trivial extension of that on $BwB$).

\begin{proposition}
\label{pr:Fw} For each $w\in W$, the image of $F_w$ is equal to
$Z_w$ and the quotient morphism
\begin{equation}
\label{eq:quotient Fw}
U \backslash F_w/U:U \backslash(X_B*X_B)_w/U\to U \backslash Z_w/U
\end{equation}
 is an isomorphism.
In particular, $F_{w_0}$ is an isomorphism of free $U\times
U$-varieties $(X_B*X_B)_{w_0}\osr Z_{w_0}$.
\end{proposition}

\begin{proof} First, note that all $U\times U$-orbits in $X_B*X_B$ are free. And each of these orbits intersects
the subset $\overline w_0T*B\overline w_0\subset X_B*X_B$ at exactly one point.

Let us study the subset $(\overline w_0T*B\overline w_0)_w$ of
$(X_B*X_B)_w$.
By definition, this is the set of all $\overline w_0t*t'u\overline w_0$,
$u\in U, t,t'\in T$ such that $\overline w_0tt'u\overline w_0\in BwB$ or,
equivalently, $tt'u\in \overline w_0^{ \, -1}BwB\overline w_0^{ \,-1}=B^-w_0ww_0B^-$,
that is, $u\in U^{w_0ww_0}$. Therefore,
$$(\overline w_0T*B\overline w_0)_w=\overline w_0T*TU^{w_0ww_0}
\overline w_0=T\overline w_0*U^{w_0ww_0}T\overline w_0 \ .$$

This characterization of $U\times U$-orbits in $(X_B*X_B)_w$ implies that
each element of $(X_B*X_B)_w$ can be uniquely expressed as $g*g'$, where
$g=u_1t\overline w_0$, $g'=ut'\overline w_0u_2$ for $u_1,u_2\in U$, $t,t'\in T$, $u\in U^{w_0ww_0}$.

Furthermore, by definitions \eqref{eq:pi} and \eqref{eq:pi0}, we
have
$$\pi(g*g')=u,~\pi^{w_0}(g)=t, \pi^{w_0}(g')=t', gg'=u_1t\overline w_0ut'\overline w_0u_2  \ .$$
That is,
\begin{equation}
\label{eq:special Fw}
F_w(g*g')=(u,t,t',u_1t\overline w_0ut'\overline w_0u_2)\in U^{w_0ww_0}\times T\times T\times BwB
\end{equation}
for $g=u_1t\overline w_0$, $g'=ut'\overline w_0u_2$ for $u_1,u_2\in U$, $t,t'\in T$, $u\in U^{w_0ww_0}$.
Next, we will prove that, in fact, $F_w(g*g')\in Z_w$. First, note that
$$\pi^w(gg')=\pi^w(u_1\overline w_0tu\overline w_0t'u_2)=\pi^w(t\overline w_0ut'
\overline w_0)\ .$$

Furthermore, using Lemma \ref{le:piw iota}, we obtain for $u\in U^{w_0ww_0}$, $t,t'\in T$:
$$\pi^w({\bf m}(u,t,t'))=\pi^w(t\sigma(u\cdot t'))=\pi^w(\iota\left(t\overline{w_0}^{\,-1}ut'\overline{w_0}\right)^{-1})
=t_w\pi^w(t\overline{w_0}^{\,-1}ut'\overline{w_0})$$
$$=t_w\pi^w(t\overline{w_0}^{\,-2}\overline{w_0}ut'\overline{w_0})=t_w\overline{w_0}^{\,-2}\pi^w(t\overline{w_0}ut'\overline{w_0})
=t_w\overline{w_0}^{\,-2}\pi^w(gg') \ .$$
Here we used the fact that $\overline {w_0}^{\, -2}\in Z(G)$.
Putting it together, we obtain:
\begin{equation}
\label{eq:pi mu}
\pi^w(gg')=\overline{w_0}^{\,2}t_w^{-1}\pi^w({\bf m}(u,t,t'))=\mu_w(u,t,t') \ .
\end{equation}
This proves that the image of $F_w$ belongs to $Z_w$.

Let us prove that the morphism \eqref{eq:quotient Fw} is an isomorphism. Consider the diagram:
$$\begin{CD}
(X_B*X_B)_w@>F_w>>Z_w\\
@AAA@AAA\\
T\overline w_0U^{w_0ww_0}* T\overline w_0@>>>U^{w_0ww_0}\times T\times T
\end{CD}
$$
where  \begin{enumerate}\item[1.]  The left vertical arrow is the
natural inclusion $T\overline w_0U^{w_0ww_0}* T\overline w_0
\subset$\linebreak $ (X_B*X_B)_w$. \item[2.] The right vertical
arrow is the injective morphism $U^{w_0ww_0}\times T\times T\osr
Z_w$ defined by $(u,t,t')\mapsto (u,t,t', \mu_w(u,t,t')\cdot
\overline w)$. \item[3.] The bottom horizontal arrow the
isomorphism $T\overline w_0U^{w_0ww_0}* T\overline w_0\to
U^{w_0ww_0}\times T\times T$ defined by $t\overline
w_0u*t'\overline w_0\mapsto (u,t,t')$.
\end{enumerate}
\smallskip

It follows from \eqref{eq:special Fw} and \eqref{eq:pi mu} that this
diagram commutes up to the $U\times U$-action on $Z_w$. Therefore,
passing to the quotient, we obtain the following commutative diagram:
$$\begin{CD}
U \backslash(X_B*X_B)_w/U@>U \backslash F_w/U>>U \backslash Z_w/U\\
@AAA@AAA\\
T\overline w_0U^{w_0ww_0}* T\overline w_0@>>>U^{w_0ww_0}\times T\times T
\end{CD}
$$
Clearly, the vertical arrows in these diagram (as well as the
bottom horizontal one) are isomorphisms. This implies the top
horizontal arrow $U \backslash F_w/U$ is also an isomorphism.

The proposition is proved. \end{proof}

Clearly, $(X_B*X_B)_{w_0}$ is a dense $U\times U$-invariant subset
in $X_B*X_B$. According to Proposition \ref{pr:Fw},  $F_{w_0}$ is
an isomorphism of free $U\times U$-varieties $({\bf X}_B*{\bf
X}_B)_{w_0}\osr {\bf Z}_{w_0}$.

Note that as a $U\times U$-variety, $Z_{w_0}=(U^{w_0}\times T\times T)\times_{T} X_B$, where
the fiber product is taken with respect to
$\mu_{w_0}:U^{w_0}\times T\times T\to T$
and  $\pi^{w_0}:X_B\to T$
 , where $\mu_{w_0}$ is given by
\eqref{eq:mu w}, In our case it simplifies to
$$\mu_{w_0}(u,t,t')=\overline {w_0}^{\,\,4}\pi^{w_0}({\bf m}(u,t,t'))=\pi^{w_0}
({\bf m}(u,t,t')) \ .$$

Using this characterization, define a morphism
$${\bf p}_Z:Z_{w_0}=(U^{w_0}\times T\times T)\times_{T} X_B\to X_B$$
to be  the projection to the second fibered factor
and a function $f_Z:Z_{w_0}\to \AA^1$ by
$$f_Z((u,t,t'),g)=f_B({\bf m}(u,t,t'))+\chi(u)+f_B(g)$$
for $u\in U^{w_0}$, $t,t'\in T$, and  $g\in X_B=Bw_0B$.

\begin{claim}
\label{cl:Zwnot} The morphism ${\bf p}_Z$ commutes with the
$U\times U$-action and the function $f_Z$ is $(U\times U,\chi^{st})$-linear, that is,
 the triple
\begin{equation}
\label{eq:Zwnot}({\bf Z}_{w_0},{\bf p}_Z,f_Z)
\end{equation}
is a $(U\times U,\chi^{st})$-linear bicrystal.

\end{claim}


Note that $(X_B*X_B)_{w_0}$ is open and dense in $X_B*X_B$. Then
Proposition \ref{pr:central charge of product}, Proposition
\ref{pr:Fw}, and Claim \ref{cl:Zwnot}  imply the following
corollary.

\begin{corollary} The inverse of the isomorphism $F_{w_0}:(X_B*X_B)_{w_0}\osr Z_{w_0}$ from Proposition \ref{pr:Fw}
defines an open embedding of $(U\times U,\chi^{st})$-linear bicrystals
\begin{equation}
\label{eq:F} ({\bf Z}_{w_0},{\bf p}_Z,f_Z) \hookrightarrow ({\bf
X}_B,\id,f_B)*({\bf X}_B,\id,f_B)\ .
\end{equation}

\end{corollary}


\section{Positive geometric crystals and unipotent bicrystals}

\subsection{Toric charts and positive  varieties}
\label{subsect:toric charts and positive structures}

Let $S$ be a split algebraic  torus defined over $\QQ$. Denote by
$X_\star(S)=\Hom(\GG_m,S)$ the lattice of {\it co-characters} of
$S$ and by $X^\star(S)=\Hom(S,\GG_m)$ the lattice of {\it
characters} of $S$. These lattices are dual: the integer pairing
$\left<\mu,\lambda\right>$ of $\lambda\in X_\star(S)$ and $\mu\in
X^\star(S)$ is defined via
$$\mu(\lambda(c))=c^{\left<\mu,\lambda\right>}$$
for $c\in \GG_m$.

The coordinate ring of $S$ is the group algebra $\QQ[X^\star(S)]$.

\begin{definition}
\noindent For any split algebraic torus, a {\it positive} regular
function on $S$ is an element of $\QQ[X^\star(S)]$ of the form
$f=\sum_\mu c_\mu \cdot \mu$, where all $c_\mu$ are non-negative
integers. A rational function on $S$ is said to be {\it positive} if it
can be expressed as a ratio $\frac {f}{g}$ where both $f$ and $g$
are positive regular functions on $S$ and $g\ne 0$.

\end{definition}

\begin{definition}
\noindent For any split algebraic tori $S$ and $S'$, we say that a
rational morphism $f$ is {\it positive} if for any character
$\mu':S'\to \GG_m$ the composition $\mu'\circ f$ is a positive
rational function on $S$.

\end{definition}

\begin{remark}
\label{rem:positive map} Equivalently (see e.g., \cite{polya} or \cite[Section 6]{reznick-hilbert}),
a rational function $f$ on $(\GG_m)^\ell$
is positive if and only if the restriction of $f$ to the positive octant
$(\QQ_{>0})^l$ is a well-defined  function $(\QQ_{>0})^\ell\to \QQ_{>0}$. Therefore, a rational morphism
$f:(\GG_m)^\ell\to (\GG_m)^k$ is positive if and only if the restriction of $f$ to the positive octant
$(\QQ_{>0})^l$ is a well-defined map $(\QQ_{>0})^\ell\to (\QQ_{>0})^k$. Taking into account that each split
algebraic torus $S$ is isomorphic to $(\GG_m)^\ell$ for some $\ell$ and the positive octant
$S(\QQ_{>0})\cong (\QQ_{>0})^l$ in $S(\QQ)$ does not depend on the choice of the isomorphism
$S\cong (\GG_m)^\ell$, the above arguments guarantee that a rational morphism $f:S\to S'$ is
positive if and only if the restriction of $f$ to $S(\QQ_{>0})$ is a well-defined  map $S(\QQ_{>0})\to S'(\QQ_{>0})$.

\end{remark}

\begin{remark} Obviously,  for any positive rational functions $f,g:S\to \AA^1$, the functions $f+g$, $fg$, and $\frac{f}{g}$ are also positive.
\end{remark}

\begin{remark}
\label{rem:non-positive inverse} Not every positive birational isomorphism $S\to S$ has positive inverse.
For instance, $f:\GG_m\times \GG_m\to \GG_m\times \GG_m$ given by
$f(x,y)=(x,x+y)$ is a positive isomorphism, but
$f^{-1}(x,y)=(x,y-x)$ is not positive.

\end{remark}

\begin{definition}
A {\it toric chart} on a variety $X$ is a birational isomorphism
$\theta:S\osr X$, where $S$ is a split algebraic torus.

\end{definition}

For any sequence $\ii=(i_1,\ldots,i_\ell)\in I^\ell$, define a
morphism $\theta^-_\ii:(\GG_m)^\ell\to B^-$ by
\begin{equation}
\label{eq:thetaii}
\theta^-_\ii(c_1,\ldots,c_\ell):=x_{-i_1}(c_1)\cdot x_{-i_2}(c_2)\cdots x_{-i_\ell}(c_\ell)
\end{equation}
for any $c_1,\ldots,c_\ell\in \GG_m$,  where
$x_{-i}:\GG_m\to B^-$ is given by the formula
\begin{equation}
\label{eq:thetai}
x_{-i}(c):= \phi_i\begin{pmatrix}
c^{-1}& 0 \\
1 & c
\end{pmatrix}.
\end{equation}

In a similar manner, following \cite{lu1}, for any sequence $\ii=(i_1,\ldots,i_\ell)\in
I^\ell,$ we define a morphism $\theta_\ii^+:(\GG_m)^\ell\to U$ by
\begin{equation}
\label{eq:thetaii plus}
\theta_\ii^+(c_1,\ldots,c_\ell):=x_{i_1}(c_1)\cdot x_{i_2}(c_2)\cdots x_{i_\ell}(c_{\ell})
\end{equation}
for any $c_1,\ldots,c_\ell\in \GG_m$, where each
$x_i:\GG_m\to U$ is a generator of $U$ (defined in Section \ref{subsect:notation}).

\begin{claim} \cite[Proposition 4.5]{bz-invent}
\label{cl:open chart} For any reduced decomposition
$\ii=(i_1,\ldots,i_\ell)$  of an element $w\in W$, one has:
\begin{enumerate}\item[(a)] The morphism $\theta^-_\ii$ is an open
embedding (hence a toric chart)
$$(\GG_m)^\ell\hookrightarrow B^-_w=B^-\cap U\overline w U \ .$$
\item[(b)]  The morphism $\theta_{\ii}^+$ is an open embedding
(hence a toric chart)
$$(\GG_m)^{l(w)}\hookrightarrow U^{w}=U\cap B^-wB^- \ .$$
\end{enumerate}
\end{claim}

\begin{example}
\label{ex:GL3 charts} For $G=GL_3$ and $\ii=(1,2,1),$ we have
 \begin{align*}\theta_\ii^-(c_1,c_2,c_3)&=\begin{pmatrix}
c_1^{-1}&0&0\\
1&c_1&0\\
0&0&1\\
\end{pmatrix}
\begin{pmatrix}
1&0&0\\
0&c_2^{-1}&0\\
0&1&c_2\\
\end{pmatrix}
\begin{pmatrix}
c_3^{-1}&0&0\\
1&c_3&0\\
0&0&1\\
\end{pmatrix}\\ &=\begin{pmatrix}
\frac{1}{c_1c_3}&0 &  0\\
\frac{c_1}{c_2}+\frac{1}{c_3}& \frac{c_1c_3}{c_2}&0\\
1&c_3&c_2
\end{pmatrix},\end{align*}

$$\theta_\ii^+(c_1,c_2,c_3)
=\begin{pmatrix}
1&c_1&0\\
0&1&0\\
0&0&1\\
\end{pmatrix}
\begin{pmatrix}
1&0&0\\
0&1&c_2\\
0&0&1\\
\end{pmatrix}
\begin{pmatrix}
1&c_3&0\\
0&1&0\\
0&0&1\\
\end{pmatrix}
=\begin{pmatrix}
1&c_1+c_3&c_1c_2\\
0&1&c_2\\
0&0&1
\end{pmatrix}.$$

\end{example}


We say that two toric charts $\theta:S\osr X$ and $\theta':S'\osr
X$ are {\it positively equivalent} if both birational isomorphisms
${\theta'}^{-1}\circ \theta:S\to S'$ and  $\theta^{-1}\circ
\theta': S'\to S$ are positive (see also Definition
\ref{def:positive equivalence}).

\begin{claim} \cite{bz-invent, lu1}
\label{cl:positively equivalent charts U B} For any reduced
decompositions $\ii$ and $\ii'$ of $w\in W$, one has:
\begin{enumerate}\item[(a)]  The toric charts  $\theta^-_\ii$ and $\theta_{\ii'}^-$ on
$B^-_w$ are positively equivalent. \item[(b)] The toric charts
$\theta^+_\ii$ and $\theta_{\ii'}^+$ on $U^w$ are positively
equivalent.
\end{enumerate}
\end{claim}

\begin{definition}
\label{def:positive structure} A {\it positive structure}
$\Theta_X$ on a  variety $X$ is a positive equivalence class of
its toric charts (if $X$ is  a split algebraic torus $S$,   then
we will always denote by $\Theta_S$  that positive structure on
$S$ which contains the identity isomorphism $\id:S\to S$); we will
call a pair $(X,\Theta_X)$ a {\it positive variety}.

\end{definition}

Define $supp(X,\Theta)$ of the positive variety $(X,\Theta)$ by
\begin{equation}
\label{eq:support positive} \Supp(X,\Theta)=\bigcup_{\theta\in
\Theta} \dom(\theta^{-1})
\end{equation}
(see Definition \ref{def:rational morphism}). That is,
$\Supp(X,\Theta)$ is a dense open subset of $X$.

For each rational morphism $\theta:S\to X$, where $S$ is a split
algebraic torus, we denote by $\theta(\QQ_{>0})\subset X(\QQ)$ the
image of the restriction of $\theta$ to $\dom(\theta)(\QQ)\cap
S(\QQ_{>0})$ (see Remark \ref{rem:positive map}).

\begin{claim} Let $(X,\Theta)$ be a positive variety. Then for each $\theta,\theta'\in \Theta$,
one has $\theta(\QQ_{>0})=\theta'(\QQ_{>0})$.

\end{claim}

Based on this, for each positive variety $(X,\Theta)$ define the set of {\it positive points} $X(\QQ_{>0})\subset X(\QQ)$ by
$X(\QQ_{>0}):=\theta(\QQ_{>0})$ for any $\theta\in \Theta$.

\begin{definition}
\label{def:positive morphism} A {\it morphism of positive
varieties} $(X,\Theta_X)\to (Y,\Theta_Y)$ is any rational morphism
$f:X\to Y$ such that one has (in the notation of
Definition~\ref{def:rational morphism}):
\begin{enumerate}\item[(1)] The intersection $\ran(f)\cap
\Supp(Y,\Theta_Y)$ is non-empty (hence it  is open in $\ran(f)$).
\item[(2)] For any toric chart $\theta:S\osr X$  in $\Theta_X$ and
any toric chart $\theta':S'\osr Y$ in $\Theta_Y$ such that
$\dom({\theta'}^{-1})\cap \ran(f)\not= \emptyset$, the composition
${\theta'}^{-1}\circ f\circ \theta$ is a positive morphism $S\to
S'$.\end{enumerate}

\end{definition}

We will sometimes refer to a morphism  of positive varieties $(X,\Theta_X)\to (Y,\Theta_Y)$ as a
{\it $(\Theta_X,\Theta_Y)$-positive morphism} $X\to Y$.

\begin{claim}
\label{cl:positive map varieties} Let $(X,\Theta_X)$ and
$(Y,\Theta_Y)$ be positive varieties. Then a rational
morphism  $f:X\to Y$ is a morphism of positive varieties
$(X,\Theta_X)\to (Y,\Theta_Y)$, if and only if the restriction of
$f$ to $X(\QQ_{>0})$ is a well-defined map $X(\QQ_{>0})\to
Y(\QQ_{>0})$. In the latter case, $\dom(\theta_Y^{-1})\cap
\ran(f)\not= \emptyset$ for any $\theta_Y\in \Theta_Y$.

\end{claim}

\begin{claim}
\label{cl:category positive varieties} Positive varieties and their morphisms form a category ${\mathcal V}_+$.
This category is monoidal with respect to the operation
$$(X,\Theta_X)\times (Y,\Theta_Y)=(X\times Y,\Theta_{X\times Y})\ ,$$
where $\Theta_{X\times Y}$ is the positive structure on $X\times Y$ such that $\theta\times \theta'\in \Theta_{X\times Y}$
for any $\theta\in \Theta_X$ and $\theta'\in \Theta_Y$.

\end{claim}

Following \cite{bk}, denote by ${\mathcal T}_+$ the category whose
objects are split algebraic tori (defined over $\QQ$), and arrows
are positive rational morphisms.

One has a natural functor ${\mathcal T}_+\hookrightarrow {\mathcal V}_+$ given by $S\mapsto (S,\Theta_S)$,
where $\Theta_S$ is the natural positive structure on $S$.

\begin{lemma}
\label{le:natural equivalence+}
The natural functor ${\mathcal T}_+\hookrightarrow {\mathcal V}_+$ is  an equivalence of monoidal categories ${\mathcal T}_+$ and  ${\mathcal V}_+$.

\end{lemma}
\begin{proof}
Clearly,  ${\mathcal T}_+\hookrightarrow {\mathcal V}_+$ is both a
full and faithful functor. This functor preserves a monoidal
structure.

\begin{claim}
\label{cl:toric chart positive} Given a positive variety
$(X,\Theta)$, each toric chart $\theta:S\osr X$ from $\Theta$ is
an isomorphism of positive varieties $(S,\Theta_S)\osr (X,\Theta)$
(where $\Theta_S$ is the natural positive structure on the
algebraic torus $S$).

\end{claim}

Claim \ref{cl:toric chart positive} guarantees that each object $(X,\Theta)$ is isomorphic to $(S,\Theta_S)$. The lemma is proved.
\end{proof}

Define the monoidal category ${\mathcal V}_{++}$ whose objects are
triples $(X,\Theta,\theta)$, where $(X,\Theta)$ is a positive
variety and $\theta\in \Theta$ is a chosen toric chart $S\osr X$
and each morphism $f:(X,\Theta,\theta)\to (Y,\Theta',\theta')$ in
${\mathcal V}_{++}$ is simply a morphism $f:(X,\Theta)\to
(Y,\Theta')$ of positive varieties. By definition, one has a full
forgetful monoidal functor ${\mathcal V}_{++}\to {\mathcal V}_+$
via $(X,\Theta,\theta)\mapsto (X,\Theta)$. We will refer to each
object of ${\mathcal V}_{++}$ as a {\it decorated positive
variety}.

\begin{claim}
\label{cl:V++}

{}\quad

 \begin{enumerate}\item[(a)] The forgetful functor
${\mathcal V}_{++}\to {\mathcal V}_+$ is an equivalence of
monoidal categories. \item[(b)] A simultaneous choice of the toric
chart $\theta:S\to X$ for each object
 $(X,\Theta_X)$ defines a functor
${\mathcal G}^*:{\mathcal V}_+\to {\mathcal V}_{++}$ adjoint to ${\mathcal G}$,
 and each adjoint to ${\mathcal G}$ functor is of this form.
\end{enumerate}
\end{claim}

Let us also define a functor $\tau:{\mathcal V}_{++}\to {\mathcal T}_+$ as follows:
\begin{enumerate}\item[(i)] $\tau(X,\Theta_X,\theta_X)=S$, where $\theta_X:S\osr  X$;
 \item[(ii)]  for each morphism $f:(X,\Theta_X,\theta_X)\to (Y,\Theta_Y,\theta_Y)$
 in ${\mathcal V}_{++}$, we define
$\tau(f):S\to S'$ to be the positive morphism of the form
$$\tau(f)=\theta_Y^{-1}\circ f\circ \theta_X \ .$$
According to Claim \ref{cl:positive map varieties}, $\tau(f)$ is well-defined.
\end{enumerate}

\begin{claim}

\label{cl:adjoint functor to positive tori}

{}\quad

 \begin{enumerate}\item[(a)]The functor $\tau:{\mathcal V}_{++}\to {\mathcal T}_+$ is an
equivalence of monoidal categories. \item[(b)] All  functors
${\mathcal G}^*:{\mathcal V}_+\to {\mathcal V}_{++}$ from Claim
\ref{cl:V++} are isomorphic to each other. \item[(c)] For each
functor ${\mathcal G}^*:{\mathcal V}_+\to {\mathcal V}_{++}$, the
composition $\tau \circ {\mathcal G}^*$ is a functor adjoint to
the natural inclusion ${\mathcal T}_+\hookrightarrow {\mathcal V}_+$ (see
Lemma \ref{le:natural equivalence+}).
\end{enumerate}
\end{claim}

Denote by  $\Theta_{\AA^1}$ the positive structure on $\AA^1$
containing the natural
inclusion $\GG_m\hookrightarrow \AA^1=\GG_m\sqcup \{0\}$; therefore, $(\AA^1,\Theta_{\AA^1})^n=(\AA^n,\Theta_{\AA^n})$,
where $\Theta_{\AA^n}$ is the positive structure containing the natural inclusion $(\GG_m)^n\hookrightarrow \AA^n$.

\begin{lemma}
\label{le:support codimension 1} For each $n>0$, we have
$\Supp(\AA^n,\Theta_{\AA^n})=(\GG_m)^n$.

\end{lemma}

\renewcommand{\qedsymbol}{}
\begin{proof} We need the following result.

\begin{claim}
\label{cl:positive rational bijection}
Let $\theta:\AA^n\to \AA^n$ be a positive birational isomorphism such that its inverse $\theta^{-1}$ is also positive.
Let $\theta_\RR:\RR^n\to \RR^n$ be the rational map obtained by the specialization of $\theta$ to the real points. Then:
 \begin{enumerate}\item[(a)] The restriction of $\theta_\RR$ to the positive octant $(\RR_{>0})^n$ is a well-defined homeomorphism
$(\RR_{>0})^n\osr (\RR_{>0})^n$. \item[(b)] If a boundary point
$x\in (\RR_{\ge 0})^n\setminus (\RR_{>0})^n$ belongs to the domain
of $\theta_\RR$, then $\theta_\RR(x)\in (\RR_{\ge 0})^n\setminus
(\RR_{>0})^n$.
\end{enumerate}
\end{claim}

To prove the lemma, we show that the assumption that
$\Supp(\AA^n,\Theta_{\AA^n})\setminus (\GG_m)^n$ is non-empty
leads to a  contradiction. Indeed, this assumption implies that
there exists a positive birational isomorphism $\theta:\AA^n\to
(\GG_m)^n$ such that $\theta^{-1}$ is also positive and
$\dom(\theta)\cap (\AA^n\setminus (\GG_m)^n)\ne \emptyset$.
Passing to the real points, and taking into account that
$(\AA^n\setminus (\GG_m)^n)(\RR)=\RR^n\setminus (\RR_{\ne 0})^n$,
we see that the restriction of $\theta_\RR:\RR^n\to (\RR_{\ne
0})^n$ to $ \RR^n\setminus (\RR_{\ne 0})^n$ is a real rational map
$\RR^n\setminus (\RR_{\ne 0})^n\to (\RR_{\ne 0})^n$. Since
$(\RR_{\ge 0})^n\setminus (\RR_{>0})^n$ contains a non-empty open
subset of $\RR^n\setminus (\RR_{\ne 0})^n$, the restriction of
$\theta_\RR$ to $(\RR_{\ge 0})^n\setminus (\RR_{>0})^n$ is a
rational map $(\RR_{\ge 0})^n\setminus (\RR_{>0})^n\to (\RR_{\ne
0})^n$. But the latter fact contradicts   Claim \ref{cl:positive
rational bijection}(b). This contradiction proves the lemma.
\end{proof}

Given a positive variety $(X,\Theta)$, we simply say that a  non-zero rational or regular function  $f:X\to \AA^1$ is
$\Theta$-positive if $f$ is positive $(\Theta,\Theta_{\AA^1})$-positive.

Based on Claim \ref{cl:positively equivalent charts U B}, for any
$w\in W$ we denote by $\Theta^w$ (resp.~by $\Theta_w^-$) the
positive structure  on $U^w=U\cap B^-wB^-$ (resp.~on
$B^-_w=B^-\cap U\overline wB$) containing each toric chart
$\theta^+_\ii$ (resp.~$\theta^-_\ii$)  for $\ii\in R(w)$.
Therefore, $(B^-_w,\Theta^-_w)$ and $(U^w,\Theta^w)$ are positive
varieties. The following result complements Lemma \ref{le:support
codimension 1}.

\begin{claim} \cite[Lemma 2.13]{bfz} The support of  each  positive variety $(B^-_w,\Theta^-_w)$ and $(U^w,\Theta^w)$
has co-dimension at least $2$ in $B^-_w$ and $U^w$, respectively.

\end{claim}

\begin{remark}
\label{rem:positive support GL_3} It is easy to show (see Example
\ref{ex:GL3 charts}) that for $G=GL_3$, the\linebreak complement
$U^{w_0}\setminus \Supp(U^{w_0},\Theta^{w_0})$ is precisely
$\begin{pmatrix}
1&0& \GG_m\\
0&1&0\\
0&0&1\\
\end{pmatrix}$. In particular,\linebreak the morphism $f:\GG_m\to U^{w_0}$ given by
$c\mapsto \begin{pmatrix}
1&0& c\\
0&1&0\\
0&0&1\\
\end{pmatrix}$ satisfies  $\ran(f)\cap$\linebreak $ \Supp(U^{w_0},\Theta^{w_0})=\emptyset$.
 Therefore, $f$ is not a morphism of positive varieties\linebreak
(because it fails the  requirement (1) of Definition \ref{def:positive morphism}).

\end{remark}

\begin{claim}
\label{cl:product standard positive structures} Let $w,w'\in W$ be
such that $l(ww')=l(w)+l(w')$. Then the natural  factorizations
$U^{ww'}\osr U^{w}\times U^{w'}$, $B^-_{ww'}=B^-_{w}\times
B^-_{w'}$
 are, respectively, the
isomorphisms of positive varieties
$$(U^{ww'},\Theta^{ww'})\osr (U^{w},\Theta^{w})\times (U^{w'},\Theta^{w'}),~(B^-_{ww'},\Theta^-_{ww'})\osr (B^-_{w},\Theta^-_w)\times (B^-_{w'}\times \Theta^-_{w'}) \ .$$

\end{claim}

Note that for any sub-torus $T'\subset T$ and $w\in W$ we have
canonical isomorphisms $T'\cdot U^w\cong T'\times U^w$ and
$T'\cdot B^-_w\cong T'\times B^-_w$. We will denote the positive
structures on $T'\cdot U^w$ and $T'\cdot B^-_w$ compatible with
these isomorphisms by $\Theta_{T'}\cdot \Theta^w$ and
$\Theta_{T'}\cdot \Theta^-_w$, respectively. Finally, for each
standard parabolic subgroup $P$, we abbreviate
$$\Theta^P:=\Theta_{Z(L_P)}\cdot \Theta^{w_P},
\Theta^-_P:=\Theta_{Z(L_P)}\cdot \Theta^-_{w_P} \ ,$$ where $L_P$
is the Levi factor of $P$.

Note that for $P=B$, the intersection $B^-\cap X_B$ is open and dense in $B^-$. Therefore, we regard each $\Theta^-_B$
as a positive structure on the entire $B^-$.

\begin{remark} We expect that $\Supp(B^-,\Theta^-_B)=\Supp(B^-_{w_0},\Theta_{w_0}^-)$.
\end{remark}

Recall that $\pi^+$ and $\pi^-$ are the regular projections
$B^-\cdot U\to U$ and $B^-\cdot U\to B^-$, respectively (as
defined in \eqref{eq:piplus}). Since $B^-\,w B^-\cdot
w^{\,-1}=B^-\, w U^-\cdot
 w^{\,-1} \subset B^-\cdot U$, we
can define a  morphism
$\eta^w:U^w\to B^-$ by
\begin{equation}
\label{eq:etaw}
\eta^w(u)=\iota(\pi^-(u\overline {w^{\,-1}}))
\end{equation}
for $u\in U^w$, where $\iota:B^-\to B^-$ is the positive inverse defined in \eqref{eq:iota}).
The following fact was proved in  \cite[Section 4.1]{bk}.

\begin{claim}
\label{cl:commutative diagram bz3} The morphism $\eta^w$ is a
biregular isomorphism $U^w\osr B^-_w$. The inverse isomorphism
$\eta_w=(\eta^w)^{-1}:B^-_w\osr U^w$ is given by
\begin{equation}
\label{eq:inverse eta}
\eta_w(b)=\iota(\pi^+(\overline w^{\,-1} b)) \ .
\end{equation}

\end{claim}

\begin{claim} \cite[Theorem 4.7]{bz-invent}
\label{cl:positive eta}  For each $w\in W$, the morphism
$\eta^w:U^w\to B^-_w$ is an isomorphism of positive varieties
$(B^-_w,\Theta^-_w)\osr (U^w,\Theta^w)$.

\end{claim}

The following fact was proved (in a slightly different form) in
\cite[Proposition~4.7]{bk}.

\begin{claim}

\label{cl:positive B-}  The multiplication $B^-\times B^-\to B^-$ is a morphism of positive varieties
$(B^-,\Theta_B^-)\times (B^-,\Theta_B^-)\to (B^-,\Theta_B^-)$.

\end{claim}

\begin{remark} Generalizing Claim \ref{cl:positive B-} and a part of Claim
\ref{cl:product standard positive structures},
one can easily show (see Claim \ref{cl:star} above and \cite[Proposition 4.7]{bk})
 that for any $w,w'\in W$, the multiplication morphisms
$U^w\times U^{w'}\to U^{w\star w'}$ and $B^-_w\times B^-_{w'}\to
T\cdot B^-_{w\star w'}$ are, respectively, the morphisms of
positive varieties $(U^w,\Theta^{w})\times (U^{w'},\Theta^{w'})\to
(U^{w\star w'},\Theta^{w'\star w''})$ and
$(B^-_w,\Theta^{w})\times  (B^-_{w'},\Theta^-_{w'})\to (T\cdot
B^-_{w\star w'},\Theta_T\cdot \Theta^-_{w'\star w''})$.

\end{remark}

\subsection{Positive geometric crystals and positive unipotent bicrystals}
\label{subsect:Positive structures on geometric crystals and
unipotent bicrystals} Let ${\mathcal
X}=(X,\gamma,\varphi_i,\varepsilon_i,e_i^\cdot|i\in I)$ be a
geometric pre-crystal and let $\Theta$ be a positive structure on
the variety $X$.  We say that  $({\mathcal X},\Theta)$ is a {\it
positive geometric pre-crystal} if:

\noindent $\bullet$ $\gamma:X\to T$ is a morphism of positive varieties $(X,\Theta_X)\to (T,\Theta_T)$.

\noindent $\bullet$ The functions $\varphi_i, \varepsilon_i:S\to \AA^1$ are $\Theta$-positive.

\noindent $\bullet$ For each $i\in I$, the action
$e_i^\cdot:\GG_m\times X\to X$ is ($\Theta_{\GG_m}\times \Theta$,
$\Theta$)-positive.
\smallskip

Recall that decorated geometric (pre)crystals are defined in Section
\ref{subsect:From unipotent chi-linear bicrystals to decorated
geometric crystals}. We say that a triple $(\XX,f,\Theta)$ is a {\it
positive} decorated geometric pre-crystal if $(\XX,f)$ is a
decorated pre-crystal, $(\XX,\Theta)$ is a positive pre-crystal, and
$f:X\to \AA^1$ is $\Theta$-positive.

Based on Definitions  \ref{def:product geometric crystal}, \ref{def:product decorated geometric pre-crystals} and
Claim \ref{cl:category positive varieties}, we can define product of
positive geometric pre-crystals and the product of decorated
positive geometric pre-crystals. Both products are associative.

\begin{definition} A {\it positive  $U$-bicrystal} of type $w$ is a triple $({\bf X},{\bf p},\Theta)$, where $({\bf X},{\bf p})$ is a
$U$-bicrystal of type $w$ and $\Theta$ is a positive structure on the variety $X^-={\bf p}^{-1}(B^-)$ such that:

\noindent (P1) The rational morphism ${\bf p}|_{X^-}:X^-\to TB^-_w$ is
 $(\Theta,\Theta_T\cdot \Theta^-_w)$-positive.

\noindent (P2) The pair $({\mathcal F}({\bf X},{\bf p}),\Theta)$
(see \eqref{eq:unipotent to geometric}) is a positive geometric crystal.

\smallskip

If  $f$ is a $(U\times U,\chi^{st})$-linear function on $X$ satisfying:

\smallskip

\noindent (P3) The function $f|_{X^-}:X^-\to \AA^1$ is $\Theta$-positive,

\smallskip

\noindent then we say that the quadruple $({\bf X},{\bf p},f,\Theta)$ is a positive $(U\times U,\chi^{st})$-linear bicrystal.

\end{definition}

In the notation of \eqref {eq:from unipotent linear decorated
geometric}, for any positive unipotent bicrystal $({\bf X},{\bf
p},\Theta)$, we denote
$${\mathcal F}({\bf X},{\bf p},\Theta):=({\mathcal F}({\bf X},{\bf p}),\Theta);$$
and for any positive $(U\times U,\chi^{st})$-linear bicrystal
$({\bf X},{\bf p},\Theta)$ we denote
$${\mathcal F}({\bf X},{\bf p},f,\Theta):=({\mathcal F}({\bf X},{\bf p},f),\Theta) \ .$$

\begin{lemma}
\label{le:from positive unipotent to positive geometric}

{}\quad

\begin{enumerate} \item[(a)] ${\mathcal F}({\bf X},{\bf p},\Theta)$  is a positive
geometric crystal for any positive  unipotent bicrystal $({\bf
X},{\bf p},\Theta)$. \item[(b)] ${\mathcal F}({\bf X},{\bf
p},f,\Theta)$  is a positive decorated geometric crystal for any
positive\linebreak $(U\times U,\chi^{st})$-linear bicrystal $({\bf
X},{\bf p},f,\Theta)$.\end{enumerate}
\end{lemma}

\begin{proof} The condition (P1) implies both the $(\Theta,\Theta_T)$-positivity of the morphism $\gamma$
and the $\Theta$-positivity of the functions $\varphi_i,
\varepsilon_i$ involved in the geometric crystal ${\mathcal F}({\bf
X},{\bf p})$. This proves (a). Part (b) also follows. The lemma is
proved.
\end{proof}

\begin{claim}
\label{cl:positive hw}
Let  $({\bf X},{\bf p},\Theta)$  be a positive unipotent bicrystal. Then
both rational morphisms $hw_X,lw_X:X^-\to T$
(defined in \eqref{eq:hw_X}) are $(\Theta,\Theta_T)$-positive.

\end{claim}

Furthermore, we will investigate the behavior of positive bicrystals under convolution products.

\begin{definition}
\label{def:star positive structure} Let $({\bf X},{\bf p},\Theta)$ and $({\bf Y},{\bf p}',\Theta')$ be positive unipotent bicrystals.  Denote
\begin{equation}
\label{eq:product positive unipotent}
({\bf X},{\bf p},\Theta)*({\bf Y},{\bf p}',\Theta'):=(({\bf X},{\bf p})*({\bf Y},{\bf p}'),\Theta* \Theta')\ ,
\end{equation}
where  $\Theta*\Theta'$ is that positive structure on the variety
$(X*Y)^-$ which is obtained from $\Theta\times \Theta'$ via the
birational isomorphism  $j_{X,Y}:X^-\times Y^-\to (X*Y)^-$ defined
in  \eqref{eq:jXY} (see also Lemma \ref{le:birational minus}).

Let $({\bf X},{\bf p},f,\Theta)$ and $({\bf Y},{\bf p}',f',\Theta')$
be positive $(U\times U,\chi^{st})$-linear bicrystals. Denote
\begin{equation}
\label{eq:product positive unipotent linear}({\bf X},{\bf p},f,\Theta)*({\bf Y},{\bf p}',f',\Theta'):=
(({\bf X},{\bf p},f)*({\bf Y},{\bf p}',f'),\Theta* \Theta')\ .
\end{equation}
\end{definition}

We can consider the category whose objects are positive unipotent
(resp.~positive $(U\times U,\chi^{st})$-linear)) bicrystals and
arrows are structure-preserving morphisms of positive varieties
(see Claim \ref{cl:category positive varieties}).

\begin{claim}
\label{cl:positive unipotent monoidal} The formula
\eqref{eq:product positive unipotent} (resp.~\eqref{eq:product
positive unipotent linear}) defines a monoidal structure  on the
category of positive unipotent (resp.~positive $(U\times
U,\chi^{st})$-linear) bicrystals.

\end{claim}


Recall from Example  \ref{ex:standard bicrystal} that the pair
$({\bf X}_w,\id)$ is a unipotent bicrystal.

\begin{claim}
\label{cl:standard bicrystal} \cite[Section 4.4]{bk} For each
$w\in W$, the triple  $({\bf X}_w,\id,\Theta_w^-)$ is a positive
unipotent bicrystal.

\end{claim}

\begin{definition} Let $({\bf X},{\bf p},f,\Theta)$  be a  positive parabolic $(U\times U,\chi^{st})$-linear bicrystal of type $w_P$.
We say that $({\bf X},{\bf p},f,\Theta)$ is  {\it strongly positive}
if the rational function $\Delta_X|_{X^-}:X^-\to \AA^1$ is
$\Theta$-positive (where  $\Delta_X:X\to \AA^1$ is defined in
\eqref{eq:Delta}).

\end{definition}

Next, we describe an important class of strongly positive parabolic $(U\times U,\chi^{st})$-linear bicrystals.

\begin{lemma}
\label{le:positive thetai} For any standard parabolic subgroup
$P$ of $G$, the quadruple $({\bf X}_P,\id,f_P,\Theta^-_P)$ is a strongly
positive parabolic $(U\times U,\chi^{st})$-linear bicrystal.

\end{lemma}

\renewcommand{\qedsymbol}{\openbox}

\begin{proof}  One can easily show (based on results of \cite{bz-invent} and \cite{bk}) that the restriction of the function
$f_P:X_P\to \AA^1$ to $X_P^-$ is a $\Theta_P^-$-positive function
$X_P^-\to \AA^1$. Therefore,  $({\bf X}_P,\id,f_P,\Theta_P^-)$ is
a positive parabolic $(U\times U,\chi^{st})$-linear bicrystal. At
the same time, it is  strongly positive since $\Delta_{X_P}$ is
zero. This proves the lemma. \end{proof}

   Our  main result on strongly positive unipotent bicrystals is the following.

\begin{maintheorem}

\label{th:positive unipotent monoidal} For any  strongly positive
parabolic $(U\times U,\chi^{st})$-linear bicrystals $({\bf X},{\bf
p},f,\Theta)$ and $({\bf Y},{\bf p}',f',\Theta')$ of type $w_0$,
their convolution product  is also a strongly positive $(U\times
U,\chi^{st})$-linear bicrystal of type $w_0$.

\end{maintheorem}

\begin{proof} Denote $({\bf Z},{\bf p}'',f'',\Theta'')=({\bf X},{\bf p},f,\Theta)*({\bf Y},{\bf p}',f',\Theta')$.
That is, $Z=X*Y$, $Z^-=(X*Y)^-=({\bf p}'')^{-1}(B^-)$, and $\Theta''=\Theta*\Theta'$. All we
have to show is that the restriction of ${\bf p}''$ to $Z^-$ is  a
positive morphism $Z^-\to B^-$ and the restriction of $\Delta_Z$
to $Z^-$ is a positive rational function  $Z^-\to \AA^1$.

First, the morphism ${\bf p}''|_{Z^-}:Z^-\to B^-$ factors as a composition
of positive  morphisms:
$$Z^-\buildrel {j_{X,Y}^{-1}} \over \longrightarrow X^-\times Y^-\buildrel {{\bf p}\times {\bf p}'} \over \longrightarrow B^-\times B^-\to B^- \ .$$
Therefore, ${\bf p}''|_{Z^-}$ is a also positive.

Furthermore,  \eqref{eq:inverse eta} and Claim \ref{cl:positive
eta} imply that the morphism $B^-\cap Bw_0B\to U$ defined by
$b\mapsto \pi^+(\overline w_0^{\,-1}b)$ (where $\pi^+:B^-\cdot
U\to U$ is the projection to the second factor defined above in
\eqref{eq:piplus}) are $(\Theta_B^-,\Theta^{w_0})$-positive. It is
also easy to see that the morphism and $B^-\cap Bw_0B\to T$
defined by $b\mapsto \pi^{w_0}(b)$  is $(\Theta^-_B,T)$ positive.
The multiplication morphisms $B^-\times B^-\to B^-$, $U\times U\to
U$, and the involution $\sigma$ are also positive. Therefore, both
$\pi$ and ${\bf m}$ (defined in \eqref{eq:pi} and \eqref{eq:m})
are positive. Finally, the functions $\chi^{st}:U\to \AA^1$ and
the restriction of $f_B$ to $B^-\cap Bw_0B$ are positive
functions. These arguments and \eqref{eq:central charge of product} imply
that the restriction of $\Delta_Z$ to $Z^-$  is a $\Theta''$-positive function.

This proves Theorem \ref{th:positive unipotent monoidal}.
\end{proof}

\begin{remark} Theorem  \ref{th:positive unipotent monoidal} implies that  the category of all strongly positive  parabolic
$(U\times U,\chi^{st})$-linear bicrystals of type $w_0$ is closed
under the convolution product. However, this category is not
monoidal because it has no unit object.

\end{remark}

\begin{remark} We expect that the result of
Theorem \ref{th:positive unipotent monoidal} holds for any parabolic
bicrystals.
The  difficulty  is to prove that $\Delta_{X*Y}$ is a positive function.

\end{remark}

We finish the section with the construction of a positive structure $\Theta_Z$ on the unipotent bicrystal $({\bf Z}_{w_0},{\bf p}_Z,f_Z)$ given by
\eqref{eq:Zwnot}.

Note that, by definition \eqref{eq:Zw}
$$Z_{w_0}^-=(U^{w_0}\times T\times T)\times_{T} (Bw_0B)^-\cong U^{w_0}\times T\times T \times B^-_{w_0} \ .$$
Using this isomorphism, we define a positive structure $\Theta_Z$  on $Z_{w_0}^-$ in such a way that
$$(Z_{w_0}^-,\Theta_Z)\cong (U^{w_0},\Theta^{w_0})\times (T,\Theta_T)\times (T,\Theta_T)\times (B_{w_0}^-,\Theta^-_{w_0}) \ .$$


\begin{claim}
\label{cl:strongly positive Zwnot}

The quadruple $({\bf Z}_{w_0},{\bf p}_Z,f_Z,\Theta_Z)$ is a strongly  positive $(U\times U,\chi^{st})$-linear bicrystal (see \eqref{eq:Zwnot}).

\end{claim}

Denote  $\XX_B={\mathcal F}({\bf X}_B,\id,f_B)$ and ${\mathcal
Z}_{w_0}={\mathcal F}({\bf Z}_{w_0},{\bf p}_Z,f_Z)$. Then the
inverse of the open embedding \eqref{eq:F} defines a birational
isomorphism of decorated geometric crystals
 \begin{equation}
 \label{eq:birational F}
\XX_B\times \XX_B \osr {\mathcal Z}_{w_0} \ .
\end{equation}

Furthermore, Lemma \ref{le:positive thetai} and Theorem
\ref{th:positive unipotent monoidal} imply that
$(\XX_B,\Theta_B)\times (\XX_B,\Theta_B)$ is a positive decorated
geometric crystal. Also Claim \ref{cl:strongly positive Zwnot}
implies that  $({\mathcal Z}_{w_0},\Theta_Z)$ is a positive
geometric crystal.

\begin{claim}
\label{cl:geometric tensor product} The isomorphism
\eqref{eq:birational F} defines a morphism of positive decorated
geometric crystals
\begin{equation}
\label{eq:positive F} (\XX_B,\Theta_B)\times (\XX_B,\Theta_B) \osr
({\mathcal Z}_{w_0},\Theta_Z) \ .
\end{equation}

\end{claim}

We expect that the inverse of \eqref{eq:positive F} is also positive (see Conjecture \ref{con:geometric tensor product}).

\section{Semi-field of polytopes  and tropicalization}
\label{sect:tropicalization}

\subsection{Semi-field of polytopes} Let $E$ be a finite-dimensional vector space over $\RR$. Denote by ${\mathcal P}_E$ the set
of all convex polytopes in $E$.

The {\it Minkovski sum} $P+Q=\{p+q|p\in P,q\in Q\}$ defines the
structure of an abelian monoid on ${\mathcal P}_E$ (with the unit
element - the single point $0$).


For any  $P\in {\mathcal P}_E$, define the {\it support function}
$\chi_P:E^*\to \RR$ by
$$\chi_P(\xi)=\min_{p\in P} \xi(p)  $$
for any $\xi\in E^*$.
Note that
$$\chi_P(\xi)=\min_{v\in \Vert(P)} \xi(v)  $$
where $\Vert(P)$ is the set of vertices of $P$.

Denote by $\Fun(E^*,\RR)$ the set of all functions $\tilde
f:E^*\to \RR$. Clearly, $\Fun(E^*,\RR)$ is an abelian monoid under
the pointwise addition of functions $(\tilde f,\tilde g)\mapsto
\tilde f+\tilde g$.

\begin{claim}
\label{cl:piecewise-linear realization} The correspondence
$P\mapsto \chi_P$ is an injective homomorphism of  monoids
$\chi:{\mathcal P}_E\to Fun(E^*,\RR)$, i.e., for any $P,Q\in
{\mathcal P}_E$, one has: \begin{enumerate} \item[(a)]
$\chi_{P+Q}=\chi_P+\chi_Q$. \item[(b)] $\chi_P=\chi_Q$ if and only
if  $P=Q$.\end{enumerate}
\end{claim}

This immediately implies the following well-known result.

\begin{corollary} The monoid ${\mathcal P}_E$ has a cancellation property:
\begin{equation}
\label{eq:cancelation}
P+R=Q+R<=>P=Q \ .
\end{equation}
\end{corollary}

For any two functions $\tilde f,  \tilde g\in Fun(E^*,\RR)$, we
write $\tilde f\le \tilde g$ if $\tilde f(\xi)\le \tilde g(\xi)$
for all $\xi\in E^*$. Clearly, this defines a partial order on
$\Fun(E^*,\RR)$.

\begin{claim} We have for any convex polytopes $P,Q$:
$$P\supseteq Q <=>\chi_P\le \chi_Q \ ,$$
that is, the two partial orders on ${\mathcal P}_E$ are equal.
\end{claim}

For any convex polytopes $P,Q\in {\mathcal P}_E$, let $P\vee Q$ be
the convex hull of  $P\cup Q$.
Clearly,
$$ \chi_{P\vee Q}=\min(\chi_P,\chi_Q)$$
for any  convex polytopes $P,Q\in {\mathcal P}_E$.
This implies the following  well-known
fact (see e.g., \cite{Glasauer}).

\begin{corollary} The operation $\vee$ satisfies
$$(P\vee Q)+R=(P+R)\vee(Q+R)$$
 for any  $P,Q,R\in {\mathcal P}_E$ (i.e., ${\mathcal P}_E$ is a semi-ring with the
  ``multiplication'' $+$ and the
 ``addition" $\vee$).

\end{corollary}

Therefore, if we make the monoid $\Fun(E^*,\RR)$ into a semi-ring
with the operation of ``addition'' $(\tilde f,\tilde g)\mapsto
\min(\tilde f,\tilde g)$ and the ``multiplication'' $(\tilde
f,\tilde g)\mapsto \tilde f+\tilde g$, the correspondence
$P\mapsto \chi_P$ is an order-preserving homomorphism of monoids.

Let us define the {\it Grothendieck group} ${\mathcal P}_E^+$ of the
monoid ${\mathcal P}_E$ by the standard construction: ${\mathcal
P}_E^+$ is the quotient of the free abelian group generated by
elements of the form $[P]$, $P\in {\mathcal P}_E$ by the relations
of the form $[P+Q]-[P]-[Q]$, $P,Q\in {\mathcal P}_E$. The group
${\mathcal P}_E^+$ is universal in the sense that one has the
canonical homomorphism of monoids ${\bf j}:{\mathcal P}_E\to
{\mathcal P}_E^+$ ($P\mapsto [P]$) and ${\bf j}({\mathcal P}_E)$
generates the group ${\mathcal P}_E^+$. The homomorphism ${\bf j}$
is injective due to Claim \ref{cl:piecewise-linear realization}. The
elements $[P]-[Q]$ of the group are referred to as {\it virtual
polytopes}  (see e.g.,  \cite{McMullen},
\cite{Panina}, \cite{Khovanski}).

The above operations and homomorphisms extend to the group
${\mathcal P}_E^+$ as  follows.

\smallskip

\begin{enumerate} \item[1.] The operation $\vee$ uniquely extends to the group
${\mathcal P}_E^+$ via the ``quotient'' rule:
$$([P]-[Q])\vee ([P']-[Q'])=[(P+Q')\vee (P'+Q)]-[Q+Q']$$
and, therefore, turns ${\mathcal P}_E^+$ into a semi-field.

 \item[2.] The homomorphism $\chi:{\mathcal P}_E \to Fun(E^*,\RR)$
defines an injective semi-field homomorphism ${\mathcal P}_E^+\to
Fun(E^*,\RR)$ via
$$[P]-[Q]\mapsto \chi_{[P]-[Q]}=\chi_P-\chi_Q \ .$$

 \item[3.] The containment relation $\supset$  extends to a
relation $\prec$ on ${\mathcal P}_E^+$  via
$[P]-[Q]\preceq [P']-[Q']$ if and only if $P+Q'\supseteq P'+Q$.
This relation satisfies
\begin{equation}
\label{eq:partial order}
[P]-[Q]\preceq [P']-[Q'] <=> \chi_{[P]-[Q]}\le \chi_{[P']-[Q']}
\end{equation}
for any $[P]-[Q], [P']-[Q']\in {\mathcal P}_E^+$.

\end{enumerate}
\medskip

\subsection{Newton polytopes and tropicalization}
\label{subsect:Newton polytopes and tropicalization} Let $S$ be a
split algebraic  torus defined over $\QQ$. For a regular function
$f$ on $S$ of the form $f=\sum _{\mu\in X^\star(S)} c_\mu \cdot
\mu$ define the {\it Newton polytope} of $f$ in $E=\RR\otimes
X^\star(S)$ to be the convex hull of all those $\mu\in X^\star(S)$
such that $c_\mu\ne 0$.

\begin{claim}
For any non-zero regular functions $f,g:S\to \AA^1$, one has
$$N(fg)=N(f)+N(g),~ N(f+g)\subseteq N(f)\vee N(g),$$
i.e, the association $f\mapsto N(f)$ is a homomorphism of monoids
 $\QQ[X^\star(S)]^\times\to {\mathcal P}_E$.

\end{claim}

\begin{definition} For any non-zero rational function $h:S\to \AA^1$, define
the {\it virtual Newton polytope} $[N](h)$ by
$$[N](h)=[N(f)]-[N(g)] \ ,$$
where   $f,g:S\to \AA^1$ are non-zero regular functions such that $h=\frac{f}{g}$.

\end{definition}

Clearly, $[N](h)$ is well-defined.

\begin{claim}
\label{cl:Newton} The correspondence $h\mapsto [N](h)$ is a
homomorphism of abelian groups $[N]:Frac(S)^\times\to {\mathcal
P}_E^+$ (where $Frac(S)$ is the field of rational functions
on~$S$), i.e.,
$$[N](fg)=[N](f)+[N](g) $$
for any non-zero rational  functions $f,g:S\to \AA^1$. This
homomorphism also\linebreak satisfies
$$[N](f)\vee [N](g)\preceq [N](f+g)$$
for any non-zero rational  functions $f,g:S\to \AA^1$ such that $f+g\ne 0$.
\end{claim}



For each rational function $h:S\to \AA^1$, define a  function
$\Trop(h):X_\star(S)\to \ZZ$  by
$$\Trop(h):=\chi_{[N](h)}  $$
and refer to $\Trop(h)$ as the {\it
tropicalization} of the rational function~$h$.

In particular, for any non-zero
regular function $f:S\to \AA^1$, we have
$$\Trop(f)(\lambda):=\chi_{N(f)}(\lambda)=\min_{\mu\in N(f)} \left<\mu,\lambda\right> \ $$
(here we used the identification of the lattice $ X_\star(S)$ with $( X^\star(S))^*$
via $\lambda\mapsto \left<\cdot,\lambda\right>$ for each
co-character $\lambda\in X_\star(S)$).

\begin{lemma}
\label{le:tropical function} For any non-zero rational functions:
$f,g:S\to \AA^1$, we have: \begin{enumerate}\item[(a)]
$\Trop(fg)=\Trop(f)+\Trop(g)$,
$\Trop(\frac{f}{g})=\Trop(f)-\Trop(g)$. \item[(b)] If $f+g\ne 0$
then $\Trop(f+g)\ge \min(\Trop(f),\Trop(g))$;  the equality
\begin{equation}
\label{eq:trop of sum} \Trop(f+g)=\min(\Trop(f),\Trop(g))
\end{equation}
is achieved if and only if $[N](f+g)=[N](f)\vee [N](g)$.
\end{enumerate}

\end{lemma}

\begin{proof} Follows from \eqref{eq:partial order} and Claim \ref{cl:Newton}.
\end{proof}

Below we will list some sufficient conditions for \eqref{eq:trop of sum} to hold.

\begin{definition} A rational function $f:S\to \AA^1$ on a split algebraic torus $S$ is
called {\it positive} if it can be written as a ratio $f=\frac{f'}{f''}$, where $f'$
and $f''$ are linear combinations of characters with positive integer coefficients.

\end{definition}

\begin{corollary}
\label{cor:trop of sum}  Let $f_1,f_2,\ldots,f_k$ be non-zero
rational functions on $S$. Then: \begin{enumerate}\item[(a)]
Assume additionally that the functions $f_1,f_2,\ldots,f_k$ are
regular and the vertex sets $\Vert(N(f_1))$, \ldots,
$\Vert(N(f_k))$ of their Newton polytopes are pairwise disjoint,
i.e., $\Vert(N(f_i))\cap \Vert(N(f_j))=\emptyset$ for all $i\ne
j$. Then
$$N(f_1+f_2+\cdots+f_k)=N(f_1)\vee N(f_2)\vee \cdots\vee N(f_k) $$
and
$\Trop(f_1+f_2+\cdots+f_k)=\min(\Trop(f_1),\Trop(f_2),\ldots,\Trop(f_k))$.
\item[(b)] Let and $n_1,n_2,\ldots,n_k$ be pairwise distinct
integers. Let  $f$ be the rational function $\GG_m\times S\to
\AA^1$ of the form $$f(c,s)=c^{n_1}f_1(s)+ c^{n_2}f_2(s)+\cdots+
c^{n_k}f_k(s) \ .$$
Then
$$[N](f)=(n_1,[N](f_1))\vee (n_2,[N](f_2))\vee \cdots\vee (n_k,[N](f_k)) $$
and $\Trop(f)=\min(n_1\cdot +\Trop(f_1),n_2\cdot
+\Trop(f_2),\ldots,n_k\cdot +\Trop(f_k))$. \item[(c)]  Let
$\mu_1,\mu_2,\ldots,\mu_k$ be pairwise distinct characters of an
algebraic torus $T'$. Let $f$ be the rational function $T'\times
S\to \AA^1$ of the form
$$f(t,s)=\mu_1(t)f_1(s)+
\mu_2(t)f_2(s)+\cdots \mu_n(t)f_k(s) \ .$$ Then
$$[N](f)=(\mu_1,[N](f_1))\vee (\mu_2,[N](f_2))\vee \cdots\vee (\mu_k,[N](f_k)) $$
and $$\Trop(f)=\min(\langle \mu_1,\bullet\rangle+\Trop(f_1),\langle
\mu_2,\bullet\rangle+\Trop(f_2),\ldots,\langle
\mu_k,\bullet\rangle+\Trop(f_k))\ .$$ \item[(d)] Assume additionally
that  $f_1,f_2,\ldots,f_k$ are positive (see Section
\ref{subsect:toric charts and positive structures}). Then
$$[N](f_1+f_2+\cdots+f_k)=[N](f_1)\vee [N](f_2)\vee \cdots\vee [N](f_k) $$
and
$\Trop(f_1+f_2+\cdots+f_k)=\min(\Trop(f_1),\Trop(f_2),\ldots,\Trop(f_k))$.
\end{enumerate}
\end{corollary}

Next, we extend the correspondence $f\mapsto \Trop(f)$ to the case
when $f$ is a rational morphism of algebraic tori.

Let now $S$ and $S'$ be split algebraic tori defined over $\QQ$,
and let $f:S\to S'$\linebreak be a rational morphism (note that
all regular morphisms $S\to S'$ are group homomorphisms up to
translations). Define the map $\Trop(f):X_\star(S)\to X_\star(S')$
by the formula
$$\left\langle\mu',\Trop(f)(\lambda)\right\rangle=\Trop(\mu'\circ f)(\lambda) $$
for any $\lambda\in X_\star(S)$, $\mu'\in X^\star(S')$, where the
rational function $\mu'\circ f:S\to \GG_m$ is considered as a
rational function $S\to \AA^1$.

Next, we will consider an appropriate category of split algebraic
tori for which the correspondence $f\mapsto \Trop(f)$ is a
functor.

\begin{definition}
\label{def:marked set} A {\it  marked set} is a pair   $(A,{\bf
0})$ where $A$ is a set,  ${\bf 0}\in A$ is a {\it marked point}.
Denote by ${\bf Set}_{\bf 0}$ the category whose objects are
marked sets and morphisms are structure-preserving maps.
\end{definition}

Recall from Section \ref{subsect:toric charts and positive
structures} that ${\mathcal T}_+$ is the category whose objects
are split algebraic tori (defined over $\QQ$), and arrows are
positive rational morphisms.

The following is a slightly modified form of the result of \cite[Section 2.4]{bk}.

\begin{theorem} [\cite{bk}]
\label{th:Trop} The correspondence $S\mapsto X_\star(S)$,
$f\mapsto \Trop(f)$ is a functor
$$\Trop:{\mathcal T}_+\to {\bf Set}_{\bf 0}\ .$$

\end{theorem}

\begin{remark} The theorem implies that for any positive birational isomorphism $f:S\osr S'$ such that
$f^{-1}$ is also positive the tropicalization $\Trop(f)$ is a
bijection $X_\star(S)\osr X_\star(S')$. The converse is not true:
if $f:\GG_m\times \GG_m\osr \GG_m\times \GG_m$ is as  in Remark
\ref{rem:non-positive inverse}, i.e., of $f^{-1}$ is not positive,
then  $\tilde f=\Trop(f):\ZZ^2\to \ZZ^2$ is given by $\tilde
f(\tilde x,\tilde y)=(\tilde x,\min(\tilde x,\tilde y))$. Clearly,
$\tilde f$ is not a bijection.
\end{remark}

We conclude this section with a discussion of the tropicalization of those (non-positive) functions which are compositions of
positive morphisms with certain rational functions.

\begin{definition} We say that a rational function $f:S\to \AA^1$ is {\it half-positive} if it can be presented as a
difference of two positive functions.

\end{definition}

Clearly, positive and each regular functions on $S$ are
half-positive. It is also clear that half-positive functions on
$S$ are closed under addition, subtraction, and multiplication --
they form a sub-ring of the field $Frac(S)$. Also, positive
functions act on half-positive ones by multiplication and
composition of each half-positive function with any positive
morphism is also a half-positive function.  The restriction of
each half-positive function on $S=(\GG_m)^\ell$ to
$(\RR_{>0})^\ell$ is  a well-defined map $(\RR_{>0})^\ell\to \RR$.
One can expect that this property is characteristic of
half-positive functions.

\begin{lemma}
\label{le:positive negative composition} Let $f:S'\to \AA^1$ be a
non-zero half-positive function. Then for any positive rational
morphism $h:S\to S'$, we have
\begin{equation}
\label{eq:inequality non-positive} \Trop(f\circ h)\ge
\Trop(f)\circ \Trop(h)\ .
\end{equation}

\end{lemma}

\begin{proof} Let us write $f$ as $f=f_+-f_-$ where both $f_+$ and $f_-$ are positive rational functions on $S'$.
If $f_+=0$ or $f_-=0$  then we have nothing to prove because Theorem \ref{th:Trop} guarantees the equality in \eqref{eq:inequality non-positive}.

Note that, by Theorem \ref{th:Trop}, we have $\Trop(f_+\circ
h)=\Trop(f_+)\circ \Trop(h)$.

Then by Lemma \ref{le:tropical function}(b),
\begin{align*}\Trop(f\circ h)&=\Trop(f_+\circ h-f_-\circ h)\ge
\min(\Trop(f_+\circ h),\Trop(f_-\circ h))\\ &=
\min(\Trop(f_+)\circ \Trop(h),\Trop(f_-)\circ \Trop(h))\\ &=
\min(\Trop(f_+),\Trop(f_-))\circ \Trop(h)=\Trop(f)\circ
\Trop(h)\end{align*} because
$\Trop(f)=\min(Trop(f_+),\Trop(f_-))$.

The lemma is proved.
\end{proof}

\begin{definition}
\label{def:positive equivalence}
We say that a  birational isomorphism $h:S\to S'$ is  a {\it positive equivalence} if both $h$ and $h^{-1}$ are positive.

\end{definition}

\begin{proposition}
\label{pr:positive equivalence composition} Let $f:S\to \AA^1$ be
a non-zero rational function on a split algebraic torus $S$. Then
for any positive equivalence $h:S'\to S$ of split algebraic tori,
one has
\begin{equation}
\label{eq:positive equivalence composition} \Trop(f\circ h)=
\Trop(f)\circ \Trop(h)\ .
\end{equation}

\end{proposition}

\begin{proof} First, let us prove the assertion for any half-positive function.
Since in this case $f\circ h$ is also half-positive, we obtain by
Lemma \ref{le:positive negative composition}. Therefore,
$$\Trop(f)=\Trop(f\circ h\circ h^{-1})\ge \Trop(f\circ h)\circ \Trop(h^{-1})  $$
$$=\Trop(f\circ h)\circ \Trop(h)^{-1}\ge \Trop(f)\circ \Trop(h)\circ \Trop(h)^{-1}=\Trop(f)\ , $$
that is, this chain of inequalities becomes a chain of equalities,
i.e.,  $\Trop(f)=\Trop(f\circ h)\circ \Trop(h)^{-1}$ which implies
\eqref{eq:positive equivalence composition}. Finally, note that
each non-zero rational function $f:S'\to \AA^1$ can be expressed
as a ratio $f=\frac{f'}{f''}$ where $f'$ and $f''$ are regular
hence half-positive and, therefore, $f'$ and $f''$ already satisfy
\eqref{eq:positive equivalence composition}. Then
\eqref{eq:positive equivalence composition} follows for $f$ as
well by Lemma \ref{le:tropical function}(a).

The lemma is proved.
\end{proof}

We say that two non-zero rational functions $f:S\to \AA^1$ and
$f':S'\to \AA^1$ are {\it positively equivalent} if there exists a
positive equivalence $h:S'\to S$ such that $f'=f\circ h$.

\begin{corollary}
\label{cor:positive class of functions} Let $f:S\to \AA^1$ be a
non-zero rational function. Then the isomorphism class of the
function $\Trop(f):X_\star(S)\to \ZZ$ in ${\bf Set}_{\bf 0}$
depends only the positive equivalence class of $f$.
\end{corollary}

\section{Kashiwara crystals, perfect bases, and associated crystals}

\label{sect:Kashiwara crystals, perfect bases, and crystal bases}
\subsection{Kashiwara crystals and normal crystals}
\label{subsect:Normal Kashiwara crystals}

First, we recall some definitions from \cite[Appendix]{bk}.

\begin{definition}
\label{def:partial bijection} A {\it partial bijection} of sets
$\tilde f:\tilde A\to \tilde B$ is a bijection $\tilde A'\to
\tilde B'$ of subsets $\tilde A'\subset \tilde A$, $\tilde
B'\subset \tilde B$. We denote the subset $\tilde A'$ by $\dom(f)$
and  the subset $\tilde B'$ by $\ran(\tilde f)$.
\end{definition}

\smallskip

The inverse $\tilde f^{-1}$ of a partial bijection $\tilde
f:\tilde A\to \tilde B$ is the inverse bijection $\ran(\tilde
f)\to \dom(\tilde f)$. The composition  $\tilde g\circ \tilde f$
of partial bijections $\tilde f:\tilde A\to \tilde B,$\linebreak $
\tilde g: \tilde B\to \tilde C$ is  a partial bijection with
$\dom(\tilde g\circ \tilde f)=\dom(\tilde f)\cap \tilde
f^{-1}(\ran(\tilde f)\cap \dom(\tilde g))$ and $\ran(\tilde g\circ
f)=\tilde g(\ran(\tilde f)\cap\ dom(\tilde g))$. In particular,
for any partial bijection $f:\tilde B\to \tilde B$ and $n\in \ZZ$
the $n$-th power $\tilde f^n$ is a partial bijection $\tilde B\to
\tilde B$.

For a partial bijection $\tilde f:\tilde A\to \tilde B$, we will
sometimes use a notation: $\tilde f(\tilde a)\in \tilde B$ if and
only if $\tilde a\in \dom(\tilde f)$.

\begin {definition}
\label{def:Kashiwara crystals} Following  \cite{k93} and
\cite[Appendix]{bk}, we say that a {\it Kashiwara crystal} is a
$5$-tuple $\BB=(\tilde B,\tilde \gamma,\tilde \varphi_i,\tilde
\varepsilon_i,\tilde e_i| i\in I)$, where $\tilde B$ is a set,
$\tilde \gamma: \tilde B\to X_\star(T)$ is a map,  $\tilde
\varphi_i,\tilde \varepsilon_i:\tilde B\to \ZZ$ are functions, and
each $\tilde e_i:\tilde B\to \tilde B$, $i\in I$, is a partial
bijection such that either $\tilde e_i=\id_{\tilde B}$, $\tilde
\varepsilon_i=\tilde \varphi_i=-\infty$  or:
$$\tilde \varphi_i(\tilde b)-\tilde \varepsilon_i(\tilde b)=\left<\alpha_i,\tilde \gamma(\tilde b)\right> $$
for all $\tilde b\in \tilde B$  and
$$\tilde \gamma(\tilde e_i^{\,n}(\tilde b))=\tilde \gamma(\tilde b)+n\alpha_i^\vee$$
whenever $\tilde e_i^{\,n}(\tilde b)\in \tilde B$.
\end{definition}

We define $\Supp \BB=\{i\in I:\tilde e_i\ne \id\}$ and call it the
{\it support} of $\BB$.

\begin{example}
\label{ex:trivial combinatorial crystal} Along the lines of Example \ref{ex:trivial geometric crystal},  $X_\star(T)$ is a {\it
trivial} Kashiwara crystal with $\tilde \gamma=\id$, $\tilde
\varepsilon_i=\tilde \varphi_i=-\infty$, $\tilde e_i=\id$ for all
$i\in I$. The support of this trivial crystal is the empty set
$\emptyset$. Another example of a trivial Kashiwara crystal is any
subset of $X_\star(T)$, in particular, the single point $0\in
X_\star(T)$.  Yet another example of a Kashiwara crystal on $X_\star(T)$ is a $5$-tuple ${\mathcal B}_T=(X_\star(T),id_{X_\star(T)},\tilde \varphi_i,\tilde \varepsilon_i,\tilde e_i|i \in I)$,
where $\tilde e_i^n(\mu)=\mu+n\alpha_i^\vee$ for all $n\in \ZZ_{\ge 0}$, $\mu\in X_\star(T)$, $i\in I$ and $\widetilde \varepsilon_i,\widetilde \varphi_i\in X^\star(T)$ are such that $\langle \tilde \varepsilon_i,\alpha_i^\vee\rangle=-1$, $ \tilde \varphi_i=\tilde \varepsilon_i+\alpha_i^\vee$ for all $i\in I$.
\end{example}

\begin{definition}
\label{def:homomorphism Kashiwara crystals} A {\it homomorphism of
Kashiwara crystals} $\tilde h:\BB\to \BB'$ is a pair $(h,J)$,
where $h:\tilde B\to \tilde B'$ is a map and  $J\subset \Supp
\BB\cap \Supp \BB'$ such that $h\circ \tilde e_i^{\,n}=\tilde e_i^{\,n}\circ h$
for $i\in \Supp \BB$, $n\in \ZZ$ and
$$\tilde \varphi'_j\circ h=\varphi_j,~\tilde \varepsilon'_j\circ h=\varepsilon_j $$
for all $j\in J$ (we  will refer to this $J$ as the {\it support} of $\tilde h$ and denote by $\Supp \tilde h$).

The composition of homomorphisms $(\tilde h,J):\BB\to \BB'$ and
$(\tilde {h'},J'):\BB'\to \BB''$ is defined by $(\tilde
{h'},J')\circ (\tilde h,J):=(\tilde {h'}\circ \tilde h,J\cap J')$
(i.e., $\Supp \tilde
{h'}\circ \tilde h:= \Supp \tilde h\cap \Supp \tilde {h'}$).

\end{definition}

\begin{remark} In the case when $\Supp \BB=\Supp \BB'=\Supp \tilde
h$, our definition of a homomorphism of Kashiwara crystals $\tilde h:\BB\to \BB'$
 coincides with the definition of {\it strict} homomorphism of crystals
from the original work \cite{k93}.
\end{remark}

In what follows, we consider the category whose objects are Kashiwara crystals and morphisms are homomorphisms of Kashiwara crystals.

\begin{definition}
\label{def:product combinatorial} Given Kashiwara crystals
$\BB=(\tilde B,\tilde \gamma,\tilde \varphi_i,\tilde
\varepsilon_i,\tilde e_i|i\in I)$ and $\BB'=(\tilde B,\tilde \gamma',\tilde
e_i,\tilde \varphi'_i,\tilde \varepsilon'_i|i\in I)$, the {\it
product} $\BB\times \BB'$ of Kashiwara crystals is a $5$-tuple
$$(\tilde B\times \tilde B',\tilde \gamma'',\tilde \varphi''_i,\tilde \varepsilon''_i,\tilde e_i| i\in I)\ ,$$
where $\Supp \BB\times \BB':=\Supp \BB\cup \Supp \BB'$, $\tilde
\gamma''(\tilde b,\tilde b')=\tilde \gamma(\tilde b)+\tilde
\gamma(\tilde b')$ and for each $i\in \Supp \BB\times \BB',$ one
has
$$\tilde \varphi_i(\tilde b,\tilde b')=\max(\tilde \varphi_i(\tilde b),\tilde \varphi_i(\tilde b')+\langle\alpha_i,
\tilde \gamma(\tilde b)\rangle),\, \tilde \varepsilon_i(\tilde
b,\tilde b')= \max(\tilde \varepsilon'_i(\tilde b'), \tilde
\varepsilon_i(\tilde b)-\langle\alpha_i,\tilde \gamma'(\tilde
b')\rangle)\ ,$$ and each partial bijection $e_i^{\,n}:\tilde
B\times \tilde B'\to \tilde B\times \tilde B'$, $n\in \ZZ$, $i\in
I$ is given by
$$\tilde e_i^{\,n}(\tilde b,\tilde b')=(\tilde e_i^{\,n_1}(\tilde b),\tilde e_i^{\,n_2}(\tilde b'))\ ,$$
for $(\tilde b,\tilde b')\in \tilde B\times \tilde B'$, where
\begin{equation}
\label{eq:tensor product of Kashiwara crystals1}
n_1=\max(\tilde \varepsilon_i(\tilde b),\tilde \varphi'_i(\tilde b'))-\max(\tilde \varepsilon_i(\tilde b)-n,\tilde \varphi'_i(\tilde b')),
\end{equation}
\begin{equation}
\label{eq:tensor product of Kashiwara crystals2}
~n_2=\max(\tilde \varepsilon_i(\tilde b),\tilde \varphi'_i(\tilde b')+n)-\max(\tilde \varepsilon_i(\tilde b),\tilde \varphi'_i(\tilde b')) \ .
\end{equation}
\end{definition}

It is easy to see that and the product is associative.

\begin{remark}
\label{rem:kashiwara product inverted} The above definition agrees
with   Kashiwara's original definition up to the permutation of
factors, i.e., the product $\BB\times \BB'$ equals the {\it tensor
product} $\BB'\otimes \BB$ in the notation of \cite{k93}. The
reason for this is functoriality  of the transition from the
geometric crystals to Kashiwara crystals (see \eqref{eq:positive
product geometric decorated} below).
\end{remark}

\begin{example}
\label{ex:projection product to lattice} For a Kashiwara crystal
$\BB$, the product $X_\star(T)\times \BB$ (where $X_\star(T)$ is
considered the trivial Kashiwara crystal as in Example
\ref{ex:trivial combinatorial crystal}) is a Kashiwara crystal
with $\Supp(X_\star(T)\times \BB)=\Supp \BB$. The projection to
the first factor is a homomorphism of Kashiwara crystals
$X_\star(T)\times \BB\to X_\star(T)$. The support of the
homomorphism is $\emptyset$.
\end{example}

\begin{claim}
\label{cl:sub-normal product} Let $\BB$ and $\BB'$ be Kashiwara
crystals and let   $i\in \Supp \BB\cup \Supp \BB'$, and  $n\in
\ZZ\setminus \{0\}$. Then for any $\tilde b\in \BB$ such that
$\tilde e_i^{\,n}(\tilde b)\in \BB$ and any $\tilde b'\in \BB',$
one has
\begin{enumerate}
\item[(a)] $\tilde e_i^{\,n}(\tilde b',\tilde b)=(\tilde b',\tilde
e_i^{\,n}(\tilde b))$ if and only if $\tilde \varphi_i(\tilde
b)\ge \tilde \varepsilon'_i(\tilde b')$, $n\ge \tilde
\varepsilon'_i(\tilde b')-\tilde \varphi_i(\tilde b)$. \item[(b)]
$\tilde e_i^{\,n}(\tilde b,\tilde b')=(\tilde e_i^{\,n}(\tilde
b),\tilde b')$ if and only if  $ \tilde \varepsilon_i(\tilde b)\ge
\tilde \varphi_i'(\tilde b')$, $n\le \tilde \varepsilon_i(\tilde
b)-\tilde \varphi_i'(\tilde b')$.
\end{enumerate}
\end{claim}


Following \cite{k93}, for each Kashiwara crystal $\BB=(\tilde
B,\tilde \gamma,\tilde \varphi_i,\tilde
\varepsilon_i,\tilde e_i|i\in I)$, define the {\it opposite} Kashiwara
crystals $\BB^{op}=(\tilde B,-\tilde \gamma,\tilde \varepsilon_i,\tilde
\varphi_i,\tilde e_i^{-1}|i \in I)$.

The following is a combinatorial analogue of Claim \ref {cl:op geometric}.

\begin{claim} \cite{k93}
\label{cl:op combinatorial} The correspondence $\BB\mapsto
\BB^{op}$ defines an involutive covariant functor from the
category of Kashiwara crystals into itself. This functor reverses
the product, i.e., $(\BB\times \BB')^{op}$ is naturally isomorphic
to ${\BB'}^{op}\times \BB^{op}$. On the underlying sets this
isomorphism is the permutation of factors $\tilde B\times \tilde
B'\osr \tilde B'\times \tilde B$.

\end{claim}

\begin{definition}
\label{def:restriction of Kasiwara crystal} For each Kashiwara
crystal $\BB=(\tilde B,\tilde \gamma,\tilde
\varphi_i,\tilde \varepsilon_i,\tilde e_i|i\in I)$ and a subset $\tilde
B'\subset \tilde B$, we define  the {\it sub-crystal}
$\BB|_{\tilde B'}$  of $\BB$ as follows:
$$\BB|_{\tilde B'}=(\tilde B',\tilde \gamma|_{\tilde B'},\tilde \varphi_i|_{\tilde B'},
\tilde \varepsilon_i|_{\tilde B'},\tilde e_i|_{\tilde B'}~|i\in I) \ ,$$ where $\tilde
e_i|_{\tilde B'}$ is the partial bijection $\tilde B'\to \tilde
B'$ obtained by the restriction of the partial bijection $\tilde
e_i:\tilde B\to \tilde B$ to $\tilde B'$ via $\dom(\tilde
e_i|_{\tilde B'})=\{\tilde b'\in \tilde B'\cap \dom(\tilde e_i):
\tilde e_i(\tilde b)\in \tilde B'\}$.

\end{definition}

Clearly, the natural embedding $\tilde B'\subset \tilde B$ defines
an injective homomorphism of Kashiwara crystals
$\BB|_{\tilde B'} \hookrightarrow \BB$. If $\tilde B'$ is not
empty, then $\Supp \BB'=\Supp \BB$ and the support of the
homomorphism also equals $\Supp \BB$ (see Definition
\ref{def:homomorphism Kashiwara crystals}).

\begin{definition} We say that a Kashiwara crystal
$\BB=(\tilde B,\tilde \gamma,\tilde \varphi_i,\tilde
\varepsilon_i,\tilde e_i|i\in I)$ is {\it connected} if for any $\tilde
b,\tilde b'\in \tilde B$ there exist a sequence
$(i_1,\ldots,i_\ell)\in I^\ell$ and a sequence
$(n_1,\ldots,n_\ell)\in \ZZ^\ell$ such that
$$\tilde b'=\tilde e_{i_1}^{\,n_1}\cdots \tilde e_{i_\ell}^{\,n_\ell}(\tilde b) \ .$$

\end{definition}

\begin{claim}
\label{cl:disjoint}
Each Kashiwara crystal $\BB$ equals to the disjoint union of its  connected sub-crystal.

\end{claim}

For each Kashiwara crystal $\BB$, we define functions
$\ell_i,\ell_{-i}:\tilde B\to \ZZ_{\ge 0}\sqcup \{+\infty\}$ by
\begin{equation}
\label{eq:eplus minus i}
\ell_i(\tilde b)=\max\{n\ge 0:\tilde e_i^{\,n}(\tilde b)\in \tilde B\},~\ell_{-i}(\tilde b)=\max\{n\ge 0:\tilde e_i^{\,-n}(\tilde b)\in \tilde B\}
\end{equation}
for $\tilde b\in \tilde B$.

\begin{definition}
\label{def:highest weight element} Given a Kashiwara crystal
$\BB$, we say that an element $\tilde b\in \tilde B$ is a {\it
highest} (resp.~{\it lowest}) {\it weight element}  if
$\ell_i(\tilde b)=0$ (resp.~$\ell_{-i}(\tilde b)=0$) for all $i\in
\Supp \BB$. Denote by $\BB^+$ (resp.~by $\BB^-$) the set of
highest (resp.~lowest) weight elements of $\BB$.
\end{definition}

\begin{claim}  $(\BB^{op})^+=\BB^-$ and $(\BB^{op})^-=\BB^+$ for each
 Kashiwara crystal $\BB$ and
$(\BB\times \BB')^+\subset \BB^+\times \BB{'^+}$,
 $(\BB\times \BB')^-\subset \BB^-\times \BB{'^-}$ for any Kashiwara
 crystals $\BB$ and~$\BB'$.

\end{claim}

\begin{definition}
\label{def:upper sub-normal} We say that  $\BB=(\tilde B,\tilde
\gamma,\tilde \varphi_i,\tilde \varepsilon_i,\tilde e_i|i\in I)$
is {\it upper subnormal} (resp.~{\it lower subnormal}) if
$\varepsilon_i(\tilde b)\ge 0$ (resp.~$\varphi_i(\tilde b)\ge 0$)
for all $i\in \Supp \BB$, $\tilde b\in \tilde B$. If $\BB$ is both
upper and lower subnormal, we will refer to it simply as {\it
subnormal}.

\end{definition}

\begin{claim} A Kashiwara  crystal $\BB$ is upper subnormal (resp.~lower subnormal) if
and only if
$\ell_i\le \tilde \varepsilon_i$ (resp.~$\ell_{-i}\le \tilde
\varphi_i$) for each $i\in \Supp \BB$.
\end{claim}

\begin{claim} Let $\BB$, $\BB'$ be Kashiwara crystals such that $\Supp \BB\supseteq
\Supp \BB'$. Then \begin{enumerate}\item[(a)] If $\BB$ is lower
sub-normal, then $\BB\times \BB'$ is also lower subnormal.
\item[(b)] If $\BB$ is upper sub-normal,  then $\BB'\times \BB$ is
also upper subnormal.\end{enumerate}
\end{claim}

In particular, subnormal Kashiwara crystals of a given support $J\subset I$ from a
 monoidal subcategory in the category of all Kashiwara crystals.

For each  Kashiwara crystal $\BB=(\tilde B,\tilde \gamma,\tilde
e_i,\tilde \varphi_i,\tilde \varepsilon_i|i\in I),$ define the
subset
$$\overline {\tilde B}=\{\tilde b\in \tilde B:\tilde \varphi_i(\tilde b)\ge 0, \tilde \varepsilon_i(\tilde b)\ge 0~\forall~i\in \Supp \BB,
~\tilde b\in \tilde B\}\ .$$
We denote by $\overline \BB$ the restriction of $\BB$ to the subset
$\overline {\tilde B}$ and call $\overline \BB$ the {\it
sub-normalization} of $\BB$. Clearly, $\overline \BB$ is always
subnormal. Also, a Kashiwara  crystal $\BB$ is subnormal if and only
if $\overline \BB=\BB$.

%

\begin{lemma}
\label{le:product trivial} Let $\BB$ and $\BB'$ be Kashiwara
crystals and let ${\bf 0}\in \BB'$ be such a point that $\tilde
\varepsilon'_i({\bf 0})=\tilde \varphi'_i({\bf 0})=0$ for $i\in
\Supp \BB'$.   Then: \begin{enumerate} \item[(a)] The
correspondence $\tilde b\mapsto ({\bf 0},\tilde b)$ (resp.~$\tilde
b\mapsto (\tilde b,{\bf 0})$) defines an injective\linebreak
homomorphism of subnormal Kashiwara crystals $\overline \BB
\hookrightarrow \overline {\BB'}\times \overline \BB$
(resp.\linebreak $\overline \BB \hookrightarrow \overline
\BB\times \overline {\BB'}$). The support of each of these
homomorphisms is $\Supp \BB$.
 \item[(b)] If $\BB$ is upper subnormal, then the correspondence
$\tilde b\mapsto (\tilde b,{\bf 0})$ defines an injective
homomorphism of  Kashiwara crystals  $\BB \hookrightarrow
\BB\times \BB'$. The support of this homomorphism is $\Supp \BB$.
 \item[(c)] If $\BB$ is lower subnormal, the correspondence
$\tilde b\mapsto ({\bf 0},\tilde b)$ defines an injective
homomorphism of  Kashiwara crystals $\BB \hookrightarrow
\BB'\times \BB$. The support of this homomorphism is $\Supp \BB$.
\end{enumerate}
\end{lemma}

\begin{proof} We need the following special case of
Claim \ref{cl:sub-normal product}.

\begin{claim}
\label{cl:sub-normal product zero} Let $\BB$ and $\BB'$ be
Kashiwara crystals and let ${\bf 0}\in \BB'$ be such a point that
$\tilde \varepsilon'_i({\bf 0})=\tilde \varphi'_i({\bf 0})=0$ for
all $i\in \Supp \BB'$. Then for each  $\tilde b\in \BB$,  $i\in
\Supp \BB\cap \Supp \BB'$, one has
$$\tilde \varphi_i({\bf 0},\tilde b)=\max(0,\tilde \varphi_i(\tilde b)),
\tilde \varepsilon_i({\bf 0},\tilde b)= \max(\tilde \varepsilon_i(\tilde b),
\tilde \varepsilon_i(\tilde b)-\tilde \varphi_i(\tilde b))\ ,$$
$$\tilde \varphi_i(\tilde b,{\bf 0})=\max(\tilde \varphi_i(\tilde b),\tilde \varphi_i(\tilde b)-\tilde \varepsilon_i(\tilde b)),
~\tilde \varepsilon_i(\tilde b,{\bf 0})= \max(0,\tilde
\varepsilon_i(\tilde b))\ ,$$ and for $n\in \ZZ\setminus \{0\}$,
$\tilde b\in \BB$, $i\in \Supp \BB$ such that $\tilde
e_i^{\,n}(\tilde b)\in \BB$, one has: \begin{enumerate}\item[(i)]
$\tilde e_i^{\,n}({\bf 0},\tilde b)=({\bf 0},\tilde
e_i^{\,n}(\tilde b))$ if and only if $\tilde \varphi_i(\tilde
b)\ge \tilde 0$ and $n\ge -\tilde \varphi_i(\tilde b)$.
\item[(ii)] $\tilde e_i^{\,n}(\tilde b,{\bf 0})=(\tilde
e_i^{\,n}(\tilde b),{\bf 0})$ if and only if $\tilde
\varepsilon_i(\tilde b)\ge 0$ and $n\le \tilde
\varepsilon_i(\tilde b)$\end{enumerate}
\end{claim}

Prove (a) now. For each $\tilde b\in \overline \BB$, one has
$$\tilde \varphi_i({\bf 0},\tilde b)=\tilde \varphi_i(\tilde b),
\tilde \varepsilon_i({\bf 0},\tilde b)= \tilde \varepsilon_i(\tilde b),\tilde \varphi_i(\tilde b,{\bf 0})
=\tilde \varphi_i(\tilde b),~\tilde \varepsilon_i(\tilde b,{\bf 0})=\tilde \varepsilon_i(\tilde b)\ ,$$

Using Claim \ref{cl:sub-normal product zero}(i), we obtain
$$\tilde e_i^{\,n}({\bf 0},\tilde b)=({\bf 0},\tilde e_i^{\,n}(\tilde b))$$
for  $\tilde b\in \overline \BB$, $i\in \Supp \BB$, $n\in \ZZ$ if
and only if $\tilde e_i^{\,n}(\tilde b)\in \overline \BB$ and
$n\ge - \tilde \varphi_i(\tilde b)$. But the latter inequality
holds automatically because $\tilde \varphi_i(e_i^{\,n}(\tilde
b))=\tilde \varphi_i(\tilde b)+n\ge 0$ if  $\tilde
e_i^{\,n}(\tilde b)\in \overline \BB$.

Similarly, using Claim \ref{cl:sub-normal product zero}(ii), we obtain
$$\tilde e_i^{\,n}(\tilde b,{\bf 0})=(\tilde e_i^{\,n}(\tilde b),{\bf 0})$$
for  $\tilde b\in \overline \BB$, $i\in \Supp \BB$, $n\in \ZZ$ if
and only if $\tilde e_i^{\,n}(\tilde b)\in \overline \BB$.
Finally, $\tilde \gamma(\tilde b,{\bf 0})=\gamma({\bf 0},\tilde
b)=\tilde \gamma(\tilde b)$. This proves (a).

Prove (b). Indeed, if $\BB$ is upper normal, then Claim
\ref{cl:sub-normal product zero} guarantees that $\tilde
\varphi_i(\tilde b,{\bf 0})=\tilde \varphi_i(\tilde b)$, $\tilde
\varepsilon_i(\tilde b,{\bf 0})=\tilde \varepsilon_i(\tilde b)\ge
0$,  and $\tilde e_i^{\,n}(\tilde b,{\bf 0})=(\tilde
e_i^{\,n}(\tilde b),{\bf 0})$ for all $n\le \tilde
\varepsilon_i(\tilde b)$. Finally, $\tilde \gamma(\tilde b,{\bf
0})=\gamma({\bf 0},\tilde b)=\tilde \gamma(\tilde b)$. This proves
(b).

Part (c) follows.
The lemma is proved.
\end{proof}

We conclude the section with definitions and results related to
normal\linebreak Kashiwara crystals.

\begin{definition}
\label{def:normal crystals} Following \cite{k93}, we say that
$\BB=(\tilde B,\tilde \gamma,\tilde \varphi_i,\tilde
\varepsilon_i,\tilde e_i| i\in I)$ is {\it upper normal}
(resp.~{\it lower normal}) if $\tilde \varepsilon_i=\ell_i$ (resp.
$\tilde \varphi_i=\ell_{-i}$) for all $i\in \Supp \BB$ (in the
notation \eqref{eq:eplus minus i}).

A Kashiwara crystal $\BB$ is  {\it normal} if  it is both lower
normal and upper normal and  $\Supp \BB=I$.

\end{definition}

\begin{lemma}
\label{le:normal crystals act on upper normals}
If $\BB$ is an upper (resp.~lower) normal Kashiwara
crystal and $\BB'$ is normal, then $\BB\times B'$ is upper normal (resp. $\BB'\times \BB$ is lower normal).
\end{lemma}

\begin{proof} For the upper normal crystals, the assertion follows from the formula
$$\tilde e_i^{\,\tilde \varepsilon_i(\tilde b,\tilde b')+n}(\tilde b,\tilde b')=
(\tilde e_i^{\,n_1}(\tilde b),\tilde e_i^{\,n_2}(\tilde b'))\ ,$$
for $(\tilde b,\tilde b')\in \BB\times \BB'$, $n\in \ZZ$, $i\in \Supp \BB\cap
\Supp \tilde \BB'$, where
$$n_1=\min( \max(0,\tilde \varepsilon_i(\tilde b)-\tilde \varphi'_i(\tilde b')),\varepsilon'_i(\tilde b')+|\tilde \varepsilon_i(\tilde b)-\tilde \varphi'_i(\tilde b')|+n),$$
$$n_2=\max(\tilde \varepsilon'_i(\tilde b')+n, \min(0, \tilde
\varepsilon_i(\tilde b)-\tilde \varphi'_i(\tilde b')) \ .$$ This
implies that if  $\BB$ is upper normal and $\BB'$ is normal, then
$n_1=\max(0,\tilde \varepsilon_i(\tilde b)-\tilde
\varphi'_i(\tilde b'))\le \varepsilon_i(\tilde b)$, $n_2=\tilde
\varepsilon'_i(\tilde b')+n$; hence $\tilde e_i^{\,\tilde
\varepsilon_i(\tilde b,\tilde b')}(\tilde b,\tilde b')\in
\BB\times \BB'$ and $\tilde e_i^{\,\tilde \varepsilon_i(\tilde
b,\tilde b')+1}(\tilde b,\tilde b')\notin \BB\times \BB'$. This
proves the assertion for the upper normal $\BB$. If $\BB$ is lower
normal, the assertion follows from the above and Claim \ref{cl:op
combinatorial} by applying $\BB\mapsto \BB^{op}$ and
$\BB'\mapsto{\BB'}^{op}$.
\end{proof}

Thus   normal Kashiwara crystals and
their homomorphisms form a monoidal category (see also \cite{k93}). This category acts from the right (resp. from the left) on the category of upper (resp. lower) normal crystals.

The following result demonstrates that each upper normal (resp. lower) crystal is semisimple.

\begin{claim} Let  $\BB$ be an upper normal crystal, and let $\tilde B'$ be a subset of $\tilde B$
such that the sub-crystal $\BB|_{\tilde B'}$ is also upper normal. Then the complement $\BB|_{\tilde B\setminus \tilde B'}$ is upper normal as well.
\end{claim}

Note that if $\{\BB_k\}$ is any family of upper normal
sub-crystals of an upper normal crystal, then both the
intersection $\bigcap_k \BB_k$ and the union $\bigcup_k \BB_k$ are
also upper normal sub-crystals. This prompts the following
definition.

\begin{definition}
Given an upper normal crystal $\BB$, for any subset $\tilde C\subset
\BB$ denote by $\BB[\tilde C]$ the intersection of all upper normal
sub-crystals of $\BB$ containing $\tilde C$, i.e., $\BB[\tilde C]$
is the smallest upper normal sub-crystal of $\BB$ containing $\tilde
C$. We will refer to $\BB[\tilde C]$ as the upper normal sub-crystal
of $\BB$ {\it generated upward} by $\tilde C$.

\end{definition}

\begin{lemma} For any $\tilde C\subset \BB$, the upper normal sub-crystal $\BB[\tilde C]$ of $\BB$ consists of all elements of the form
\begin{equation}
\label{eq:upper sub-crystal generation}
\tilde e_{i_\ell}^{\,n_\ell}\tilde e_{i_{\ell-1}}^{\,n_{\ell-1}} \cdots \tilde e_{i_1}^{\,n_1}(\tilde c)
\end{equation}
for all $\tilde c_0\in \tilde C$ and any $i_1,\ldots,i_\ell\in
\Supp \BB$, $n_1,\ldots,n_\ell\in \ZZ_{\ge 0}$ are such that
$$n_k\le \tilde \varepsilon_{i_k}(\tilde e_{i_{k-1}}^{\,n_{k-1}} \cdots \tilde e_{i_1}^{\,n_1}(\tilde c))$$
for $k=1,2,\ldots, \ell$.
\end{lemma}

\begin{proof} Denote by $\BB_0$ the sub-crystal of $\BB$ which consists of all elements of the form \eqref{eq:upper sub-crystal generation}.
Clearly, if $\BB'$ is an upper normal sub-crystal of $\BB$
containing each $\tilde c\in \tilde C$, then $\BB'$ also contains
$\BB_0$. It remains to show that $\BB_0$ is upper normal. But, by
definition of $\BB_0$ for each $\tilde b_0\in \BB_0$, one has
$$\tilde e_i^{\,n}(\tilde b_0)\in \BB_0$$
for each $i\in \Supp \BB$, $n\in \ZZ_{\ge 0}$ such that $n\le
\varepsilon_i(\tilde b_0)$. This verifies Definition
\ref{def:normal crystals} and proves the upper normality of
$\BB_0$. Thus, $\BB_0$ is the smallest upper normal crystal
containing $\tilde C$. The lemma is proved.
\end{proof}

Note that if $\tilde C=\{\tilde b\}$ is a single element, then
$\tilde c$ is the {\it lowest} weight element in $\BB[\tilde c]$.
Moreover, $\BB[\tilde C]=\bigcup_{\tilde c\in \tilde C}
\BB[\{\tilde c\}]$.

\begin{claim}
\label{cl:normal sub-crystals are equal} Let $\BB$ be a Kashiwara
crystal and let $\BB_1$ and $\BB_2$ be its normal sub-crystals.
Then $\BB_1\cap \BB_2$ is also normal. If, in addition, $\BB_1$
and $\BB_2$ are connected, then either  $\BB_1\cap
\BB_2=\emptyset$ or $\BB_1=\BB_2$.
\end{claim}

\begin{claim}
\label{cl:connected normal crystal}
Any  normal  crystal is a disjoint union of its normal connected sub-crystals.

\end{claim}

\subsection{Perfect bases, upper normal crystals, and  associated crystals}
\label{subsect:Perfect bases, upper normal crystals, and  crystal
bases} Let $I$ be a finite set of indices, let $\Lambda$  be a
lattice, and let $\{\alpha_i, i\in I\}$ be a subset of $\Lambda$.
Let $\{\alpha_j^\vee,j\in I\}$ be a subset of the dual lattice
$\Lambda^\vee=\Lambda^*$ (e.g., $\Lambda=X^\star(T)$ is the weight
lattice of $G$, $\Lambda^\vee=X_\star(T)$ is the co-weight lattice
of $G$,  and $\alpha_j^\vee$, $\alpha_i$ are respectively co-roots
and roots of $G$). Denote by $\hat \bb$ the Lie algebra generated
by $h_i,e_i$, $i\in I$ subject to the relations
$$[h_i,e_j]=\langle \alpha_i, \alpha_j^\vee\rangle  e_j$$
for all $i,j\in I$, where $\langle\bullet,\bullet\rangle:
\Lambda\times \Lambda^\vee \to \ZZ$ is the evaluation pairing. By
definition, $\hat \bb$ is a Borel sub-algebra of a generalized
Kac-Moody Lie algebra and $\hat \bb=\hh^\vee \ltimes \hat \nn$,
where $\hh^\vee$ is an Abelian Lie algebra with the basis
$\{h_i,i\in I\}$ and $\hat \nn$ is the free Lie algebra generated by
all $e_i$.

Denote by $\BBB$ the universal enveloping algebra $U(\bb)$, that is,
$\BBB$ is an associative algebra generated over $\CC$ by $h_i,e_i$,
$i\in I$ subject to the relations
$$h_ie_j-e_j h_i=\langle \alpha_i, \alpha_j^\vee\rangle e_j$$
for all $i,j\in I$.


We say that a $\BBB$-module $V$ is {\it locally finite} if:

\noindent $\bullet$  $V$  has a $\Lambda^\vee$-weight decomposition:
\begin{equation}
\label{eq:general weight decomposition}
V=\bigoplus\limits_{\mu\in \Lambda^\vee} V(\mu)\ .
\end{equation}

\noindent $\bullet$ For any $v\in V$, the cyclic sub-module
$\BBB(v)$ of $V$ is finite-dimensional.

\medskip

For any non-zero vector $v\in V$ and $i\in I$, denote by
$\ell_i(v)$ the smallest positive integer $\ell$ such that
$e_i^{\,\ell+1}(v)=0$. For $v=0$ we will use the convention
$\ell_i(0)=-\infty$ for all $i$.

Given a locally finite $\BBB$-module $V$, for each sequence
$\ii=(i_1,\ldots,i_{m})\in I^m$, $m\ge 1$, we define a binary
relation $\preceq_\ii$ on $V\setminus \{0\}$  as follows.

If $m=1$,  $\ii=(i)$, write $v\preceq_\ii v'$ if and only if
$\ell_i(v)\le \ell_i(v')$. For $\ii=(i;\ii')$, write
$v\preceq_\ii v'$ if and only if either $\ell_i(v)<\ell_i(v')$ or
$\ell_i(v)=\ell_i(v')$, $e_i^{\ell_i(v)}(v)\preceq_{\ii'}
e_i^{\ell_i(v)}(v')$.

\begin{claim} For each $\ii$, the relation $\preceq_\ii$ is a pre-order on $V\setminus \{0\}$.

\end{claim}

Therefore, we define the equivalence relation $\equiv_\ii$ on $V\setminus \{0\}$ by setting
$v\equiv_\ii v'$ if  $v\preceq_\ii v'$ and $v'\preceq_\ii v$. Also we will write $v\prec_\ii v'$ if $v\preceq_\ii v'$ and $v\not \equiv_\ii v'$.

\begin{claim}
\label{cl:general equiv} For any sequence $\ii\in I^m$, one has:
\begin{enumerate}\item[(a)] For any non-zero $v\in V$, the set $\{v'\in V:v'\prec_\ii
v\}$ is a subspace of $V$. \item[(b)] If $v \not \equiv_\ii v'$,
then $v+v' \equiv_\ii
\begin{cases} v & \text{if $v' \prec_\ii v$}\\
v' & \text{if $v \prec_\ii v'$}
\end{cases}$\ .
\end{enumerate}
\end{claim}


For each $i\in I$ and $\ell\ge 0$,  define the subspace
$$V_i^{<\ell}:=\{v\in V:\ell_i(v)<\ell\}=\{v\in V:e_i^\ell(v)=0\} \ .$$

We say that a basis ${\bf B}$ of a locally finite $\BBB$-module $V$
is a {\it weight basis} if ${\bf B}$ is compatible with the weight
decomposition \eqref{eq:general weight decomposition}, i.e., ${\bf
B}(\mu):=V(\mu)\cap {\bf B}$ is a basis of $V(\mu)$ for any $\mu\in
\Lambda^\vee$.

\begin{definition}
\label{def:general upper perfect basis} We say that a weight basis
${\bf B}$ in a locally finite $\BBB$-module $V$ is {\it perfect}
if for each $i\in I$ there is a partial bijection $\tilde e_i:{\bf
B}\to {\bf B}$ (i.e., a bijection of a subset of ${\bf B}$ onto
another subset of ${\bf B}$, see also Definition \ref{def:partial
bijection}  above) such that $\tilde e_i(b)\in {\bf B}$ if and
only if  $e_i(b)\ne 0$, and in the latter case one has
\begin{equation}
\label{eq:general upper asociated crystal}
e_i(b)\in \CC^\times\cdot \tilde e_i(b)+ V_i^{<\ell_i(b)-1} \ .
\end{equation}

Following \cite{lu}, we refer to a pair $(V,{\bf B})$, where $V$ is
a locally finite $\BBB$-module and ${\bf B}$ is a perfect basis of
$V$, as a {\it based} $\BBB$-module.

\end{definition}

Replacing the lattice $X_\star(T)$ by $\Lambda^\vee$, we extend Definition \ref{def:Kashiwara crystals}
of Kashiwara crystals to the setup of this section.

Recall from  Definition \ref{def:normal crystals} that a Kashiwara crystal $\BB$ is upper normal if
$$\varepsilon_i(\tilde b)=\ell_i(\tilde b)=\max\{n:\tilde e_i^{\,n}(\tilde b)\in B\} $$
for all $\tilde b\in \tilde B$.

\begin{claim}
\label{cl:perfect asociated crystal} Each based $\BBB$-module
$(V,{\bf B})$  defines an upper normal Kashiwara crystal
$\BB(V,{\bf B})=({\bf B},\tilde \gamma,\tilde
\varphi_i,\tilde \varepsilon_i,\tilde e_i|i\in I)$, where \begin{enumerate}
\item[$\bullet$] $\tilde \gamma: {\bf B}\to \Lambda^\vee$ is given
by $\tilde \gamma (b)=\mu$ whenever $b\in {\bf B}(\mu)$.
\item[$\bullet$] the functions $\tilde \varepsilon_i,\tilde
\varphi_i:{\bf B}\to \ZZ$ are given by $\tilde
\varepsilon_i=\ell_i$, $\tilde \varphi_i=\tilde
\varepsilon_i+\left<\alpha_i,\tilde \gamma(\bullet)\right>$ for
all $i\in I$.\end{enumerate}

\end{claim}

In order to establish the functoriality of the association
$(V,{\bf B})\mapsto \BB(V,{\bf B})$, we will define morphisms of
based $\BBB$-modules. Following \cite{bz-string}, for a locally
finite $\BBB$-module and any sequence $\ii\in I^m$, define a map
$e_\ii^{\rm top}:V\setminus\{0\}\to V\setminus\{0\}$   as follows. First, for $i\in I$ and $v\ne 0$, we set
$$e_i^{\rm top}(v):=e_i^{\ell_i(v)}(v) \ ,$$ and for $\ii=(i_1,\ldots,i_m)$, $m>1$, define
\begin{equation}
\label{eq:e top def}
e_{\ii}^{\rm top}=e_{i_m}^{\rm top}\circ e_{i_{m-1}}^{\rm top}\circ \cdots \circ
 e_{i_1}^{\rm top} \ .
\end{equation}

Denote by $V^+$ the space of the {\it highest weight vectors} of $V$:
\begin{equation}
\label{eq:V+}
V^+=\{v\in V:e_i(v)=0 ~\forall ~i\in I\} \ .
\end{equation}

Since $V$ is a locally finite module, for each $v\in V$, there
exists a sequence $\ii$ such that $e_\ii^{\rm top}(v)\in V^+$. In
particular, $V^+\ne \{0\}$ if $V\ne \{0\}$.

For each based $\BBB$-module $(V,{\bf B}),$  denote ${\bf
B}^+:={\bf B}\cap V^+$.

\begin{claim}
\label{cl:highest weight perfect basis} For any perfect basis
${\bf B}$ for $V$, the subset ${\bf B}^+$ is a basis  for $V^+$.
\end{claim}

Note that ${\bf B}^+=\BB(V,{\bf B})^+$, the set of highest weight elements in the sense of Definition \ref{def:highest weight element}.

\begin{definition}
\label{def:morphism of based modules} A {\it homomorphism of based
$\BBB$-modules} $(V,{\bf B})\to (V',{\bf B}')$ is any
$\BBB$-linear map $f:V\to V'$ such that $0\notin f({\bf B})$, and
there exist a map $\tilde f:{\bf B}\to {\bf B}'$ and a function
$c:{\bf B}\setminus {\bf B}^+\to \CC^\times$ ($b\mapsto c_b$)
satisfying $\tilde f(b)=f(b)$ for each $b\in {\bf B}^+$ and
\begin{equation}
\label{eq:general transition}
 f(b)-c_b\tilde f(b)\prec_\ii f(b)
\end{equation}
for each $b\in {\bf B}\setminus {\bf B}^+$ and any sequence $\ii$ such
that $e_\ii^{\,\rm top}(b)\in V^+$.

\end{definition}

For example, for each based $\BBB$-module $(V,{\bf B})$, the map
$f:V\oplus V\to V$ given by $f(v,v')=v+v'$ is a morphism of
$\BBB$-modules $(V\oplus V,{\bf B}\sqcup {\bf B})\to (V,{\bf B})$.

\begin{lemma}
\label{le:uniquenes of tilde f} The map $\tilde f$ and the
function $b\mapsto c_b$ in Definition \ref{def:morphism of based
modules} are unique. More precisely, for each $b\in {\bf
B}\setminus {\bf B}^+$ and $\ii$ such that $e_\ii^{\,\rm
top}(b)\in V^+$, the element $\tilde f(b)\in {\bf B}'$ and the
number $c_b$ satisfying \eqref{eq:general transition} are unique.

\end{lemma}

\begin{proof} We need some notation. Similarly to \eqref{eq:e top def}, for a based
$\BBB$-module
$(V,{\bf B})$ and any sequence $\ii=(i_1,\ldots,i_m)\in I^m$, $m\ge 1$ define the maps
 $\tilde e_\ii^{\rm \, top}:{\bf B}\to {\bf B}$ as follows.
First, define $\tilde e_i^{\rm \,top}:{\bf B}\to {\bf B}$, $i\in
I,$ by
$$\tilde e_i^{\,\rm top}(b):=\tilde e_i^{\,\ell_i(b)}(b) \ ,$$ and for
$\ii=(i_1,\ldots,i_m)$, define
$$\tilde e_{\ii}^{\,\rm top}=\tilde e_{i_m}^{\,\rm top}\circ
\tilde e_{i_{m-1}}^{\,\rm top}\circ \cdots \circ  \tilde e_{i_1}^{\,\rm top} \ .$$

\begin{claim}
\label{cl:tilde top} Let $(V,{\bf B})$ be a based $\BBB$-module.
Then  for any sequence $\ii\in I^m$, $m\ge 1$, we have:
\begin{enumerate}\item[(a)] $e_\ii^{\,\rm top}(b)\in \CC^\times\cdot \tilde
e_\ii^{\,\rm top}(b)$ for any $b\in {\bf B}$. \item[(b)] If
$e_\ii^{\,\rm top}(b)\in V^+$ for some $b\in {\bf B}$, then
$\tilde e_\ii^{\,\rm top}(b)\in {\bf B}^+$. \item[(c)] If $b
\equiv_\ii b'$ and $\tilde e_\ii^{\,\rm top}(b)=\tilde
e_\ii^{\,\rm top}(b')$, then  $b=b'$ (for any $b,b'\in {\bf B}$).
\end{enumerate}
\end{claim}

In the notation of Lemma \ref{le:uniquenes of tilde f},  let $b\in {\bf B}$,
$b',b''\in {\bf B}'\setminus{\bf B}^{'+}$ be such that $f(b)
-c'b'\prec_\ii f(b)$ and $f(b) -c''b''\prec_\ii f(b)$ for some and
$c',c''\in \CC^\times$ and some  $\ii$  such that $e_\ii^{\rm
top}(b)\in V^+$. We will show that $b'=b''$ and $c'=c''$. Indeed,
the above inequalities imply that $c'b'\equiv_\ii c''b''\equiv_\ii
f(b)$ and  $c'b'-c''b''\prec_\ii f(b)$. Therefore, $c'e_\ii^{\,\rm
top}(b')=c''e_\ii^{\,\rm top}(b'')$. This and  Claim
\ref{cl:tilde top}(a) imply that  $\tilde e_\ii^{\,\rm
top}(b')=\tilde e_\ii^{\,\rm top}(b'')$. Finally,
Claim~\ref{cl:tilde top}(c) implies that $b'=b''$. In  turn, this
implies that $c'=c''$.

The lemma is proved.
\end{proof}

Clearly, based $\BBB$-modules with such defined morphisms form a
category. This category possesses direct sums  $(V,{\bf B})\oplus
(V',{\bf B}')=(V\oplus V',{\bf B}\sqcup {\bf B}')$. However, the
category is not additive, because $\Hom((V,{\bf B}),(V',{\bf
B}'))$ is not an abelian group.

\begin{lemma}
\label{le:functor crystal}
The association $(V,{\bf B})\mapsto \BB(V,{\bf B})$ is a functor from the category of based $\BBB$-modules to the category of upper normal crystals.
More precisely, each morphism $f:(V,{\bf B})\to (V',{\bf B}')$ of based $\BBB$-modules defines the homomorphism
$\tilde f=\BB(f):\BB(V,{\bf B})\to \BB(V',{\bf B}')$ of upper normal crystals.

\end{lemma}

\begin{proof} By definition, a morphism of based modules $f:(V,{\bf B})\to (V',{\bf B}')$ defines a map
$\tilde f:{\bf B}\to {\bf B}'$ satisfying \eqref{eq:general
transition}. In particular, taking $\ii=(i_1,i_2,\ldots,i_m)$ in
\eqref{eq:general transition} with $i_1=i$ proves that
$\ell_i(\tilde f(b))=\ell_i(b)$ for all $b\in {\bf B}$, $i\in I$.
Therefore, in order to prove that $\tilde f$ is a homomorphism of
Kashiwara crystals, it suffices to show that $\tilde f$  is $\tilde
e_i$-equivariant for each $i\in I$. Let us choose any sequence
$\ii=(i_1,\ldots,i_m)\in I^m$ such that $i_1=i$ and $e_\ii^{\rm
top}(b)\in V^+$.

Applying $e_i$ to \eqref{eq:general transition} yields (provided that $\ell_i(b)\ne 0$):
$$e_i(\tilde f(b))\equiv_\ii e_i(f(b)) \ .$$
At the same time,  \eqref{eq:general upper asociated crystal} implies that
$$e_i(b)\equiv_\ii \tilde e_i(b),~ e_i(\tilde f(b))\equiv_\ii  \tilde e_i(\tilde f(b)) \ .$$
Therefore, taking into  account that
$e_i(f(b))=\varphi(e_i(b))\equiv_\ii f(\tilde e_i(b))$, we obtain:
$$\tilde e_i(\tilde f(b))\equiv_\ii f(\tilde e_i(b)) \ .$$
Therefore, $\tilde e_i(\tilde f(b))=\tilde f(\tilde e_i(b))$ because
according to \eqref{eq:general transition} and Lemma
\ref{le:uniquenes of tilde f}, $\tilde f(\tilde e_i(b))$ is the only
element of ${\bf B}'$ such that $\tilde f(\tilde e_i(b))\equiv_\ii
f(\tilde e_i(b))$. This proves the lemma.
\end{proof}

The following is the main result of this section.
\begin{maintheorem}
\label{th:the upper explicit crystal bijection}  Let ${\bf B}$ and
${\bf B}'$ be perfect bases of a locally finite $\BBB$-module such
that ${{\bf B}'}^+={\bf B}^+$. Then the identity map $\id:V\to V$
defines an isomorphism of based $\BBB$ modules $(V,{\bf B})\osr
(V,{\bf B}')$, that is, there is a unique map $\tilde f:{\bf B}\to
{\bf B}'$ and a unique function $c:{\bf B}\setminus {\bf B}^+\to
\CC^\times$ satisfying $\tilde f(b)=b$ for each $b\in {\bf
B}\setminus {\bf B}^+$, and
\begin{equation}
\label{eq:transition}
\tilde f(b)-c_b b\prec_\ii b
\end{equation}
for each $b\in {\bf B}\setminus {\bf B}^+$ and any sequence $\ii$ such that $\tilde e_\ii^{\,\rm top}(b)\in {\bf B}^+$.

In particular, one has   an isomorphism of (upper normal)
Kashiwara crystals $\BB(V,{\bf B})\osr \BB(V,{\bf B}')$.

\end{maintheorem}

\begin{proof} We need the following result.

\begin{lemma}
\label{le:top} Let $(V,{\bf B})$ be a based $\BBB$-module. Then
for any sequence $\ii\in I^m$, $m\ge 1,$ we have:
\begin{enumerate}
\item[(a)] Each $v\in V\setminus \{0\}$ is in the $\CC$-linear
span of $\{b\in {\bf B}:b\preceq_\ii v\}$. \item[(b)]  For each
$v\in V\setminus \{0\},$ there exists a unique  vector $v_0$ in
the linear span of $\{ b \in {\bf B}: b\equiv_\ii v\}$ such that
$v-v_0\prec_\ii v$.

\item[(c)] For each $v\in V\setminus \{0\}$ such that $ce_\ii^{\rm
top}(v)\in {\bf B}$ (for some $c\in \CC^\times$), there exists a
unique element $b\in {\bf B}$ such that $v\equiv_\ii b$. Moreover,
$v-c b \prec_\ii b$ for some $c\in \CC^\times$.
\end{enumerate}
\end{lemma}

\begin{proof} Prove (a) first. We will proceed by induction in the length $m$ of $\ii=(i_1,\ldots,i_m)$.
Let us write the expansion of $v\in V\setminus \{0\}$:
\begin{equation}
\label{eq:expansion}
v=\sum_{b\in {\bf B}} c_b b \ .
\end{equation}
Denote by ${\bf B}_0$ the set of all $b\in {\bf B}$ such that
$c_b\ne 0$;  and denote $n_1:=\max\limits_{b\in {\bf B}_0} \ell_{i_1}(b)$. We have to show that ${\bf B}_0\preceq v$.
Clearly, $\ell_{i_1}(v)\le n_1$. Next, we will show that $n_1=\ell_{i_1}(v)$.
Indeed, assume, by   contradiction, that $\ell_{i_1}(v)<n_1$. Applying
$e_{i_1}^{n_1}$ to \eqref{eq:expansion}, we obtain
$$0=\sum_{b\in {\bf B}_1} c_{b} e_{i_1}^{n_1}(b)\ ,$$
where ${\bf B}_1:=\{b\in {\bf B}_0:\ell_{i_1}(b)=n_1\}$. In other
words, taking into account Claim~\ref{cl:tilde top}(a), we obtain
a non-trivial linear dependence
$$0=\sum_{b'\in \tilde e_{i_1}^{\,n_1}({\bf B}_1)} c'_{b}b'\ .$$
This is in contradiction   to     $\tilde e_{i_1}^{\,n_1}({\bf
B}_1)\subset {\bf B}$ being  linearly independent. Therefore,
$n_1=\ell_{i_1}(v)$. This, in particular, proves that ${\bf B}_0\setminus {\bf B}_1\prec_\ii v$. In order
to show that ${\bf B}_1\preceq v$, let us
again apply $e_{i_1}^{n_1}$ to \eqref{eq:expansion}:
$$e_{i_1}^{n_1}(v)=\sum_{b\in {\bf B}_1} c_b e_{i_1}^{n_1}(b)\ .$$
Using Claim \ref{cl:tilde top}(a) again and the inductive assumption with $\ii'=(i_2,\ldots,i_m)$, we obtain:
$\tilde e_{i_1}^{\,n_1}({\bf B}_1)\preceq_{\ii'} e_{i_1}^{n_1}(v)$. Therefore,  ${\bf B}_1\preceq_{\ii} v$.
This proves  (a).

Prove (b) now. Denote by ${\bf B}'_0$ the set of all $b\in {\bf B}$
such that $c_b\ne 0$ and $b\equiv_\ii v$. Then we set
$v_0:=\sum\limits_{b\in {\bf B}'_0} c_b b$ in the notation of
\eqref{eq:expansion}. Clearly, $v-v_0$ is in the span of all $b\in
{\bf B}$ such that $b\prec_\ii v$. Using Claim \ref{cl:general equiv}(a), we see that an inequality $v-v'_0\prec_\ii v$ implies that $v_0-v'_0\prec_\ii v$, i.e.,  $v_0-v'_0$ is in the span of ${\bf B}\setminus {\bf B}'_0$. Therefore, if $v'_0$ is in the span of ${\bf B}'_0$, this implies that $v'_0=v_0$, which proves the uniqueness of $v_0$ in the linear span of ${\bf B}'_0$ with the desired property. Part (b) is proved.

To prove (c), note that, under the assumptions of (c), the above set
${\bf B}'_0$ consists of a single element $b$. Therefore, $v_0=c_bb$
and $v-v_0\prec_\ii v\equiv_\ii v_0\equiv _\ii b$.  Part (c) is also
proved.

Lemma \ref{le:top} is proved. \end{proof}

Let $b\in {\bf B}\setminus {\bf B}^+$. We now apply Lemma
\ref{le:top}(c) with $v=b$ relative to the perfect basis ${\bf
B}'$. Then for any sequence $\ii$  such that $\tilde e_\ii^{\,\rm
top}(b)\in {\bf B}^+={{\bf B}'}^+$, there exists a unique element
$b'\in {\bf B}'$  such that $b\equiv_\ii b'$ and
\begin{equation}
\label{eq:general b bpirme}
b=v+c b'
\end{equation}
for some $c\in \CC^\times$, where  either $v=0$ or $v\prec_\ii b$,
$v\prec_\ii b'$. If $v=0$ then we set $\tilde f(b)=b'$ and end the
proof here. So we assume that $v\ne 0$. Then, according to Lemma
\ref{le:top}(a), $v$ is in the span of $\{b_1\in {\bf B}:
b_1\preceq_\ii v\}$; hence $v$ is the   span of ${\bf
B}\setminus\{b\}$.  By the same reasoning, $v$ is in the span of
${\bf B}'\setminus\{b'\}$.

Let now $\ii'$ be any other sequence such that $\tilde
e_{\ii'}^{\,\rm top}(b)\in {\bf B}^+$, i.e., $c'e_{\ii'}^{\rm
top}(b)\in {\bf B}^+$ for some $c'\in \CC^\times$. Next, we will
show that $v\not \equiv_{\ii'} b$. Indeed, if $v \equiv_{\ii'} b,$
then we would have (again by Lemma \ref{le:top}(c)) that $v-c'b$
is in the span of  ${\bf B}\setminus\{b\}$. But this contradicts
the above observation that $v$ itself is in the span of  ${\bf
B}\setminus\{b\}$. Similarly, one shows that $v\not \equiv_{\ii'}
b'$. Finally, applying Claim \ref{cl:general equiv}(b) with $\ii'$
for $\ii$, $v'=cb'$, and $v+v'=b$, we obtain $v\prec_{\ii'} b$,
$v\prec_{\ii'} b'$, and $b\equiv_{\ii'} b'$. Thus, we proved that
$b\equiv_{\ii'} b'$ for any $\ii'$ such that $\tilde
e_{\ii'}^{\,\rm top}(b)\in {\bf B}^+$.

Therefore, the identity map $V\to V$ defines a homomorphism $(V,{\bf
B}) \to (V,{\bf B}')$ of based $\BBB$-modules, which is, obviously,
an isomorphism. Applying Lemma \ref{le:functor crystal}, we obtain
an isomorphism of the corresponding upper normal crystals.

Theorem \ref{th:the upper explicit crystal bijection} is proved.
\end{proof}

\begin{corollary} \label{cor:isomorphic normal crystals}
Let $V$ be any locally finite $\BBB$-module such that $\dim V^+=1$.
Then for any perfect bases ${\bf B}$ and ${\bf B}'$ for $V$ the
upper normal crystals $\BB(V,{\bf B})$ and $\BB(V,{\bf B}')$ are
isomorphic.

\end{corollary}

This  allows for the following definition.

\begin{definition}
\label{def: upper asociated crystal} For each locally finite
$\BBB$-module $V$ such that $V^+\cong \CC$ and such that $V$
admits a perfect basis, define the {\it associated crystal}
$\BB(V)$ to be isomorphic to each of the Kashiwara crystals
$\BB(V,{\bf B})$, where ${\bf B}$ is any perfect basis of $V$.
Furthermore, if $V=\oplus_k V_k$, where each $V_k$ is a locally
finite $\BBB$-module such that $V_k^+\cong \CC$ and $V_k$ admits a
perfect basis, then we define the {\it associated crystal}
$\BB(V)$ by  $\BB(V):=\bigsqcup_k \BB(V_k)$.

\end{definition}

\subsection{Perfect bases of $\bb^\vee$-modules and associated crystals}
\label{subsect:Perfect bases of b-modules} Let $G^\vee$ be the
Langlands dual group of a reductive algebraic group $G$, and
$\gg^\vee$ be the Lie algebra of $G^\vee$, and let
$\gg^\vee=\uu_-^\vee\oplus \hh^\vee \oplus \uu^\vee$ be the Cartan
decomposition. Denote by $\bb^\vee:=\hh^\vee\oplus \uu^\vee$ the
Borel sub-algebra of $\gg^\vee$.

Recall that  in the setup of Section \ref{subsect:Perfect bases,
upper normal crystals, and  crystal bases}, $I$ is the vertex set
of Dynkin diagram of $G$, $\Lambda=X^\star(T)$,
$\Lambda=X_\star(T)$, and $\langle\bullet,\bullet\rangle:
\Lambda\times \Lambda^\vee\to \ZZ$ is the canonical pairing.  By
definition, in the notation of Section \ref{subsect:Perfect bases,
upper normal crystals, and crystal bases}, we have a natural
surjective homomorphism $\hat \bb\twoheadrightarrow \bb^\vee$ via
$h\mapsto h$, $e_i\mapsto e_i$ for $i\in I$, where the latter $e_i$
are Chevalley generators of $\uu^\vee$. Therefore, each
$\bb^\vee$-module is a $\BBB$-module, and all definitions and
results of Section \ref{subsect:Perfect bases, upper normal
crystals, and  crystal bases} are valid for locally finite
$\bb^\vee$-modules.

In this section we will construct some important examples of based $\bb^\vee$-modules and the associated crystals.

Denote by  $T^\vee$, $B^\vee$, $U^\vee$  the Lie groups of
$\hh^\vee$, $\bb^\vee$, and $\uu^\vee$, respectively. By
definition, $B^\vee$ is a Borel subgroup in $G^\vee$ such that
$B^\vee=T^\vee \cdot U^\vee=U^\vee \cdot  T^\vee$.

 Let us define the right $B^\vee$-action on $U^\vee$  via the isomorphism  $U^\vee=T^\vee\setminus B^\vee$.
 This defines a $\BB^\vee$-action on $U$ and, therefore, structure of the left $B^\vee$-module on the  algebra $\CC[U^\vee]$.
This locally finite $B^\vee$-module is naturally a  $\bb^\vee$-module.

\begin{claim} The $\bb^\vee$-module  $\CC[U^\vee]$ is locally finite.

\end{claim}

For each $\lambda\in X_\star(T)^+$, denote by $V^\lambda$ the
restriction to $\bb^\vee$ of the finite-dimensional
$\gg^\vee$-module with the highest weight $\lambda$.

It is well-known (see e.g., \cite{lu}) that
\begin{equation}
\label{eq:Vlambda quotient}
V^\lambda=U(\uu^\vee)/I_{\lambda}\ ,
\end{equation}
where $I_\lambda=\sum U(\uu^\vee)\cdot e_i^{\langle \alpha_i,\lambda\rangle +1}$.

Therefore, taking the graded dual of the $\bb^\vee$-equivariant
quotient map $U(\uu^\vee)\twoheadrightarrow V^\lambda$, we obtain
an embedding of $\bb$-modules:
\begin{equation}
\label{eq:jlambda}
{\bf j}_\lambda:V^\lambda\hookrightarrow \CC[U^\vee] \ .
\end{equation}
for each $\lambda\in X_\star(T)^+$.

Denote by $\iota^*:\CC[U^\vee]\to \CC[U^\vee]$ the pullback of the
positive inverse $\iota:U^\vee\to U^\vee$ (defined in
\eqref{eq:iota}) above. By twisting the $\bb^\vee$-action by
$\iota^*$, we obtain another $\bb^\vee$-action on $\CC[U^\vee]$.
Denote
\begin{equation}
\label{eq:ei*}
e_i^*:=\iota^*\circ e_i\circ \iota^* \ .
\end{equation}

\begin{claim}
\label{cl:PRV} Each embedding ${\bf j}_\lambda$ satisfies
\begin{equation}
\label{eq:PRV}
{\bf j}_\lambda(V^{\lambda})=\{f\in \CC[U^\vee]: {e_i^*}^{\langle \alpha_i,\lambda\rangle +1}(f)=0~\forall ~i\in I\} \ .
\end{equation}
Therefore,
$${\bf j}_{\lambda}(V^{\lambda})\subset {\bf j}_{\lambda+\mu}(V^{\lambda+\mu}) $$
for any $\lambda,\mu\in X_\star(T)^+$; and  ${\bf
j}_{\lambda}(V^{\lambda})={\bf j}_{\lambda'}(V^{\lambda'})$ whenever
$\lambda-\lambda'\in Z(G)$.
\end{claim}



%
%

\begin{corollary} For each $\lambda\in X_\star(T)^+$, one has
$$\CC[U^\vee]=\lim\limits_{\lambda\in  X_\star(T)^+} V^{\lambda}=\bigcup\limits_{\lambda\in  X_\star(T)^+} {\bf j}_\lambda(V^{\lambda}) \ .$$

\end{corollary}

This result can be generalized to $\CC[U_P^\vee]$, where $P$ is a
standard parabolic subgroup of $G$, $U_P\subset U$ is the unipotent
radical, i.e., $U_P=U\cap \overline {w_P}U^-\overline {w_P}^{\,-1}$
(see e.g., Example \ref{ex:Uw}), and  $U_P^\vee\subset U$ is the
{\it dual} unipotent radical, i.e., $U_P^\vee=U^\vee\cap
\overline{w_P}\, (U^-)^\vee \, \overline{w_P}^{\,-1}$.

\begin{proposition} For each $\lambda\in X_\star(Z(L_P))\cap X_\star(T)^+$,
 the $\bb^\vee$-linear embedding \eqref{eq:jlambda}  factors through $V^{\lambda}\hookrightarrow \CC[U_P^\vee]$.
With respect to these embeddings,
$$\CC[U_P^\vee]=\lim\limits_{\lambda\in X_\star(Z(L_P))\cap X_\star(T)^+}
V^{\lambda}=\bigcup\limits_{\lambda\in X_\star(Z(L_P))\cap
X_\star(T)^+} {\bf j}_\lambda(V^{\lambda}) \ .$$

\end{proposition}

\begin{proof}
Taking $\lambda\in X_\star(Z(L_P))$, we see that $\langle
\alpha_i,\lambda\rangle=0$ for all $i\in J(P)$ (i.e.,  such that
$e_i\in \ll^\vee$). That is, according to \eqref{eq:CUP}, ${\bf
j}_\lambda(V^{\lambda})\subset \CC[U_P^\vee]$. Finally, taking
$\lambda\to \infty$ in $X_\star(Z(L_P))\cap X_\star(T)^+$, Claim
\ref{cl:PRV} implies that ${\bf j}_\lambda(V^{\lambda})\subset
\CC[U_P^\vee]$.
\end{proof}

Note that the dual unipotent radical $U^\vee_P$  is being acted on
by $B^\vee$. Therefore, the coordinate algebra $\CC[U_P^\vee]$ is
a locally finite $\bb^\vee$-module.

Denote by ${\bf B}^{\dual}$ the {\it dual canonical basis} for
$\CC[U^\vee]$ which is the specialization at $q=1$ of the dual
canonical basis of the quantized coordinate algebra
$\CC_q[U^\vee]$ (see, e.g.,  \cite{bz-cpl}).

\begin{claim} The basis ${\bf B}^{\dual}$ is a perfect basis of the locally finite $\bb^\vee$-module $\CC[U^\vee]$.
Moreover, for each $\lambda\in X_\star(T)^+$, the intersection
${\bf B}^\lambda={\bf B}^{\dual}\cap {\bf j}_\lambda(V^{\lambda})$
is a perfect basis in ${\bf j}_\lambda(V^{\lambda})$.

\end{claim}

Therefore, according to Corollary \ref{cor:isomorphic normal
crystals} and Definition \ref{def: upper asociated crystal}, each
$\bb^\vee$-module of the form $V^\lambda$ defines the associated
crystal $\BB(V^\lambda)$ and one has the crystal
$\BB(\CC[U^\vee])$  associated with the $\bb^\vee$-module
$\CC[U^\vee]$. These crystals are related  by a family of
injective homomorphisms of upper Kashiwara crystals
\begin{equation}
\label{eq:tilde jlambda}\tilde {\bf j}_\lambda:\BB(V^\lambda)\hookrightarrow \BB(\CC[U^\vee]) \ .
\end{equation}

Moreover,
$$\BB(\CC[U^\vee])=\lim\limits_{\lambda\in X_\star(T)^+} \BB(V^{\lambda})=\bigcup\limits_{\lambda\in X_\star(T)^+}
\tilde {\bf j}_\lambda(\BB(V^{\lambda}))\ .$$

\begin{proposition}
\label{pr:perfect basis parabolic} For each parabolic subgroup
$P^\vee\subset G^\vee$, the intersection ${\bf B}^{\dual}_P={\bf
B}^{\dual}\cap \CC[U_P^\vee]$ is a perfect basis in the locally
finite $\bb^\vee$-module $\CC[U_P^\vee]$.

\end{proposition}

\begin{proof} Denote by  $U_L^\vee$ the intersection $U^\vee\cap \overline{w_P}\, U^\vee\, \overline{w_P}^{\,-1}$
(see e.g., Example \ref{ex:Uw}). By definition, one has a $B^\vee$-equivariant isomorphism $U_P^\vee\cong U_L^\vee \backslash U^\vee$,
where the $U^\vee$-action on the quotient $U_L^\vee \backslash U^\vee$ is the right multiplication, and the $T^\vee$-action is the conjugation.
Therefore, one has a surjective map $U^\vee\to U_P^\vee$ and hence,  a $\bb^\vee$-equivariant embedding of algebras:
$\CC[U_P^\vee]\hookrightarrow \CC[U^\vee]$. In fact, the image of $\CC[U_P^\vee]$ is the algebra $\CC[U^\vee]^{U_L^\vee}$ of $U_L^\vee$-invariants.

\begin{lemma} The intersection of the dual canonical basis ${\bf B}^{\dual}$ for $\CC[U^\vee]$ with
$\CC[U^\vee]^{U_{L^\vee}}$ is a basis of
$\CC[U^\vee]^{U_{L^\vee}}$.

\end{lemma}

\begin{proof}   Note that
\begin{equation}
\label{eq:CUP}
\CC[U^\vee]^{U_{L^\vee}}=\{f\in \CC[U^\vee]:e_j^*(f)=0~\forall j\in J\}\ ,
\end{equation}
where $J=J(P^{op})=\{j\in I:U^\vee_j\subset U_{L^\vee}\}$ and
$e_i^*$ is defined above in \eqref{eq:ei*}. It is well-known that
$\iota^*({\bf B}^{\dual})={\bf B}^{\dual}$. This implies that the
basis ${\bf B}^{\dual}$ is perfect under the twisted
$\bb^\vee$-action. Replacing $B^\vee$ with $B_P^\vee=T^\vee \cdot
U_{L^\vee}$, we see that under the twisted $B_P^\vee$-action the
basis ${\bf B}^{\dual}$ is still a perfect basis of the locally
finite $B_P^\vee$-module $\CC[U^\vee]$ and that
$\CC[U^\vee]^{U_{L^\vee}}$ is the set of all highest weight
vectors in this module. Therefore, the assertion follows from
Claim \ref{cl:highest weight perfect basis}. The lemma is proved.
\end{proof}

Therefore, the basis ${\bf B}^{\dual}_P:={\bf B}^{\dual}\cap
\CC[U^\vee]^{U_{L^\vee}}$ is the desired perfect basis of the
$\bb^\vee$-module $\CC[U_P^\vee]$. Proposition \ref{pr:perfect
basis parabolic} is proved.
\end{proof}

Corollary \ref{cor:isomorphic normal crystals} and Definition
\ref{def: upper asociated crystal} imply  existence of  the {\it
crystal basis} $\BB(\CC[U^\vee_P])$ for $\CC[U^\vee_P]$. For each
$\lambda\in X_\star(Z(L_P))\cap X_\star(T)^+$, the homomorphism
\eqref{eq:tilde jlambda} factors through $\BB(\CC[U^\vee_P])$
\begin{equation}
\label{eq:tilde jlambda P}\tilde {\bf j}_\lambda:\BB(V^\lambda)\hookrightarrow \BB(\CC[U^\vee_P]) \ .
\end{equation}
and
\begin{equation}
\label{eq:limit crystals P}
\BB(\CC[U^\vee_P])=\bigcup\limits_{\lambda\in X_\star(Z(L_P))\cap X_\star(T)^+}
 \tilde {\bf j}_\lambda(\BB(V^{\lambda}))\ .
\end{equation}

In fact, the associated crystal $\BB(\CC[U^\vee_P])$ can be
constructed as a limit of a directed family
\begin{equation}
\label{eq:directed family}
\tilde {\bf f}_{\lambda,\mu}:\BB(V^\lambda)\to \BB(V^{\lambda+\mu})
\end{equation} the category of upper normal crystals. First, define
$\tilde {\bf j}_{\lambda,\mu}: \BB(V^\lambda)\to
\BB(V^\lambda)\times \BB(V^{\mu})$  by $\tilde {\bf
j}_{\lambda,\mu}:\tilde b\mapsto (\tilde b,\tilde b_\mu)$. According
to Lemma \ref{le:product trivial}(b), each  $\tilde {\bf
j}_{\lambda,\mu}$ is an injective homomorphism of upper normal
crystals $\BB(V^\lambda)\hookrightarrow \BB(V^{\mu})\times
\BB(V^\lambda)$. Moreover, the range of $\tilde {\bf
f}_{\lambda,\mu}$  belongs to the range of the embedding
$\BB(V^{\lambda+\mu})\hookrightarrow \BB(V^\lambda)\times
\BB(V^\mu)$. This defines an injective homomorphism
\eqref{eq:directed family} and, therefore, a directed family
$(\BB(V^\lambda),\tilde {\bf f}_{\lambda,\mu})$, $\lambda,\mu\in
X_\star(T)$ of upper normal crystals.

\begin{corollary}
\label{cor:directed family} For each standard parabolic subgroup
$P$ of $G$, the limit of the directed family
$(\BB(V^\lambda),\tilde {\bf f}_{\lambda,\mu})$ over
$\lambda,\mu\in X_\star(Z(L_P))\cap X_\star(T)^+$ is isomorphic to
the associated crystal $\BB(\CC[U_P^\vee])$.

\end{corollary}

\subsection{Perfect bases of $\gg^\vee$-modules and associated crystals}
\label{subsect:Perfect bases of gg-modules}
Let $\bb^\vee=\hh^\vee\oplus \uu^\vee$ and $\bb^\vee_-=\uu_-^\vee\oplus \hh^\vee$ be the opposite Borel sub-algebras in $\gg^\vee$.

Denote by $U(\gg^\vee)$ the universal enveloping algebra of $\gg^\vee$ and say that $\gg^\vee$-module $V$ is {\it locally finite} if, for
each $v\in V$, the $\gg^\vee$-module $U(\gg^\vee)(v)$ is
finite-dimensional and $V$ admits the weight decomposition
\eqref{eq:general weight decomposition}. Note that each locally
finite $\gg^\vee$-module is also locally finite as both the
$\bb_-^\vee$-module and the $\bb^\vee$-module.

Similarly to Section \ref{subsect:Perfect bases, upper normal
crystals, and  crystal bases}, for each locally finite
$\gg^\vee$-module $V$, we define:

\noindent $\bullet$
the functions $\ell_i:V\setminus\{0\}\to \ZZ_{\ge 0}$, $i\in I\sqcup -I$,

\noindent $\bullet$ the subspaces $V_i^{<\ell}:=\{v\in
V:\ell_i(v)<\ell\}$ for $i\in I\sqcup -I$, $l\ge 0$,

\noindent $\bullet$ the space $V^+=V^{\bb^\vee}$ of highest weight vectors.

Now we propose the following analogue of Definition \ref{def:general upper perfect basis}.

\begin{definition}
\label{def:upper perfect basis} We say that a weight basis ${\bf
B}$ in a locally finite $\gg^\vee$-module $V$ is $\gg^\vee$-{\it
perfect} (or simply {\it perfect}) if for each $i\in I\sqcup -I$, there is a partial
bijection $\tilde e_i:{\bf B}\to {\bf B}$ (i.e., a bijection of a
subset of ${\bf B}$ onto another subset of ${\bf B}$, see also
Definition \ref{def:partial bijection}  above) such that $\tilde
e_i(b)\in {\bf B}$ if and only if  $e_i(b)\ne 0$. And in the
latter case, one has
\begin{equation}
\label{eq:upper asociated crystal}
e_i(b)\in \CC^\times\cdot \tilde e_i(b)+ V_i^{<\ell_i(b)-1} \ .
\end{equation}

\end{definition}

That is, each  $\gg^\vee$-perfect basis  is  $\bb_-^\vee$-perfect and $\bb^\vee$-perfect at the same time.

\begin{lemma}
\label{le:submodule perfect basis} For each locally finite
$\gg^\vee$-module $V$, there exists a perfect basis. Moreover, if
$V'$ is a $\gg^\vee$-submodule of $V$, then there exists a perfect
basis ${\bf B}$ for $V$ such that $V'\cap {\bf B}$ is a (perfect)
basis of $V'$.

\end{lemma}

\begin{proof} One can think of $V$ as a specialization at $q=1$ of a locally finite
module $V_q$ over the quantized enveloping algebra
$U_q(\gg^\vee)$. G. Lusztig constructed in \cite{lub1,lub2} the {\it canonical basis} for each simple $V_q$ which defines  a basis ${\bf B}_q$ for each locally finite $V_q$. It follows from \cite[Theorem 7.5]{lub2}  that for each $b\in
{\bf B}_q$ and $i\in I\sqcup -I$ with $e_i(b)\ne 0$,
there exists a unique element $b' \in {\bf B}_q$ such that
\begin{equation}
\label{eq:normalizatiion}
e_i(b)\in [\ell_i(b)]_i\cdot b'+ V_{q,i}^{<\ell_i(b)-1} \ ,
\end{equation}
where $[n]_i$ is a $q$-analogue of the number $n$ (the formula \eqref{eq:normalizatiion} also holds for any {\it global crystal basis} of $V_q$ (see \cite[Proposition 5.3.1]{k})). 

Substituting $q=1$ into \eqref{eq:normalizatiion} we see that the weight basis ${\bf B}_q|_{q=1}$ also satisfies
\eqref{eq:upper asociated crystal} and, therefore, is a perfect basis of
$V$ (moreover, ${\bf B}_q|_{q=1}$ defines a basis in each isotypic
component of $V$).

Now let  $V'$ be a  $\gg^\vee$-submodule of $V$. Clearly, $V'$ is
also locally finite and there exists a complementary (locally
finite) submodule $V''$ of $V$  such that $V\cong V'\oplus V''$.
Then one can choose a perfect basis ${\bf B}'$ for $V'$,  a
perfect basis ${\bf B}''$ for $V''$, and let ${\bf B}={\bf
B}'\sqcup {\bf B}''$ be the desired perfect basis of $V$.

The lemma is proved.
\end{proof}

According to Claim \ref{cl:perfect asociated crystal}, a perfect
basis ${\bf B}$ for a $\gg^\vee$-module $V$ defines an upper
normal Kashiwara crystal $\BB(V,{\bf B})_+$ - when $V$ is
considered $\bb^\vee$-module and a lower normal Kashiwara crystal
$\BB(V,{\bf B})_-$ - when $V$ is considered $\bb^\vee_-$-module
(see Section~\ref{subsect:Normal Kashiwara crystals}).

\begin{lemma}
\label{le:perfect normal} The Kashiwara crystals $\BB(V,{\bf B})_+$
and $\BB(V,{\bf B})_-$ are canonically isomorphic and constitute a
normal Kashiwara crystal $\BB(V,{\bf B})$.

\end{lemma}

\begin{proof} The isomorphism between $\BB(V,{\bf B})_+$  and $\BB(V,{\bf B})_-$
follows from the fact that the partial bijections $\tilde e_{-i}$ and $\tilde e_i$ are inverse to each other.

The normality of $\BB(V,{\bf B})$ follows from  the fact that for
each  vector $v\in V(\mu)$ (e.g., for each $b\in {\bf B}(\mu)$)
and $i\in I$, one has
\begin{equation}
\label{eq:epsilon and phi}
\ell_i(b)=\ell_{-i}(b)-\left<\alpha_i,\mu\right> \ ,
\end{equation}
where $\left<\bullet,\bullet\right>$ is the pairing between roots and co-weights  of $G$.
\end{proof}

\begin{remark} Every $\gg^\vee$-perfect basis in $V$ is {\it good} in the sense of
Gelfand-Zelevinsky (\cite{gz}), that is,
it is compatible with each subspace $V_i^{<\ell}$ of $V$.

\end{remark}

\begin{remark} According to \cite[Lemma 4.3]{cr}, every {\it positive basis} ${\bf B}$ for $V$
(i.e., satisfying $e_i(b)\in \QQ_{\ge 0}\cdot {\bf B}$ for any $b\in {\bf B}$, $i\in I\sqcup -I$) is perfect.

\end{remark}

\begin{claim} For any perfect basis ${\bf B}$ for $V$, the subset ${\bf B}^+=V^+\cap {\bf B}$ is  a basis  for $V^+$. In particular,
\begin{equation}
\label{eq:V decomposition}
V=\bigoplus_{b\in {\bf B}^+} U(\gg^\vee)\cdot b \ .
\end{equation}
\end{claim}


The following is the main result of this section.

\begin{maintheorem}
\label{th:the crystal}  For any perfect bases ${\bf B}$ and ${\bf
B}'$ in a locally finite $\gg^\vee$-module $V$, the normal
crystals $\BB(V,{\bf B})$ and $\BB(V,{\bf B}')$ are isomorphic.

\end{maintheorem}

\begin{proof} We will construct the isomorphism explicitly.
First, we need the following fact.

\begin{claim}  Let  ${\bf B}^+$ and  ${\bf B}'^+$ be any weight bases of $V^+$. Then
\begin{enumerate}\item[(a)] There exists a weight-preserving bijection
$\varphi^+:{\bf B}^+\hookrightarrow {{\bf B}'}^+$. \item[(b)] Each
such weight-preserving bijection $\varphi^+:{\bf B}^+\osr {{\bf
B}'}^+$ uniquely extends to an isomorphism of the
$\gg^\vee$-modules $\varphi:V\osr V$.\end{enumerate}
\end{claim}

Indeed, for any $\gg^\vee$-module automorphism $\varphi:V\osr V$
the set $\varphi({\bf B})$ is a perfect basis of $\varphi(V)=V$.
Therefore, we may assume, without loss of generality, that
$\varphi$ is the identity and hence ${\bf B}^+= {{\bf B}'}^+$.
This and Theorem \ref{th:the upper explicit crystal bijection}
guarantee an  isomorphism of normal crystals $\BB(V,{\bf B})\osr
\BB(V,{\bf B}')$.

Theorem \ref{th:the crystal} is proved. \end{proof}

\begin{remark} In the case when  ${\bf B}$ is any perfect basis and ${\bf B}'={\bf B}_q|_{q=1}$ is the specialization of  a canonical basis (\cite{lub1,lub2}) or a  global crystal basis (\cite{k}),
the normalized perfect basis ${\bf B}_{norm}=\{c_b\cdot b\,|\,b\in
{\bf B}\}$ satisfies the specialization of \eqref{eq:normalizatiion} at $q=1$, where the  function
${\bf B}\to \CC^\times$: $b\mapsto c_b$ is defined in
\eqref{eq:transition} and $c_b:=1$ for $b\in {\bf B}^+$.

\end{remark}

\begin{remark}
\label{rem:perfect canonical} Theorem \ref{th:the crystal} implies
that the canonical basis, the dual canonical basis, the (upper and lower) global crystal
bases (see e.g., \cite{k}), and the semi-canonical basis
(\cite{L00}) in each irreducible $\gg$-module $V_\lambda$,
$\lambda\in X_\star(T)^+$ induce the same associated crystal
$\BB(V_\lambda)$.

\end{remark}


Theorem \ref{th:the crystal}  leads to the following definition.

\begin{definition} For each locally finite $\gg$-module $V,$ define the {\it associated crystal} $\BB(V)$ to be isomorphic to each of the
Kashiwara crystals of the form $\BB(V, {\bf B})$, where ${\bf B}$ is any perfect basis of $V$.

\end{definition}

Furthermore, Theorem \ref{th:the crystal} and Lemma
\ref{le:submodule perfect basis}  imply the following result.

\begin{corollary}

\label{cor:injective homomorphism of crystal bases} Any injective
homomorphism of locally finite $\gg^\vee$-modules $V'\hookrightarrow
V$ defines an injective homomorphism $\BB(V')\hookrightarrow \BB(V)$
of associated crystals.

\end{corollary}



\section{Tropicalization of geometric crystals and unipotent bicrystals}

\subsection{From decorated geometric crystals to normal Kashiwara crystals}
\label{subsect:from decorated geometric crystals to combinatorial
crystals} Recall from Claim \ref{cl:category positive varieties}
that ${\mathcal V}_+$ is the category of positive varieties and from
Claim \ref{cl:V++} that ${\mathcal V}_{++}$ is the category of
decorated positive varieties. Recall also from Section
\ref{subsect:Newton polytopes and tropicalization} that ${\bf
Set_0}$ is the category of marked sets.

\begin{claim}
\label{cl:trop on positive varieties} In the notation of Claim
\ref{cl:adjoint functor to positive tori} and Theorem
\ref{th:Trop}, the composition of functors $\Trop\circ \tau\circ
{\mathcal G}^*$ is a monoidal functor ${\mathcal V}_+\to {\bf
Set_0}$ such that the image of each positive variety $(X,\Theta)$
under this functor is isomorphic (in ${\bf Set_0}$) to a lattice of the same
dimension.
\end{claim}
The monoidal functor ${\mathcal V}_+\to {\bf Set_0}$ from Claim
\ref{cl:trop on positive varieties} is unique up to isomorphism
due to Claim \ref{cl:adjoint functor to positive tori}(b).
Throughout the end of the paper, we fix one of these isomorphic functors
and, similarly to Theorem \ref{th:Trop},  refer to it as the
tropicalization of positive varieties, and denote it by
$\Trop:{\mathcal V}_+\to {\bf Set_0}$.

%

In fact, due to Corollary \ref{cor:positive class of functions} of
\cite[Appendix]{bk}, the tropicalization is well-defined for any
non-zero function on each positive variety. Namely, given a
positive variety $(X,\Theta)$ and  two non-zero rational functions
$f,f':X\to \AA^1$, we say that $f$ and $f'$ are {\it positively
equivalent} if there is an isomorphism of positive varieties
$h:(X,\Theta)\osr (X,\Theta)$ such that $f'=f\circ h$. Proposition
\ref{pr:positive equivalence composition} implies that for
positively equivalent functions $f$ and $f'$, one has
$$\Trop(f')=\Trop(f)\circ \tilde h \ ,$$
where $\tilde h$ is that bijection $\Trop(X,\Theta)\to
\Trop(X,\Theta)$ which is the tropicalization of the positive
equivalence $h:(X,\Theta)\to \Trop(X,\Theta)$.
Therefore, the following version of Corollary \ref{cor:positive class of functions} holds.

\begin{claim} Let $(X,\Theta)$ be a positive variety, and $f:X\to \AA^1$ be a non-zero rational function.  Then
the isomorphism class of the function $\Trop(f):\Trop(X,\Theta)\to
\ZZ$ in ${\bf Set_0}$ depends only on the positive equivalence
class of $f$.

\end{claim}

Let $({\mathcal X},\Theta)=((X,\gamma, \varphi_i,\varepsilon_i,
e_i^\cdot| i\in I),\Theta)$ be a positive a geometric pre-crystal
(see Section \ref{subsect:Positive structures on geometric
crystals and unipotent bicrystals}).  We define the $5$-tuple
$\BB_\Theta=\Trop({\mathcal X},{\Theta}):=(\widetilde X,\tilde
\gamma,\tilde \varphi_i,\tilde \varepsilon_i,\tilde e_i| i\in I)$,
where \begin{enumerate} \item[(i)] $\widetilde X=\Trop(X,\Theta)$,
\item[(ii)] the map $\tilde \gamma: \widetilde X\to  X_\star(T)$
is given by
$$\tilde \gamma=\Trop(\gamma)\ ,$$
\item[(iii)] the functions $\tilde \varphi_i,\tilde
\varepsilon_i,\tilde \alpha_i:\widetilde {X}\to \ZZ$:
$$\tilde \varphi_i=-\Trop( \varphi_i), \tilde \varepsilon_i=-\Trop(\varepsilon_i),\tilde \alpha_i=\Trop(\alpha_i\circ \gamma)=\langle \alpha_i,\tilde \gamma(\cdot)\rangle \ ,$$
\item[(iv)] the $\ZZ$-action $\tilde e_i^\cdot:\ZZ\times\widetilde
X\to \widetilde X$ is the tropicalization of the
$(\Theta_{\GG_m}\times \Theta, \Theta$)-positive $\GG_m$-action
$\tilde  e_i^\cdot:\GG_m\times X\to X$. This $\ZZ$-action defines
a bijection $\tilde  e_i:\widetilde X\to \widetilde X$ via $\tilde
e_i=\tilde e_i^{\,1}$.
\end{enumerate}

\begin{claim}\cite[Theorem 2.11]{bk}
\label{cl:tropicalization geometric crystal} For any geometric
pre-crystal ${\mathcal X}$, the $5$-tuple
$\BB_\Theta=\Trop({\mathcal X},{\Theta})$ is a Kashiwara crystal
in the category ${\bf Set}_{\bf 0}$ of marked sets. This crystal
is torsion-free in the sense that each $\tilde e_i$ is a
bijection.

\end{claim}

Given a positive  geometric pre-crystal $({\mathcal X},\Theta)$,
clearly, the pair $({\mathcal X}^{op},\Theta)$ is also positive (see
Claim \ref{cl:op geometric}).

\begin{claim}
\label{cl:op tropicalized} For any positive geometric pre-crystal
$({\mathcal X},\Theta)$, one has (in the\linebreak notation of
Claim \ref{cl:op combinatorial}):
$$\BB_{\Theta^{op}}=(\BB_{\Theta})^{op}\ , $$
where $\BB_{\Theta^{op}}:=\Trop({\mathcal X}^{op},\Theta)$.

\end{claim}

The following result explains why our Definition \ref{def:product
combinatorial} differs from the original definition of \cite{k93} by
a permutation of factors.

\begin{lemma}
\label{le:trop of product} Let $({\mathcal X},\Theta_X)$ and
$({\mathcal Y},\Theta_Y)$ be positive  geometric pre-crystals. Then
one has a canonical isomorphism of torsion-free  Kashiwara crystals:
\begin{equation}
\label{eq:positive product geometric decorated} \Trop({\mathcal
X}\times {\mathcal Y},\Theta_X\times \Theta_Y)=\Trop({\mathcal
X},\Theta_X)\times \Trop({\mathcal Y},\Theta_Y) \ .
\end{equation}

\end{lemma}

\begin{proof} Indeed, applying $\Trop$ to Definition \ref{def:product geometric crystal},
we obtain for $\tilde x\in \tilde X$, $\tilde y\in \tilde Y$:
$$\tilde \gamma''(\tilde x,\tilde y)=\tilde \gamma(\tilde x)+\tilde \gamma'(\tilde y)\ ,$$
and for each $i\in \Supp {\mathcal X}\cup \Supp {\mathcal Y}$, the
functions $\tilde \varphi''_i,\tilde \varepsilon''_i:\tilde
X\times \tilde Y\to \ZZ$ given by
$$-\tilde  \varphi''_i(\tilde x,\tilde y)=
\min(-\tilde \varphi_i(\tilde x),-\tilde \varphi'_i(\tilde y)-\langle
\alpha_i,\tilde \gamma(\tilde x)\rangle
=-\max(\tilde \varphi_i(\tilde x),\tilde \varphi'_i(\tilde y)+\langle
 \alpha_i,\tilde \gamma(\tilde x)\rangle,$$
$$-\tilde \varepsilon''_i(\tilde x,\tilde y)=\min(-\tilde \varepsilon'_i(\tilde y),
-\tilde \varepsilon_i(\tilde x)+ \langle \alpha_i,\tilde
\gamma'(\tilde y)\rangle)= -\max(\tilde \varepsilon'_i(\tilde
y),\tilde \varepsilon_i(\tilde x)- \langle \alpha_i,\tilde
\gamma'(\tilde y)\rangle);$$ the bijection $\tilde
e_i^{\,n}:\tilde X\times \tilde Y\to \tilde X\times \tilde Y$,
$n\in \ZZ$ is given by the formula
$$\tilde e_i^{\,n}(x,y)=(\tilde e_i^{\,n_1}(\tilde x),\tilde e_i^{\,n_2}(\tilde y)) \ ,$$
where
\begin{alignat*}{2}
&n_1&&=\min(n-\tilde \varepsilon_i(\tilde x),-\tilde
\varphi'_i(\tilde y))- \min(-\tilde \varepsilon_i(\tilde
x),-\tilde \varphi'_i(\tilde y))\\
& &&=\max(\tilde \varepsilon_i(\tilde x),\tilde \varphi'_i(\tilde
y))-\max(\tilde \varepsilon_i(\tilde x)-n,\tilde \varphi'_i(\tilde
y))\ ,\\
 &n_2&&=\min(-\tilde \varepsilon_i(\tilde x),-\tilde
\varphi'_i(\tilde y))- \min(-\tilde \varepsilon_i(\tilde x),-
n-\tilde \varphi'_i(\tilde y))\\ & &&=
 \max(\tilde \varepsilon_i(\tilde x), n+\tilde \varphi'_i(\tilde y))-\max
 (\tilde \varepsilon_i(\tilde x),\tilde \varphi'_i(\tilde y)) \ .
 \end{alignat*}

The above equations agree with Definition \ref{def:product combinatorial}. The lemma is proved.
\end{proof}

Furthermore, recall from Section
\ref{subsect:Positive structures on geometric crystals and unipotent
bicrystals} that a triple  $(\XX,f,\Theta)$ is a
positive decorated geometric pre-crystal if $(\XX,f)$ is a decorated geometric pre-crystal,
$(\XX,\Theta)$ is a positive geometric pre-crystal, and the function $f$
is  $\Theta$-positive.

For each positive decorated geometric pre-crystal
$(\XX,f,\Theta)$, one  defines:
\begin{enumerate}
\item[(v)] the function $\tilde f:\widetilde X\to \ZZ$  by $\tilde
f=\Trop(f)$, \item[(vi)] the set $\tilde B_{\tilde f}\subset
\widetilde X$ by $\tilde B_{\tilde f}:=\{\tilde x\in \widetilde
X:\tilde f(\tilde x)\ge 0\}$.\end{enumerate}

Denote by  $\BB_{f,\Theta}$ the Kashiwara crystal obtained by
restricting the torsion-free Kashiwara crystal
$\BB_\Theta=\Trop(\XX,\Theta)$ to the subset $\tilde B_{\tilde
f}$, i.e.,
\begin{equation}
\label{eq:truncated Kashiwara crystal decorated}
\BB_{f,\Theta}=(\tilde B_{\tilde f},\tilde \gamma|_{\tilde B_{\tilde f}},
\tilde \varphi_i|_{\tilde B_{\tilde f}},\tilde \varepsilon_i|_{\tilde B_{\tilde f}},\tilde e_i|_{\tilde B_{\tilde f}}| i\in I) \ .
\end{equation}

\begin{proposition}
\label{pr:normal Kashiwara decorated} For any positive  decorated
geometric pre-crystal $(\XX,f,\Theta),$ the Kashiwara crystal
$\BB_{f,\Theta}$ is normal.
\end{proposition}

\begin{proof} Let us rewrite \eqref{eq:f of e(x)} as follows:
\label{cl:f of e}
\begin{equation}
\label{eq:f of e decorated}
f(e_i^c(x))=f_0(x)+
\frac{c}{\varphi_i(x)}+
\frac{c^{-1}}{\varepsilon_i(x)}
\end{equation}
for $x\in X$, $c\in \GG_m$, $i\in I$, where $f_0(x)=f(x)-\frac{1}{\varphi_i(x)}-\frac{1}{\varepsilon_i(x)}$.



Denote $\widetilde{X}:=\Trop(X,{\Theta})$ and let
$\displaystyle{\tilde \varphi_i,\tilde \varepsilon_i, \tilde
f,\tilde f_0:\tilde X\to \ZZ}$ be the $\Theta$-tropicalizations of
the functions ${\frac{1}{\varphi_i}}$, $
{\frac{1}{\varepsilon_i}}$, $f$, and $f_0$, respectively (if $f_0=
0$, then $\tilde f_0=+\infty$).  Note that the  function
$\GG_m\times X\to \AA^1$ given by $(c,x)\mapsto f(e_i^c(x))$ is
$\Theta_{\GG_m}\times \Theta$-positive. It follows from  Theorem
\ref{th:Trop} that the tropicalization of this positive function
is the function  $\ZZ\times \widetilde{X}\to \ZZ$ given by
$(n,\tilde b)\mapsto \tilde f(\tilde e_i^{\,n}(\tilde b))$.
Therefore, applying Corollary~\ref{cor:trop of sum}(b) to the
identity \eqref{eq:f of e decorated}, we obtain:
\begin{equation}
\label{eq:f of tilde e decorated}
\tilde f(\tilde e_i^{\,n}(\tilde b))=\min(\tilde f_0(\tilde b),n+\tilde \varphi_i(\tilde b),-n+\tilde \varepsilon_i(\tilde b))
\end{equation}
for $\tilde b\in \widetilde{X}$, $n\in \ZZ$, $i\in \Supp \XX$. In
particular, taking $n=0$, we obtain
$$\tilde f(\tilde b)=\min(\tilde f_0(\tilde b),\tilde \varphi_i(\tilde b),\tilde \varepsilon_i(\tilde b))\ .$$

Since $\tilde b\in \tilde B_{\tilde f}$ if and only if $\tilde
f(\tilde b)\ge 0$, the above identity implies that $\tilde b\in
\tilde B_{\tilde f}$ if and only if  $\tilde f_0(\tilde b)\ge 0$,
$\tilde \varphi_i(\tilde b)\ge 0$, $\tilde \varepsilon_i(\tilde
b)\ge 0$. Therefore, according to \eqref{eq:f of tilde e decorated},
$\tilde e_i^{\,n}(\tilde b)\in \tilde B_{\tilde f}$ for $\tilde b\in
\tilde B_{\tilde f}$ if and only if $-\tilde \varphi_i(\tilde b)\le
n\le \varepsilon_i(\tilde b)$. On the other hand, by definition
\eqref{eq:eplus minus i}, $\tilde e_i^{\,n}(\tilde b)\in \tilde
B_{\tilde f}$ for some $\tilde b\in \tilde B_{\tilde f}$ if and only
if $-\ell_{-i}(\tilde b)\le n\le \ell_i(\tilde b)$. This implies
that
$$\ell_i(\tilde b)=\tilde \varepsilon_i(\tilde b) ,\ell_{-i}(\tilde b)=\tilde \varphi_i(\tilde b) $$
for all $\tilde b\in \tilde B_{\tilde f}$.
That is, $\BB_{f,\Theta}$ is normal.
The proposition is proved. \end{proof}


As we argued in Section \ref{subsect:Positive structures on geometric crystals and
unipotent bicrystals}, based on Definition \ref{def:product decorated geometric pre-crystals}, one defines product of positive
decorated geometric pre-crystals by
$$(\XX,f,\Theta)\times ({\mathcal Y},f',\Theta'):=(\XX\times {\mathcal Y},f*f',
\Theta\times \Theta') \ ,$$
where $\Theta\times \Theta'$ is the
natural positive structure on the product $X\times Y$ and $f*f':X\times Y\to \AA^1$
 is given by $(f*f')(x,y)=f(x)+f'(y)$.

\begin{proposition}
\label{pr:product normal crystals decorated} For any  positive
decorated geometric pre-crystals $(\XX,f,\Theta)$ and $({\mathcal
Y},f',\Theta')$, one has
\begin{equation}
\label{eq:product normal crystals decorated}
\BB_{f*f',\Theta\times \Theta'}=\BB_{f,\Theta}\times \BB_{f',\Theta'} \ .
\end{equation}
\end{proposition}

\begin{proof}
Due to the  canonical isomorphism \eqref{eq:positive product
geometric decorated} of torsion-free Kashiwara crystals it
suffices to show that the subset of $\widetilde {X\times
Y}=\Trop(X\times Y,\Theta\times \Theta')=\widetilde {X}\times
\widetilde {Y}$ involved in the left hand side of
\eqref{eq:product normal crystals decorated} is equal to the set
involved in the right hand side. Indeed, by definition,
$$\tilde B_{\widetilde {f* f'}}=\{(\tilde x,\tilde y)\in \widetilde {X}\times \widetilde {Y}:\widetilde {f* f'}(\tilde x,\tilde y)\ge 0\} \ ,$$
where $\widetilde {f* f'}:\widetilde X\times \widetilde Y\to \ZZ$
is the $\Theta\times \Theta'$-tropicalization of $f*f'$. Since
$f:X\to \AA^1$, $f':Y\to \AA^1$, and $f*f':X\times Y\to \AA^1$ are
positive functions, Corollary \ref{cor:trop of sum}(d) guarantees
that their tropicalizations $\tilde f:\widetilde X\to \ZZ$,
$\tilde f':\widetilde Y^\to \ZZ$, and $\widetilde {f*
f'}:\widetilde X\times \widetilde Y\to \ZZ$,  respectively,
satisfy
$$\widetilde {f*f'}(\tilde x,\tilde y)=\min(\tilde f(\tilde x),\tilde{f'}(\tilde y))$$
for all $\tilde x\in \widetilde {X}$, $\tilde y\in \widetilde {Y^-}$.
Hence, $\widetilde{f*f'}(\tilde x,\tilde y)\ge 0$ if and only if $\tilde f(\tilde x)\ge 0$ and $\tilde {f'}(\tilde y)\ge 0$. Therefore,
$$\tilde B_{\widetilde {f*f'}}=\{(\tilde x,\tilde y): \tilde x\in \tilde B_{\tilde f},\tilde y\in \tilde B_{\tilde f'}\}
=\tilde B_{\tilde f}\times \tilde B_{\tilde f'} \ .$$
The proposition is proved.
\end{proof}

\subsection{From positive unipotent bicrystals to normal Kashiwara crystals}
\label{subsect:from positive unipotent bicrystals to combinatorial
crystals} Let  $({\bf X},{\bf p},\Theta)$ be a positive unipotent
bicrystal. As in      Section \ref{subsect:from decorated
geometric crystals to combinatorial crystals}  above, we denote by
$\BB_\Theta$ the torsion-free Kashiwara crystals associated with
the positive geometric crystal of the form $(\XX,\Theta)={\mathcal
F}({\bf X},{\bf p},\Theta)$ (in the notation of Lemma~\ref{le:from
positive unipotent to positive geometric}(a)).

Clearly, all the definitions and results of the above Section
\ref{subsect:from decorated geometric crystals to combinatorial
crystals} are valid for  positive decorated geometric crystals
obtained this way from positive unipotent bicrystals.

The following result provides the tropicalization of the morphism ${\bf f}_w$ from \eqref{eq:homomorphism geometric crystals}.

\begin{claim}
\label{cl:homomorphism free crystals} Given a positive unipotent
bicrystal $({\bf X},{\bf p},\Theta)$ of type $w$, the
tropicalization of the $(\Theta,\Theta_T\cdot
\Theta_w^-)$-positive morphism ${\bf f}={\bf p}|_{X^-}:X^-\to
TB^-_w$ is a homomorphism of torsion-free Kashiwara crystals
\begin{equation}
\label{eq:homomorphism free crystals}
\tilde {\bf f}_w:\BB_{\Theta}\to \BB_{\Theta_T\cdot \Theta_w^-}=
X_\star(T)\times\BB_{\Theta_w^-} \ ,
\end{equation}
where $X_\star(T)$ is the trivial Kashiwara crystal as in Example
 \ref{ex:trivial combinatorial crystal}.  The support of $\tilde {\bf f}_w$ is $|w|$.
\end{claim}

For a positive  $(U\times U,\chi^{st})$-linear unipotent bicrystal
$({\bf X},{\bf p},f,\Theta)$, we define the map $\widetilde
{hw}:\BB_{\Theta}\to \BB_{\Theta_T\cdot \Theta_w^-}=X_\star(T)$ to
be the composition of $\tilde {\bf f}_w$ with the projection to
the first factor $X_\star(T)\times\BB_{\Theta_w^-}\to X_\star(T)$.
This is a homomorphism of Kashiwara crystals (see Example
\ref{ex:projection product to lattice}).

Note that  according  to Claim \ref{cl:geometric highest weight morphism} and Claim \ref{cl:positive hw}, $\widetilde {hw}$ is the tropicalization
of the $(\Theta,\Theta_T)$-positive rational morphism $hw_X:{\mathcal F}({\bf X},{\bf p})\to T$ defined
in \eqref{eq:hw_X}.

For each  $\lambda\in X_\star(T)$, define
$\BB_{\Theta}^\lambda:=\widetilde {hw}^{-1}(\lambda)$. In
particular, the underlying set $\BB_{\Theta}^\lambda$ is $\{\tilde
b\in \tilde B: \widetilde {hw}(\tilde b)=\lambda\}$.

\begin{claim}
\label{cl:highest weight decomposition crystal free} For a
positive unipotent bicrystal $({\bf X},{\bf p},\Theta)$ of type
$w,$ one has: \begin{enumerate}\item[(a)] Each
$\BB_{\Theta}^\lambda$ is invariant under the action of $\tilde
e_i^{\,n}$, $i\in I$, $n\in \ZZ$ and, therefore, is a
(torsion-free or empty) sub-crystal of the torsion-free Kashiwara
crystal $\BB_{\Theta}$. \item[(b)]  The direct decomposition
$\BB_{\Theta}=\bigsqcup_\lambda \BB_{\Theta}^\lambda$. \item[(c)]
$(\BB_{\Theta}^\lambda)^{op}=\BB_{\Theta^{op}}^{-w\lambda}$ for
each $\lambda\in X_\star(T)$, where
$\BB_{\Theta^{op}}=\Trop({\mathcal F}({\bf X},{\bf
p},\Theta)^{op})$.\end{enumerate}
\end{claim}


Now let  $({\bf X},{\bf p},f,\Theta)$ be a positive  $(U\times
U,\chi^{st})$-linear bicrystal. Denote by $\BB_{f,\Theta}$ the
tropicalization of the positive decorated geometric crystal
${\mathcal F}({\bf X},{\bf p},f,\Theta)$ (see Lemma \ref{le:from
positive unipotent to positive geometric} and \eqref{eq:truncated
Kashiwara crystal decorated}). According to Proposition
\ref{pr:normal Kashiwara decorated}, $\BB_{f,\Theta}$ is normal.
Furthermore, for each $\lambda\in X_\star(T)$, denote
\begin{equation}
\label{eq:Blambda}
\BB_{f,\Theta}^\lambda:=\BB_{f,\Theta}\cap \BB_\Theta^\lambda \ .
\end{equation}

\begin{claim}
\label{cl:highest weight decomposition crystal} For a positive
$(U\times U,\chi^{st})$-linear bicrystal $({\bf X}, {\bf
p},f,\Theta)$ of type $w$, one has: \begin{enumerate}\item[(a)]
Each $\BB_{f,\Theta}^\lambda$ is invariant under the action of
$\tilde e_i^{\,n}$, $i\in I$, $n\in \ZZ$ and, therefore, is a
normal sub-crystal of  $\BB_{f,\Theta}$. \item[(b)] The
direct decomposition $\BB_{f,\Theta}=\bigsqcup_\lambda
\BB_{f,\Theta}^\lambda$. \item[(c)]
$(\BB_{f,\Theta}^\lambda)^{op}=\BB_{f,\Theta^{op}}^{-w\lambda}$
for $\lambda\in X_\star(T)$, where
$\BB_{f,\Theta^{op}}=\Trop({\mathcal F}({\bf X},{\bf
p},f,\Theta)^{op})$.
\end{enumerate}

\end{claim}

\begin{remark} Each $\BB_{f,\Theta}^\lambda$ can be
thought of as an isotypic component of $\BB_{f,\Theta}$. Later
on, in Theorem \ref{th:Blambda}, we will make this analogy precise.
\end{remark}

\subsection{From strongly positive unipotent bicrystals to crystals associated with modules}
\label{subsect:From strongly positive unipotent bicrystals to
crystal bases} Given a strongly positive parabolic $(U\times
U,\chi^{st})$-linear bicrystal $({\bf X},{\bf p},f,\Theta)$ of
type $w_P$, we denote by $\tilde \Delta_X:\widetilde {X^-}\to \ZZ$
the $\Theta$-tropicalization of the positive function
$\Delta_X|_{X^-}:X^-\to \AA^1$. For each $\lambda\in
X_\star(Z(L))$ and $n\in \ZZ$, we denote
$$\BB_{f,\Theta;n}^\lambda:=\{\tilde b\in \BB_{f,\Theta}^\lambda:\tilde \Delta_X(\tilde b)=n\} \ .$$

\begin{claim}
\label{cl:tilde delta} The function $\tilde \Delta_X$ is invariant
under all crystal operators $\tilde e_i^{\,n}$. In particular,
each non-empty $\BB_{f,\Theta;n}^\lambda$ is a normal sub-crystal
of the normal crystal~$\BB_{f,\Theta}^\lambda$.
\end{claim}

Due to its $\tilde e_i^{\,n}$-invariance, will refer to the function $\tilde \Delta_X$ as the {\it combinatorial central charge}.

Recall that  $(X_P,\id,f_P,\Theta_P^-)$ is the strongly positive
$(U\times U,\chi^{st})$-linear bicrystal defined in Section
\ref{subsect:toric charts and positive structures}.

\begin{proposition}
\label{pr:homomorphism normal crystals} Let
$({\bf X},{\bf p},f,\Theta)$ be a strongly positive parabolic $(U\times U,\chi^{st})$-linear bicrystal of type $w_P$.
Then
$$\tilde {\bf f}_{w_P}(\BB_{f,\Theta})\subset\BB_{f_P,\Theta_P^-}\ , $$
where $\tilde {\bf f}_{w_P}:\BB_{\Theta}\to \BB_{\Theta_P^-}$ is given by \eqref{eq:homomorphism free crystals}. This defines
a homomorphism of normal crystals
\begin{equation}
\label{eq:homomorphism normal crystals}
\tilde {\bf m}:\BB_{f,\Theta}\to \BB_{f_P,\Theta_P^-}
\end{equation}
such that $\tilde  {\bf m}(\BB_{f,\Theta}^\lambda)\subset \BB_{f,\Theta_P^-}^\lambda$ for each $\lambda\in X_\star(Z(L))$.
\end{proposition}
\begin{proof} We need the following simple fact.

\begin{claim}
\label{cl:strongly positive to standard} In the notation of
Proposition \ref{pr:homomorphism normal crystals}, let $\tilde
f:\widetilde {X^-}\to \ZZ$, $\tilde f_{P}:\widetilde {X_P^-}\to
\ZZ$, and $\tilde \Delta_X:\widetilde {X^-}\to \ZZ$  be the
tropicalizations of positive functions $f|_{X^-}:X^-\to \AA^1$,
$f_P|_{X_P^-}:X_P^-\to \AA^1$, and  $\Delta_X|_{X^-}:X^-\to
\AA^1$, respectively. Then
$$\tilde f=\min(\tilde f_P\circ \tilde {\bf f}_{w_P},\tilde \Delta_X)\ .$$
In particular,
\begin{equation}
\label{eq:inequality f, f_P}
\tilde f(\tilde x)\le \tilde f_P ( \tilde {\bf f}_{w_P}(\tilde x))
\end{equation}
for any $\tilde x\in \widetilde {X^-}$.
\end{claim}

The  inequality \eqref{eq:inequality f, f_P} guarantees that $\tilde
{\bf f}_{w_P}(\tilde B_{\tilde f})\subset \tilde B_{\tilde f_P}$,
i.e., the  restriction of $\tilde {\bf f}_{w_P}$  to $\tilde
B_{\tilde f}$ defines a homomorphism of  normal crystals
$\BB_{f,\Theta}\to \BB_{f_P,\Theta^-_P}$.

The proposition is proved. \end{proof}


Recall that the associated crystal  $\BB(V)$ of each locally
finite $\gg^\vee$-module $V$  is defined in  Section
\ref{subsect:Perfect bases of gg-modules} above. In particular,
each $\BB(V)$ is the union of $\BB(V_\lambda)$, where  $V_\lambda$
is an irreducible finite-dimensional $\gg^\vee$-module with the
highest (co-)weight $\lambda\in X_\star(T)^+$, where
$X_\star(T)^+$ is the monoid of all dominant co-weights defined in
Section \ref{subsect:notation}.


Recall also from Claim
\ref{cl:highest weight decomposition crystal} that each Kashiwara
crystal of the form $\BB_{f,\Theta}$ is the disjoint union of
$\BB_{f,\Theta}^\lambda$, $\lambda \in X_\star(T)$.

The following is our first main result on the Kashiwara crystals
coming from strongly positive unipotent bicrystals.

\begin{maintheorem}
\label{th:Blambda} For each  $\lambda\in X_\star(T),$ we have:
$\BB_{f_B,\Theta_B^-}^{\lambda}$ is empty if $\lambda\notin
X_\star(T)^+$ and $\BB_{f_B,\Theta_B^-}^{\lambda}$ is isomorphic
to the associated crystal $\BB(V_\lambda)$ if $\lambda\in
X_\star(T)^+$.

\end{maintheorem}

\begin{proof} The proof is based on  Joseph's characterization of the associated crystals $\BB(V_\lambda)$. Following
\cite[Section 6.4.21]{Joseph}, we say that a family ${\mathcal
C}_\cdot=\{ {\mathcal C}_\lambda |\lambda\in  X_\star(T)^+\} $ of
normal crystals is {\it closed} if

\noindent $\bullet$ For each  $\lambda \in  X_\star(T)^+$, there
is a unique highest weight element $c_\lambda\in {\mathcal
C}_\lambda$ (see Definition \ref{def:highest weight element}) such
that $\tilde \gamma(c_\lambda)= \lambda$.

\noindent $\bullet$ For each $\lambda,\mu\in  X_\star(T)^+,$ the
correspondence $c_{\lambda+\mu}\mapsto (c_\lambda,c_\mu)\in
{\mathcal C}_\lambda\times {\mathcal C}_\mu$ defines an injective
homomorphism of normal crystals ${\mathcal
C}_{\lambda+\mu}\hookrightarrow {\mathcal C}_\lambda\times
{\mathcal C}_\mu$.

\begin{theorem} [\cite{Joseph}]
\label{th:closed family} If ${\mathcal C}_\cdot=\{ {\mathcal
C}_\lambda | \lambda\in X_\star(T)^+\}$ is a closed family of
crystals, then each ${\mathcal C}_\lambda$ is isomorphic  to the
associated crystal $\BB(V_\lambda)$.

\end{theorem}

We will show that the crystals $\BB_{f_B,\Theta_B^-}^{\lambda}$,
$\lambda\in X_\star(T)^+$ form a closed family. First, we  prove
that each $\BB_{\Theta_B^-}^{\lambda}$ has a unique highest weight
element. This is the most technical part of the entire proof
(Theorem \ref{th:Blambda highest}). After this we will use a
relatively simple argument based on the strong positivity of
$\Theta_B^-$ and $\Theta_{B^-}*\Theta_{B^-}$ to construct the embeddings
$\BB_{\Theta_B^-}^{\lambda+\mu}\hookrightarrow \to \BB_{\Theta_B^-}^{\lambda}\times \BB_{\Theta_B^-}^\mu$ and thus finish the proof of that
crystals $\BB_{\Theta_B^-}^{\lambda}$ form a closed family.

In fact, we will prove the existence and uniqueness of the highest
weight elements in a more general situation. Recall that
$X_P=UZ(L_P)\overline {w_P}U$, where $L_P$ is the Levi factor of $P$
and $Z(L_P)\subset T$ is the center of $L_P$. In order to formulate
the following result,   we identify $X_P^-=Z(L_P)B^-_{w_P}$ with
$Z(L_P)\times B^-_{w_P}$ and, passing to the tropicalization, we
identify $\widetilde {Z(L_P)B^-_{w_P}}$ with $X_\star(Z(L_P))\times
\widetilde {B^-_{w_P}}$.

\begin{theorem}
\label{th:Blambda highest} Let $P$ be a  standard parabolic
subgroup of $G$. Then for each  $\lambda\in X_\star(Z(L_P))$, we
have: if $\lambda$ is not dominant, then
$\BB_{f_P,\Theta_P^-}^\lambda=\emptyset$; and if $\lambda$ is
dominant, then the normal crystal $\BB_{f_P,\Theta_P^-}^\lambda$
has a unique highest weight element $\tilde b_\lambda$ of the
weight $\lambda$ or, more precisely, under the identification
$\widetilde {Z(L_P) B^-_{w_P}}=X_\star(Z(L_P))\times \widetilde
{B^-_{w_P}}$, one has
$$\tilde b_\lambda=(\lambda,{\bf 0}) \ ,$$
where ${\bf 0}\in \widetilde {B^-_{w_P}}$ is the marked point (see
Sections \ref{subsect:Newton polytopes and tropicalization} and
\ref{subsect:from decorated geometric crystals to combinatorial
crystals}).

\end{theorem}
\begin{proof}
In the notation of Theorem \ref{th:Blambda highest}, denote by
$\tilde f_P:X_\star(Z(L_P))\times \widetilde {B^-_{w_P}}\to \ZZ$ be
the $\Theta_P^-$-tropicalization of the positive function $f_P:
Z(L_P)\times B^-_{w_P}\to \AA^1$, and by $\tilde
\varepsilon_i:X_\star(Z(L_P))\times \ZZ^{l(w_P)}\to \ZZ$ -- the
$\Theta_P^-$-tropicalization of the positive function
$\varepsilon_i: Z(L_P)\times B^-_{w_P}\to \AA^1$ for $i\in I$.
Therefore, it suffices to prove the following result.

\begin{proposition}
\label{pr:lambda b, b=0}
Let $(\lambda,\tilde b) \in X_\star(Z(L_P))\times \widetilde {B^-_{w_P}}$
 be any point such that
$$\tilde f_P(\lambda,\tilde b)\ge 0, ~\tilde \varepsilon_i(\lambda,\tilde b)\le 0$$
for all $i\in I$. Then $\tilde b={\bf 0}$.
\end{proposition}

\begin{proof} For a character $\chi:U\to \AA^1$ define a regular function
 $f_{w,\chi}:U\overline w U\to \AA^1$ by
\begin{equation}
\label{eq:fwst}
f_{w,\chi}(g)=\chi(\pi^+(\overline w^{\,-1}g))
\end{equation}
for $g\in BwB$.

It follows from Lemma \ref{le:formula fL} that
$$f_P(g)=f_{w,\chi^{st}}(g)+\sum_{j\in J(P)} h_j(g)$$
for any $g\in X_P=UZ(L_P)\overline{w_P}U$, where
$h_j(g)=\chi_j(\pi^+(\overline {w_P^{\,-1}}^{\,-1}\iota(g)))$. Using
the fact that $f_{w,\chi^{st}}(tg)=f_{w,\chi^{st}}(g)$  and
$h_j(tg)= \alpha_j(t)f_j(g)$ for $g\in Bw_PB$, $t\in T$, $j\in
J(P)$,  we obtain
\begin{equation}
\label{eq:fp modified}
f_P(t\cdot b)=f_{w,\chi^{st}}(b)+\sum_{j\in J(P)}\alpha_j(t)f_j(b)
\end{equation}
for $t\in Z(L_P)$, $b\in B_{w_P}^-=U\overline{w_P}U\cap B^-$.

Therefore, applying the $\Theta_P^-$-tropicalization to \eqref{eq:fp
modified} and using Corollary \ref{cor:trop of sum}(c) with $T'=T$,
$\mu_j=\alpha_j$, we obtain:
\begin{equation}
\label{eq:tropical f_P} \tilde f_P(\lambda, \tilde b)=\min(\tilde
f_{w,\chi^{st}}(\tilde b), \min_{j\in J(P)}(
\langle\alpha_j,\lambda\rangle+\tilde h_j(\tilde b)))\le \tilde
f_{w,\chi^{st}}(\tilde b)
\end{equation}for all $\tilde b\in \widetilde {B^-_{w_P}}$,
where $\tilde f_{w,\chi^{st}}:\widetilde {B^-_{w_P}}\to \ZZ$ is the
$\Theta_{w_P}^-$-tropicalization of $f_{w,\chi^{st}}|_{B^-_{w_P}}$,
and  $\tilde f_j:X_\star(Z(L_P))\times \widetilde {B^-_{w_P}}\to
\ZZ$ is the $\Theta_P^-$-tropicalization of $f_j|_{X_P^-}$.

 In particular,
$$\{\tilde b\in \widetilde {B^-_{w_P}}:\tilde f_P((\lambda,\tilde b)\ge 0\}\subset \{\tilde b \in
\widetilde {B^-_{w_P}}: \tilde f_{w,\chi^{st}}(\tilde b)\ge 0 \}\ .$$

Therefore, taking into  account that  $\varepsilon_i(t\cdot
b)=\varepsilon_i(b)$ for any $b\in B^-_w$, $t\in T$ and hence
$\tilde \varepsilon_i(\lambda,\tilde b)=\tilde
\varepsilon_i(\tilde b)$ for any $\tilde b\in \widetilde
{B^-_{w_P}}$, $\lambda\in X_\star(T)$, Proposition \ref{pr:lambda
b, b=0} follows from the following result.

\begin{lemma}
\label{le:psiii tilde} For any $w\in W$, the tropicalizations
$\widetilde \varepsilon_i:\widetilde {B^-_w}\to \ZZ$ and
$\widetilde f_{w,\chi^{st}}:\widetilde {B^-_w}\to \ZZ$ of
$\Theta_w^-$-positive functions $\varepsilon_i:B^-_w\to \AA^1$ and
$f_{w,\chi^{st}}:B^-_w\to \ZZ$ satisfy
\begin{equation}
\label{eq:positive psi}
\{\tilde b\in \widetilde {B^-_w}:\tilde f_{w,\chi^{st}}(\tilde b)
\ge 0, \tilde \varepsilon_i(\tilde b)\le 0~\forall i\in I\}=\{{\bf 0}\} \ .
\end{equation}

\end{lemma}

\begin{proof} We need the following  recursive formula for $f_{w,\chi}$.

\begin{lemma}
\label{le:f_w factorization} For any $w_1,w_2\in W$ such that
$l(w_1w_2)=l(w_1)+l(w_2)$, one has
$$f_{w_1w_2,\chi}(g_1g_2)=f_{w_2,\chi}(u_1g_2)$$
for any $g_1\in Bw_1 B$, $g_2\in Bw_2 B$, where $u_1=\pi^+(\overline {w_1}^{\,-1}g_1)$.

\end{lemma}

\begin{proof} Indeed,
$$\pi^+(\overline{w_1w_2}^{\,-1}g_1g_2)=\pi^+(\overline{w_2}^{\,-1}\overline{w_1}^{\,-1}g_1g_2)
=\pi^+(\overline{w_2}^{\,-1}b_1
u_1g_2)=\pi^+(\overline{w_2}^{\,-1}b_1 u_1g_2)\ ,$$ where
$u_1=\pi^+(\overline {w_1}^{\,-1}g_1)$ and $b_1=\pi^-(\overline
{w_1}^{\,-1}g_1)\in T\cdot (U^-\cap
\overline{w_1}^{\,-1}U\overline{w_1})$ (we used the fact that
$\overline{w_2}^{\,-1}T\cdot (U^-\cap
\overline{w_1}^{\,-1}U\overline{w_1})\overline{w_2}\subset T\cdot
(U^-\cap \overline{w_1w_2}^{\,-1}U\overline{w_1w_2})\subset B^-$).
Therefore,
$$f_{w_1w_2,\chi}(g_1g_2)=\chi(\pi^+(\overline{w_1w_2}^{\,-1}g_1g_2))=\chi(\pi^+(\overline{w_2}^{\,-1}b_1 u_1g_2))=f_{w_2,\chi}(u_1g_2) \ .$$
The lemma is proved.
\end{proof}

\begin{lemma}
\label{le:positive i-decomposition f right}
Let  $j\in I$ be such that  $w=s_jw'$ and $l(w')=l(w)-1$. Then
\begin{equation}
\label{eq:positive i-decomposition f right}
f_{s_jw',\chi^{st}}(x_{-j}(c)\cdot b)=f_{w',\chi^{st}}(x') +\sum\limits_{k\ge 1} c^k f_k(b)
\end{equation}
for any  $b\in B^-_{w'}$, $c\in \GG_m$, where each $f_k:B^-_{w'}\to \AA^1$ is a regular function and $f_k\not = 0$ for at least one $k>1$.

\end{lemma}

\begin{proof} We need the following simple fact.

\begin{claim}
\label{cl:u-shift}
Let $Y$ be a variety equipped with a  left $U$-action  (which we denote by $(u,y)\mapsto u\cdot y$) and let $f:Y\to \AA^1$ be a regular function. Then
\begin{equation}
\label{eq:u-shift}
f(u\cdot y)=f(y)+\sum_{k\ge 1} h_k(u)f_k(y)
\end{equation}
for all $u\in U$, $y\in Y$, where the sum is finite, each $f_k$ is a
regular function on $Y$, and each $h_k$ is a regular function on $U$
such that $h_k(e)=0$.

\end{claim}

Now, using Lemma \ref{le:f_w factorization} with $w_1=s_j$,
$w_2=w'$, $g_1=x_{-j}(c)$, $g_2=b$, and taking into account that $\pi^+(\overline
{s_j}^{\,-1}x_{-j}(c))=x_j(c)$, the formula \eqref{eq:u-shift}
with $f=f_{w',\chi^{st}}$ becomes
$$f_{w,\chi^{st}}(x_{-j}(c)\cdot b)=f_{w',\chi^{st}}(x_j(c)\cdot b)=f_{w',\chi^{st}}(b) +\sum\limits_{k\ge 1} h_k(x_j(c)) f_k(b) \ .$$
Without loss of generality, we can assume that $h_k(x_j(c))=c^k$
for $k\ge 1$.

Finally, it is clear that  the function $(c,b)\mapsto f_{w',\chi^{st}}(x_j(c)\cdot b)$ depends on $c$, therefore, $f_k\not = 0$ for some $k$.
The lemma is proved.
\end{proof}

We continue the proof of Lemma \ref{le:psiii tilde} by induction
on the length $l(w)$. Indeed, if $l(w)=0$, i.e., $w=e$, we have
nothing to prove. Now let us write $w=s_jw'$ so that
$l(w')=l(w)-1$. Then, according to Claim \ref{cl:product standard
positive structures},   the positive structure $\Theta_w^-$
factors as $\Theta_w^-= \Theta_{s_j}^-\times \Theta^-_{w'}$.
Therefore, we can identify the tropicalizations  $\widetilde
{B^-_w}= \ZZ\times \widetilde {B^-_{w'}}$.

Then, tropicalizing  \eqref{eq:positive
i-decomposition f right}, we obtain (using Corollary \ref{cor:trop
of sum}(b))
\begin{equation}
\label{eq:positive i-decomposition f right tropical}
\tilde f_{w,\chi^{st}}(m,\tilde b)=\min(\tilde f_{w',\chi^{st}}(\tilde b),\min_{k>1}
  (km+ \tilde f_k(\tilde b)))
\end{equation}
for $(m,\tilde b)\in \ZZ\times \widetilde {B^-_{w'}}$, where
$\tilde f_k$ is the tropicalization of $f_k$ (e.g., $\tilde f_k=
+\infty$ whenever $f_k= 0$). In particular, the inequality $\tilde
f_{w,\chi^{st}}(m,\tilde b)\ge 0$ implies
\begin{equation}
\label{eq:positive i-inequality f right tropical}
 \tilde f_{w',\chi^{st}}(\tilde b)\ge 0, km \ge -\tilde f_k(\tilde b),
\end{equation}
for each $k>0$ such that $f_k\not = 0$ (such $k$ always exists by
 Lemma \ref{le:positive i-decomposition f right}).

On the other hand, $\widetilde {B^-_w}= \ZZ\times \widetilde
{B^-_{w'}}$ is the product of  torsion-free Kashiwara crystals. Then
we have by Definition \ref{def:product combinatorial}:
$$\tilde \varepsilon_i(m,\tilde b)= \max(\tilde \varepsilon'_i(\tilde b),
\tilde \varepsilon_i(m)-\langle \alpha_i,\tilde \gamma'(\tilde b)\rangle )
=\max(\tilde \varepsilon'_i(\tilde b),\delta_{ij}\cdot m-\langle \alpha_i,\tilde \gamma'(\tilde b)\rangle)$$
since $\tilde \varepsilon'_i(m)=\delta_{ij}\cdot m$.

In particular, the inequalities $\tilde \varepsilon_i(m,\tilde b)\le 0$ for $i\in I$ imply
 $m \le  \langle \alpha_j,\tilde \gamma'(\tilde b)\rangle$ and
\begin{equation}
\label{eq:positive i-decomposition epsilon tropical}
\tilde \varepsilon'_i(\tilde b)\le 0
\end{equation}
for all $i\in I$.

Summarizing, if $(m,\tilde b)\in \widetilde {B^-_w}= \ZZ\times
\widetilde {B^-_{w'}}$ satisfies the inequalities
\eqref{eq:positive psi}, then $(m,\tilde b)$ satisfies
\eqref{eq:positive i-decomposition f right tropical} and
\eqref{eq:positive i-decomposition epsilon tropical}. Therefore,
$\tilde b\in \widetilde {B^-_{w'}}$ satisfies the inductive
hypothesis \eqref{eq:positive psi} with $w'$ for $w$. Therefore,
by this inductive hypothesis, $\tilde b={\bf 0}$. This and the
above inequalities  imply that $m=0$. This finishes the induction.
Lemma \ref{le:psiii tilde} is proved.
\end{proof}

Therefore, Proposition  \ref{pr:lambda b, b=0} is proved. \end{proof}

Now we will show that $\BB_{f_P,\Theta_{w_P}}^\lambda$ is empty if
$\lambda$ is not dominant and that $\BB_{f_P,\Theta_{w_P}}^\lambda$
is non-empty and has a highest weight element if $\lambda$ is
dominant. We need the following result.
\begin{lemma}
\label{le:bounded weights} For each $(\lambda,\tilde b)\in
\BB_{f_P,\Theta_{w_P}}^\lambda$, the co-weight $\tilde
\gamma(\lambda,\tilde b)\in X_\star(T)$ is bounded from the above
by $\lambda$, i.e.,
\begin{equation}
\label{eq:less than dominant}
\langle \mu,\tilde \gamma(\lambda,\tilde b)-\lambda\rangle\ \le 0
\end{equation}
for all $\mu\in X^\star(T)^+$.
\end{lemma}

\begin{proof} Recall from Section \ref{subsect:toric charts and positive structures}
that for each $w\in W$ the isomorphism $\eta^w:U^w\osr B^-_w$ is
$(\Theta^w,\Theta_w^-)$-positive and its inverse $\eta_w:B^-_w\osr
U^w$ is $(\Theta_w^-,\Theta^w)$-positive. Note that, by definition
\eqref{eq:inverse eta} of $\eta_w$ one has, for each $\mu\in
X^\star(T)^+$,
$$\Delta_{\mu,w^{-1}\mu}(u)=\Delta_{\mu,\mu}(u\overline w^{\,-1})=\Delta_{\mu,\mu}(\iota(b))=\mu(\gamma(b))^{-1}$$
for all $u\in U^w$, where $b=\eta_w(u)\in B^-_w$, $\gamma:B^-\to T$
is the natural projection, and $\Delta_{\mu,w^{-1}\mu}:G\to \AA^1$
is the generalized minor (see Section \ref{subsect:notation}). One
also has
$$f_{w,\chi^{st}}(b)=\chi^{st}(u)$$
for any $u\in U^w$, where $b=\eta_w(u)\in B^-_w$.

Passing to the tropicalization, we obtain
\begin{equation}
\label{eq:the difference}
\left<\mu,\tilde \gamma(\tilde b)\right>=-\tilde \Delta_{\mu,w^{-1}\mu}(\tilde a)\ , \tilde f_{w,\chi^{st}}(\tilde b)=
\widetilde {\chi^{st}}(\tilde u)
\end{equation}
for any $\tilde b\in \widetilde {B^-_w}$, $\mu\in X^\star(T)^+$,
where $\tilde u=\widetilde {{\eta_w}}(\tilde b)\in \widetilde
{U^{w}}=\Trop(U^w,\Theta^w)$, $\tilde \Delta_{\mu,w^{-1}\mu}$
is the $\Theta^w$-tropicalization of the restriction
$\Delta_{\mu,w^{-1}\mu}|_{U^w}$, and $\widetilde {\chi^{st}}:\widetilde {U^w}\to \ZZ$ is the tropicalization of the
$\Theta^w$-positive function  $\chi^{st}|_{U^w}:U^w\to \AA^1$.

Now choose  a reduced decomposition $\ii=(i_1,\ldots,i_\ell)$ of
$w\in W$ and recall from Section \ref{subsect:toric charts and
positive structures} that $\theta_\ii^+$ is a toric chart of the
class $\Theta_w$. Note that the regular function
$\Delta_{\mu,w^{-1}\mu}(\theta_\ii^+(a_1,\ldots,a_\ell))$ is a
monomial in $a_1,\ldots, a_\ell$ (see e.g., \cite[(6.3)]{bz-invent})
and
$$\chi^{st}(\theta_\ii^+(a_1,\ldots,a_\ell))=a_1+\cdots +a_\ell\ .$$
Therefore, the tropicalizations of $\chi^{st}\circ \theta_\ii^+$ and
$\Delta_{\mu,w^{-1}\mu}\circ \theta_\ii^+$ satisfy $\widetilde
{\chi^{st}}(\tilde a_1,\ldots,\tilde a_\ell)=\min\limits_{1\le k\le
\ell} \tilde a_k$ for all $(\tilde a_1,\ldots,\tilde a_\ell)\in
\ZZ^\ell$; and $\tilde \Delta_{\mu,w^{-1}\mu}(\tilde
a_1,\ldots,\tilde a_\ell)$ is a linear function in $\tilde
a_1,\ldots,\tilde a_\ell$ with non-negative coefficients. Therefore,
$\tilde \Delta_{\mu,w^{-1}\mu}(\tilde a_1,\ldots,\tilde a_\ell)\ge
0$ whenever $\widetilde {\chi^{st}}(\tilde a_1,\ldots,\tilde a_\ell)\ge
0$.

Therefore, the first equation \eqref{eq:the difference} implies
that $\left<\mu,\tilde \gamma(\tilde b)\right>\le 0$ for all $\mu\in
X^\star(T)^+$ whenever $\tilde f_{w,\chi^{st}}(\tilde b)\ge 0$.
Finally, for $w=w_P$ we have $\tilde \gamma(\lambda,\tilde b)=
\lambda+\tilde \gamma(\tilde b)$ for all $(\lambda,\tilde b)\in
\BB_{f_P,\Theta_{w_P}}^\lambda$ and \eqref{eq:less than dominant}
follows. \end{proof}

Therefore, since   weights of all elements of each non-empty normal crystal
$\BB_{f_P,\Theta_B^-}^\lambda$ are bounded from above by $\lambda$,
this normal crystal must have at least one highest  weight element. This and Proposition \ref{pr:lambda b, b=0} imply that $(\lambda,{\bf 0})$ is, indeed, the only
highest weight element in
$\BB_{f_P,\Theta_{w_P}}^\lambda$.


Theorem \ref{th:Blambda highest} is proved.
\end{proof}

According to Theorem \ref{th:positive unipotent monoidal},
$(X_B,\id,f_B,\Theta_B^-)*(X_B,\id,f_B,\Theta_B^-)$ is a strongly
positive $(U\times U,\chi^{st})$-linear bicrystal type $w_0$.
Applying Proposition \ref{pr:product normal crystals decorated}
and Proposition \ref{pr:homomorphism normal crystals}, we obtain a
homomorphism   of Kashiwara crystals
\begin{equation}
\label{eq:crystal multiplication}
*:\BB_{f_B,\Theta_B^-}\times \BB_{f_B,\Theta_B^-}\to \BB_{f_B,\Theta_B^-} \ .
\end{equation}

Furthermore, let us denote by $c_\lambda$ the unique highest weight
element of the normal crystal ${\mathcal C}_\lambda:=\BB_{f_B,\Theta_B^-}^{\lambda}$. By
definition, the product of the combinatorial crystals ${\mathcal
C}_\lambda\times {\mathcal C}_\mu$ contains the element
$(c_\lambda,c_\mu)$. Clearly, this is a unique element of the weight
$\lambda+\mu$ in ${\mathcal C}_\lambda\times {\mathcal C}_\mu$.  And
this $(c_\lambda,c_\mu)$ is a highest weight element in ${\mathcal
C}_\lambda\times {\mathcal C}_\mu$ since weights of all other
elements in ${\mathcal C}_\lambda\times {\mathcal C}_\mu$ are less
than $\lambda+\mu$. Clearly, $c_\lambda*c_\mu=c_{\lambda+\mu}$ (in
terms of \ref{eq:crystal multiplication}). Denote by ${\mathcal
C}'_{\lambda+\mu}\subset {\mathcal C}_\lambda\times {\mathcal
C}_\mu$ the preimage of ${\mathcal C}_{\lambda+\mu}$ under the
homomorphism $*$.

\begin{lemma}
\label{le:isomorphism closed family} The restriction of $*$ to
${\mathcal C}'_{\lambda+\mu}$ is an isomorphism ${\mathcal
C}'_{\lambda+\mu}\osr {\mathcal C}_{\lambda+\mu}$.

\end{lemma}

\begin{proof} Clearly,  $(c_\lambda,c_\mu)\in {\mathcal C}'_{\lambda+\mu}$. Let us show that
${\mathcal C}'_{\lambda+\mu}$ contains no other highest weight
elements. Indeed, if  $z\in {\mathcal C}'_{\lambda+\mu}$ is
another highest weight element, then $z$ has a weight less than
$\lambda+\mu$ and the image $*(z)$ is a highest weight element of
${\mathcal C}_{\lambda+\mu}$ different from $c_{\lambda+\mu}$,
which contradicts   the uniqueness of the highest weight element
in ${\mathcal C}_{\lambda+\mu}$.

Therefore,  $(c_\lambda,c_\mu)$ is the only highest weight element
of the normal crystal ${\mathcal C}'_{\lambda+\mu}$. This and Claim
\ref{cl:normal sub-crystals are equal} imply that any homomorphism
${\mathcal C}'_{\lambda+\mu}\to {\mathcal C}_{\lambda+\mu}$ is an
isomorphism. The lemma is proved.
\end{proof}

Thus, Lemma \ref{le:isomorphism closed family} asserts that the
assignment $\tilde c_{\lambda+\mu}\mapsto \tilde c_\lambda\times
\tilde c_\mu$ defines an injective homomorphism of normal crystals
${\mathcal C}_{\lambda+\mu}\hookrightarrow {\mathcal
C}_{\lambda+\mu}\times {\mathcal C}_{\lambda+\mu}$ for all
$\lambda,\mu\in X_\star(T)^+$. This and Theorem \ref{th:Blambda
highest} prove  that the family of crystals ${\mathcal
C}_\cdot=\{\BB_{f_B,\Theta_B^-}^\lambda|\lambda\in X_\star(T)^+\}$
is closed. In  turn,  Theorem \ref{th:closed family} guarantees
that each $\BB_{f_B,\Theta_B^-}^\lambda$ is isomorphic to the
associated crystal $\BB(V_\lambda)$.

Theorem \ref{th:Blambda} is proved. \end{proof}

Note that Claim \ref{cl:standard parabolic and op}(c) and Claim
\ref{cl:highest weight decomposition crystal}(c) imply the following
corollary.

\begin{corollary} For each $\lambda\in X_\star(T)^+$, one has $(\BB^\lambda_{f_B,\Theta^-_B})^{op}\cong \BB^{-w_0\lambda}_{f_B,\Theta^-_B}$.
\end{corollary}

This result agrees with Theorem \ref{th:Blambda} and the fact that the dual of the $\gg^\vee$-module $V_\lambda$ is isomorphic to $V_{-w_0\lambda}$.

\begin{example} Let $G=GL_3$, so that $T=\{t=diag(t_1,t_2,t_3)\}\subset GL_3$.
We fix a reduced decomposition $\ii=(1,2,1)$ of $w_0\in W=S_3$ and
choose a toric chart $\theta:T\times \GG_m^3 \to TB^-_w$ of the
class $\Theta_B^-$ as follows:
$$\theta(t;c_1,c_2,c_3)= t\cdot\theta_\ii^-(c_1,c_2,c_3)= \displaystyle{\begin{pmatrix}
t_1\frac{1}{c_1c_3}&0 &  0\\
t_2(\frac{c_1}{c_2}+\frac{1}{c_3})& t_2\frac{c_1c_3}{c_2}&0\\
t_3&t_3c_3&t_3c_2
\end{pmatrix}}$$
(see Example \ref{ex:GL3 charts}).
Therefore, the restriction of $f_B$ to $TB^-_w$ is given by (in the new coordinates $(t;c_1,c_2,c_3)$):
$$f_B(t;c_1,c_2,c_3)=c_1+\frac{c_2}{c_3}+c_3+\frac{t_2}{t_3}\cdot \left (\frac{c_1}{c_2}+\frac{1}{c_3}\right )+\frac{t_1}{t_2}\cdot \frac{1}{c_1} \ .$$
And the rest of the decorated geometric crystal structure ${\mathcal
X}$ on $TB^-_{w_0}$ is given by the morphism $\gamma$, the actions
$e_i^\cdot$, and the functions $\varphi_i,\varepsilon_i$, $i=1,2$:
$$\gamma(t;c_1,c_2,c_3)=\left (t_1\frac{1}{c_1c_3},t_2\frac{c_1c_3}{c_2},t_3c_2\right )\ ,$$
$$e_1^d(t;c_1,c_2,c_3)=\left(t;c_1\frac{c_2+c_1c_3}{d\cdot c_2+c_1c_3},c_2,c_3\frac{c_2+d^{-1}\cdot c_1c_3}{c_2+c_1c_3}\right)\ ,$$
$$e_2^d(t;c_1,c_2,c_3)=(t;c_1,d^{-1}\cdot c_2,c_3)\ ,$$
$$\varphi_1(t;c_1,c_2,c_3)=\frac{t_2}{t_1}\cdot \left(\frac{c_1^2c_3}{c_2}+c_1\right), \varphi_2(t;c_1,c_2,c_3)=\frac{t_3}{t_2}\cdot \frac{c_2}{c_1}\ ,$$
$$\varepsilon_1(t;c_1,c_2,c_3)= \frac{1}{c_3}+\frac{c_2}{c_1c_3^2},~ \varepsilon_2(t;c_1,c_2,c_3)=\frac{c_3}{c_2}\ .$$

The isomorphism ${\mathcal X}\osr  {\mathcal X}^{op}$ is given by
$$(t;c_1,c_2,c_3)\mapsto \left(t^{op};  \frac{t_2}{t_3}\cdot c_3^{-1},\frac{t_1}{t_3}
\cdot c_2^{-1},\frac{t_1}{t_2}\cdot c_1^{-1}\right)\ ,$$ where
$t^{op}=diag(t_1,t_2,t_3)^{op}=diag(t_3,t_2,t_1)$.

The central charge $\Delta=\Delta_{X\times X}:{\mathcal X}\times
{\mathcal X}\to \AA^1$ in these coordinates is given for
$x=(t;c_1,c_2,c_3)$, $x'=(t';c'_1,c'_2,c'_3)$ by the formula
$$\Delta(x,x')=c_1 +\frac{c_2}{c_3}+ c_3  + \frac{t'_2}{c'_3
 t'_3} + \frac{c'_1 t'_2}{c'_2t'_3}+\Delta_1+\Delta_2\ ,$$
where \begin{align*}&\Delta_1= \frac{c'_1 (c'_1+\frac{c'_2}{c'_3})
t'_2t_2 +c'_2 (t'_1+c_3c'_1t'_2) t_3 + (c_1 +\frac{c_2}{ c_3})
c'_1 c'_2 t'_3 t_2} {t_3\cdot (\frac{c'_2}{c'_3} t'_1 + c_3 c'_1
(c'_1+\frac{c'_2}{c'_3}) t'_2 + c_2 c'_1 c'_2  t'_3) }\\
 &\Delta_2 =\frac{ c'_1 (c'_1+\frac{c'_2}{c'_3}) t'_1 t'_2 t_2 + c'_2
 (t'_1+c_3c'_1t'_2) t'_3 t_1  +
(c_1+\frac{c_2}{ c_3}) c'_1 c'_2  t'_1 t'_3 t_2}{t_2\cdot (c'_1 t'_1 t'_2 +
 (c_1+\frac{c_2}{c_3}) c'_2 t'_1 t'_3 + c_1 c_3 c'_1 c'_2 t'_2 t'_3)} \
 .\end{align*}

The tropicalization of the above structures consists of:

\noindent $\bullet$ The set  $\tilde X=\Lambda^\vee\times \ZZ^3$ where $\lambda\in \Lambda^\vee=X_\star(T)= \ZZ^3$.

\noindent $\bullet$ The functions $\tilde f_B,\tilde \varphi_i,\tilde \varepsilon_i:\Lambda^\vee\times \ZZ^3\to \ZZ$, $i=1,2$:
$$\tilde f_B(\lambda;{\bf m})=\min(m_1,m_2-m_3,m_3,\lambda_2-\lambda_3-\max(m_3,m_2-m_1),\lambda_1-\lambda_2-m_1)\ ,$$
$$\tilde \varphi_1(\lambda;{\bf m})=\lambda_1-\lambda_2-\min(m_1,2m_1+m_3-m_2),~\tilde \varphi_2(\lambda;{\bf m})=\lambda_2-\lambda_3+m_1-m_2\ ,$$
$$\tilde \varepsilon_1(\lambda;{\bf m}))= \max(m_3,m_1+2m_3-m_2),~\tilde \varepsilon_2(\lambda;{\bf m})=m_2-m_3 $$
for $(\lambda;{\bf m})=((\lambda_1,\lambda_2,\lambda_3);(m_1,m_2,m_3))\in \tilde X$.

\noindent $\bullet$ The set $\tilde B_{\tilde f_B}$   which consists
of all those $(\lambda;{\bf m})\in \tilde B$ that satisfy $\tilde
f_B(\lambda;{\bf m})\ge 0$, equivalently: $m_1\ge 0,~m_2\ge m_3\ge
0,~\lambda_1-\lambda_2\ge m_1,~\lambda_2-\lambda_3\ge m_3,
\lambda_2-\lambda_3\ge m_2-m_1$. That is, each point of
$(\lambda;{\bf m})\in \tilde B$ is a Gelfand-Tsetlin pattern:
$$
\begin{pmatrix}
\lambda_1 &               & \lambda_2    &              & \lambda_3\\
          & \lambda_2+m_1 &              &\lambda_2+m_3 & \\
          &               &\lambda_3+m_2 &              &
\end{pmatrix}
$$

\noindent $\bullet$ The bijection $e_i^n:\tilde B\to \tilde B$, $i=1,2$:
$$\tilde e_1^{\,n}(\lambda;{\bf m})=\left(\lambda;m_1+\max(\delta-n,0)-\max(\delta,0),m_2,m_3+\max(\delta,0)-\max(\delta,n)\right) \ ,$$
where $\delta=m_1+m_3-m_2$,
$$\tilde e_2^{\,n}(\lambda;{\bf m})=\left(\lambda;m_1,m_2-n,m_3\right) \ .$$
Therefore, one has the decomposition into the connected components.
$${\mathcal B}_{f_B,\Theta_{w_0}^-} \cong \bigsqcup_{\lambda=(\lambda_1\ge \lambda_2\ge \lambda_3)} {\mathcal B}(V_\lambda) \ .$$

\end{example}

The following is the generalization of Theorem \ref{th:Blambda} to all standard parabolic bicrystals.

\begin{theorem}
\label{th:BlambdaL} For each $\lambda\in X_\star(Z(L_P)),$ we
have: \begin{enumerate}\item[(a)] If $\lambda$ is not dominant,
then $\BB_{f_P,\Theta_P^-}^{\lambda}$ is empty. \item[(b)]  If
$\lambda$ is dominant, then $\BB_{f_P,\Theta_P^-}^{\lambda}\cong
\BB(V_\lambda)$.
\end{enumerate}
\end{theorem}

\begin{proof} Recall that $X_w=Bw B$ is  the  Bruhat cell for each $w\in W$. Clearly, $X_w$ has a unique factorization
$$X_w=V(w)T\overline w U \ ,$$ where $V(w)=U\cap \overline{ w} U^-\overline  w^{\,-1}$ is the Schubert cell. Recall
that for any $w,w'$ such that $l(ww')=l(w)+l(w')$ the multiplication
$B^-\times B^-\to B^-$ defines  an open embedding $B_w^-\times
B_{w'}^-\hookrightarrow B_{ww'}^-$.

Let now $P$ be a standard parabolic subgroup, $L_P$ be the Levi
factor of $P$,  $U_P=U\cap P$ be the unipotent radical of $P$, and
$U_{L_P}=U\cap L_P$ be the maximal unipotent subgroup of $L_P$.
Obviously, $V(w_P)=U_P$. That is,  $X_P=U_PZ(L_P)\overline {w_P}U$
(see Example \ref{ex:Uw})


Recall from \eqref{eq:parabolic element} that $w_0^P$ is the
longest element in the Weyl group $W_P$ of $L_P$ and
$w_P=w_0^Pw_0$ so that $w_0^Pw_P=w_0$. Then the multiplication in
$B^-$ defines an open embedding (and hence a birational
isomorphism) of geometric crystals
\begin{equation}
\label{eq:factorization partial parabolic} B_{w_0^P}^-\times Z(L_P)B_{w_P}^-\hookrightarrow Z(L_P)B_{w_0}^- \ .
\end{equation}

This isomorphism is $(\Theta_{w_0^P}^-\times
\Theta_P^-,\Theta_{Z(L_P)}\cdot \Theta_{w_P})$-positive;
therefore, applying the tropicalization to \eqref{eq:factorization
partial parabolic},  one obtains an isomorphism of the
corresponding torsion-free Kashiwara crystals (in the notation of
Section \ref{subsect:from positive unipotent bicrystals to
combinatorial crystals}):
$$\BB_{\Theta_{w_0^P}^-}\times \BB_{\Theta_P^-}\osr \BB_{\Theta_{Z(L_P)}\cdot B_{w_0}^-}\ .$$
In what follows we will simply identify $\BB_{\Theta_{Z(L_P)}\cdot B_{w_0}^-}$ with the product $\BB_{\Theta_{w_0^P}^-}\times \BB_{\Theta_P^-}$.

For any sub-torus $T'\subset T$, denote temporarily
$X_{B;T'}=UT'{w_0}U$. Clearly, the quadruple
$(X_{B;T'},\id,f_P|_{X_{B;T'}}, \Theta_{T'}\cdot \Theta_{w_0}^-)$
is a strongly positive $(U\times U,\chi^{st})$-linear bicrystal; and the natural inclusion $X'_B=UT'w_0U\subset X_B=UT{w_0}U$
extends to an embedding of $(U\times U,\chi^{st})$-linear
bicrystals, and, therefore, to the natural embedding of normal
crystals
$$\BB_{f_B,\Theta_{T'}\cdot \Theta_{w_0}^-}\subset \BB_{f_B,\Theta_T\cdot
\Theta_{w_0}^-}=\BB_{f_B,\Theta_B^-}\ ,$$
or, more precisely,
$$\BB_{f_B,\Theta_{T'}\cdot \Theta_{w_0}^-}=\bigsqcup_{\lambda\in X_\star(T')}
 \BB_{f_B,\Theta_B^-}^\lambda \ .$$
 In   turn,  Theorem \ref{th:Blambda} implies that $\BB_{f_B,\Theta_{T'}\cdot
  \Theta_{w_0}^-}^\lambda\cong \BB(V_\lambda)$ if $\lambda\in X_\star(Z(L_P))\cap X_\star(T)^+$,
 and $\BB_{f_B,\Theta_{Z(L_P)}\cdot \Theta_{w_0}^-}^\lambda$ is empty if $\lambda\in  X_\star(Z(L_P))\setminus X_\star(T)^+$. That is,
\begin{equation}
\label{eq:partial decomposition T'}
 \BB_{f_B,\Theta_{T'}\cdot \Theta_{w_0}^-}\cong \bigsqcup_{\lambda\in X_\star(T')\cap X_\star(T)^+}  \BB(V_\lambda)\ .
\end{equation}

Furthermore, taking $T'=Z(L_P)$, we have a natural inclusion of Kashiwara crystals:
$$\BB_{f_B,\Theta_{Z(L_P)}\cdot \Theta_{w_0}^-}\subset  \BB_{\Theta_{Z(L_P)}\cdot \Theta_{w_0}^-}  \ .$$
According to Theorem \ref{th:Blambda highest}, the normal crystal
$\BB_{f_B,\Theta_{Z(L_P)}\cdot \Theta_{w_0}^-}^\lambda$ has a
unique highest weight element $\tilde b_\lambda$. Under the
identifications $\BB_{\Theta_P^-}= X_\star(Z(L_P))\times
\BB_{\Theta_{w_P}^-}$ and $\BB_{\Theta_{Z(L_P)}\cdot
B_{w_0}^-}=\BB_{\Theta_{w_0^P}^-}\times \BB_{\Theta_P^-}$, we have
$\tilde b_\lambda=({\bf 0},\lambda,{\bf 0}')$, where ${\bf 0}\in
\BB_{\Theta_{w_0^P}^-}$ and  ${\bf 0}'\in
\Trop(B_{w_P}^-,\Theta_{w_P}^-)$ are the marked points (see
Sections \ref{subsect:Newton polytopes and tropicalization} and
\ref{subsect:from decorated geometric crystals to combinatorial
crystals}).

Theorem \ref{th:Blambda highest} implies that for each $\lambda\in
X_\star(Z(L_P))\cap X_\star(T)^+$  the normal crystal
$\BB_{f_P,\Theta_P^-}^\lambda$ has a unique highest weight element
$\tilde b_\lambda'=(\lambda,{\bf 0}')$ (under the above
identification $\BB_{\Theta_P^-}= X_\star(Z(L_P))\times
\Trop(B_{w_P}^-,\Theta_{w_P}^-)$). Clearly,
$b_\lambda'=(\lambda,{\bf 0}')\in \overline {\BB_{\Theta_P^-}}$.
In particular, $\BB_{f_P,\Theta_P^-}$ is a normal sub-crystal of
$\overline {\BB_{\Theta_P^-}}$.

Lemma \ref{le:product trivial}(a) guarantees that the correspondence
$\tilde b\mapsto ({\bf 0},\tilde  b)$ is an injective
homomorphism of subnormal Kashiwara crystals ${\bf j}: \overline
{\BB_{\Theta_P^-}} \hookrightarrow
\overline{\BB_{\Theta_{w_0^P}^-}}\times
\overline{\BB_{\Theta_P^-}}$.

Furthermore, Theorem \ref{th:Blambda highest} guarantees that
$\BB^\lambda_{f_B,\Theta_{Z(L_P)}\cdot
\Theta_{w_0}^-}=\BB_{f_P,\Theta_P^-}^\lambda=\emptyset$ if
$\lambda\in X_\star(Z(L_P))\setminus X_\star(T)^+$.

Therefore, in order to finish the proof of Theorem \ref{th:BlambdaL}, it suffices to show that
\begin{equation}
\label{eq:j-isomorphism}
{\bf j}(\BB_{f_P,\Theta_P^-}^\lambda)=\BB_{f_B,\Theta_{Z(L_P)}\cdot \Theta_{w_0}^-}^\lambda
\end{equation}
for each $\lambda\in X_\star(Z(L_P))\cap X_\star(T)^+$.

Indeed, ${\bf j}(\tilde b'_\lambda)={\bf j}(\lambda,{\bf 0}')=({\bf
0},\lambda,{\bf 0}')=\tilde b_\lambda$, i.e., the normal crystals
${\bf j}(\BB_{f_P,\Theta_P^-}^\lambda)$ and
$\BB_{f_B,\Theta_{Z(L_P)}\cdot \Theta_{w_0}^-}^\lambda$ share the
(unique) highest weight element $\tilde b_\lambda=({\bf
0},\lambda,{\bf 0}')$. Therefore, these  normal crystals are equal
by Claim \ref{cl:normal sub-crystals are equal} and we obtain
\eqref{eq:j-isomorphism}.

Theorem \ref{th:BlambdaL} is proved. \end{proof}

\begin{corollary}
\label{cor:projection tilde m}
 Let $({\bf X},{\bf p},f,\Theta)$ be a strongly positive parabolic $(U\times U,\chi^{st})$-linear bicrystal of type $w_P$.
 Then, in the notation of Proposition \ref{pr:homomorphism normal crystals},
  for each $\lambda\in X_\star(Z(L)),$ one has:
\begin{enumerate}\item[(a)] If $\lambda$ is not dominant, then
$\BB_{f,\Theta}^\lambda$ is empty. \item[(b)] If $\lambda$ is
dominant and  $\BB_{f,\Theta}^\lambda$ is non-empty, then the
restriction of $\tilde {\bf m}$ to $\BB_{f,\Theta}^\lambda$ is a
surjective homomorphism of normal crystals
$$\BB_{f,\Theta}^\lambda \twoheadrightarrow \BB(V_\lambda)\ .$$\end{enumerate}
\end{corollary}

\subsection{From right unipotent bicrystals to crystals associated with $\bb^\vee$-modules}
\label{subsect:upper normal crystals} In this section we construct
a number of upper normal crystals (see Definition \ref{def:upper
sub-normal} above) based on certain  unipotent bicrystals which we
will refer to as {\it right linear} unipotent bicrystals.

\begin{definition} Given $U\times U$-variety ${\bf X}$, subgroups $U',U''\subset U$,
and a character $\chi:U\to \AA^1$, we say that a regular function $f:X\to \AA^1$ is  $(U'\times U'',\chi)$-{\it linear} if
$$f(u'xu'')=\chi(u')+f(x)+\chi(u'')$$
for all $u'\in U'$, $u''\in U''$, $x\in X$. In particular, we
refer to each $(e\times U,\chi)$-linear (resp. $(U\times
e,\chi)$-linear) function $f$ as a {\it right $(U,\chi)$-linear}
(resp.~{\it left $(U,\chi)$-linear}) function.

\end{definition}

Similarly to Section \ref{subsect:unipotent chi-linear
bicrystals}, for any $U$-bicrystal $({\bf X},{\bf p})$ and any
$(U'\times U'',\chi)$-linear function $f$ on $X$ we will refer to
the triple $({\bf X},{\bf p},f)$ as to $(U'\times
U'',\chi)$-linear bicrystal. In particular, if $f$  is a right
(resp.~left) $(U,\chi)$-linear function $f$, then we  will refer
to $({\bf X},{\bf p},f)$ as a {\it right} (resp.~{\it left})
$(U,\chi)$-linear bicrystal.

\begin{claim} In the notation of Claim \ref{cl:op bicrystal},
for each right $(U'\times U'',\chi)$-linear bicrystal $({\bf
X},{\bf p},f)$, the triple $({\bf X},{\bf p},f)^{op}=(({\bf
X},{\bf p})^{op},f)$ is a $(U''\times U',\chi)$-linear bicrystal.
In particular, the correspondence  $({\bf X},{\bf p},f)\mapsto
({\bf X},{\bf p},f)^{op}$ takes   \it right $(U,\chi)$-linear
bicrystals to the left ones and vice versa.

\end{claim}

Note that for any right $(U,\chi)$-linear bicrystal $({\bf X},{\bf
p},f)$ and a $(U\times U,\chi)$-linear bicrystal $({\bf Y},{\bf
p}',f')$ the convolution  product $({\bf X},{\bf p},f)*({\bf
Y},{\bf p}',f')$ given by \eqref{eq:bicrystal product} and is a
well-defined right $(U,\chi)$-linear bicrystal. In other words,
the category of $(U\times U,\chi)$-linear bicrystals acts (from
the right) on the category of right $(U,\chi)$-linear bicrystals.
(This observation is parallel to Lemma \ref{le:normal crystals act
on upper normals}, which implies that the category of normal
Kashiwara crystals acts from the right on the category of the
upper normal ones.)

Similarly to Section \ref{subsect:Positive structures on geometric
crystals and unipotent bicrystals}, a positive right
$(U,\chi^{st})$-linear bicrystal $({\bf X},{\bf p},f,\Theta)$ is a
positive  $U$-bicrystal $({\bf X},{\bf p},\Theta)$  such that the
restriction of $f$ to $X^-={\bf p}^{-1}(B^-)$ is
$\Theta$-positive. To each such positive right
$(U,\chi^{st})$-linear bicrystal, we associate a Kashiwara crystal
$\BB_{f,\Theta}$ by the formula \eqref{eq:truncated Kashiwara
crystal decorated}.

The following result is an ``upper'' analogue of Proposition
\ref{pr:normal Kashiwara decorated}.
\begin{proposition}
\label{pr:upper normal Kashiwara} Let $({\bf X},{\bf p},f,\Theta)$
be a positive  right (resp.~left) $(U,\chi^{st})$-linear
bicrystal. Then the Kashiwara crystal $\BB_{f,\Theta}$ is upper
(resp.~lower) normal.
\end{proposition}

\begin{proof} We claim that for any  $x\in X_w$, $u\in U$ and   $i\in I$, one has
\begin{equation}
\label{eq:upper right linear}
f(x_i(a)\cdot x\cdot u)=\chi^{st}(u)+f(x)+\sum\limits_{k>0} a^k f_i^{(k)}(x) \ ,
\end{equation}
where each $f_i^{(k)}(x)$ is a regular function on $X$.

Indeed, by definition,
$f(xu)=f(x)+\chi^{st}(u)$ for any $x\in X$, $u\in U$. Since the correspondence
$(a,x)\mapsto f(x_i(a)\cdot x)$ is a regular function on $\AA^1\times X$, we obtain
$$f(x_i(a)\cdot x\cdot u)=\chi^{st}(u)+f(x_i(a)\cdot x)=\chi^{st}(u)+f(x)
+\sum\limits_{k\ge 0} a^k f_i^{(k)}(x)\ .$$ But
$f_i^{(0)}(x)=f(x)$ because $f(x_i(0)\cdot x)=f(x)$. This proves
\eqref{eq:upper right linear}.

Now let ${\mathcal X}:={\mathcal F}({\bf X},{\bf
p})=(X^-,\gamma,\varphi_i,\varepsilon_i,e_i^\cdot|i\in I)$ be the
corresponding geometric crystal (as defined in \eqref{eq:unipotent
to geometric}).

Then \eqref{eq:upper right linear} taken with
$a=\frac{c-1}{\varphi_i(x)}$, $u=x_i(a')$,
$a'=\frac{c^{-1}-1}{\varepsilon_i(x)}$, $i\in \Supp \XX$, $x\in
X^-$ gives in conjunction with \eqref{eq:simple multiplicative
generator ei}:
$$f(e_i^c(x))=\frac{c^{-1}-1}{\varepsilon_i(x)}+f(x)+ \sum\limits_{k>0}
 (c-1)^k \frac{f_i^{(k)}(x)}{\varphi_i(x)^k}\ . $$
Equivalently,
$$f(e_i^c(x))=
\frac{c^{-1}}{\varepsilon_j(x)}+ f'(x)+ \sum\limits_{k>0} c^k h_i^{(k)}(x) \ ,$$
where $f'$ and each $h_i^{(k)}$ are rational functions on $X^-$.

The rest of the proof is nearly identical to the proof of Proposition \ref{pr:normal Kashiwara decorated}.
\end{proof}

Similarly to \eqref{eq:Blambda}, for each $\lambda\in X_\star(T),$
denote $\BB_{f,\Theta}^\lambda:=\BB_{f,\Theta} \cap
\BB_{\Theta}^\lambda$. We obtain the following corollary from
Proposition \ref{pr:upper normal Kashiwara}.

\begin{corollary}
\label{cor:upper normal Kashiwara lambda} Let $({\bf X},{\bf
p},f,\Theta)$ be a  positive  right (resp.~left)
$(U,\chi^{st})$-linear bicrystal. Then each non-empty Kashiwara
crystal $\BB_{f,\Theta}^\lambda$ is upper (resp.~lower) normal.
\end{corollary}

In what follows, we will construct a number of right
$(U,\chi)$-linear bicrystals. The first observation is that one
can ``truncate'' a given upper normal crystal using some right-invariant functions. Indeed, if $f$ is a right $(U,\chi)$-linear
function on $({\bf X},{\bf p})$ and $f':X\to \AA^1$ is a right
$U$-invariant function, then $f+f'$ is also right
$(U,\chi)$-linear.

\begin{claim} For any positive right $(U,\chi)$-linear bicrystal
$({\bf X},{\bf p},f,\Theta)$ and any $\Theta$-positive right
$U$-invariant function $f':X\to \AA^1$, one has a natural
injective homomorphism of upper normal crystals:
$$\BB_{f+f',\Theta}\hookrightarrow \BB_{f,\Theta} \ .$$
Moreover, for a given $f'$, the association  $({\bf X},{\bf
p},f)\mapsto ({\bf X},{\bf p},f+f')$ is a covariant invertible
functor from the category of right $(U,\chi)$-linear bicrystals
into itself.
\end{claim}

Another approach to  constructing new right $(U,\chi)$-linear
bicrystals consists of twisting the left $U$-action on $X$ with
morphisms $u_0:X\to U$.

\begin{claim} Let $({\bf X},{\bf p},f)$ be a right $(U,\chi)$-linear bicrystal and let $u_0:X\to U$ be any morphism.
Then the function $f_{u_0}$ on $X$ given by
\begin{equation}
\label{eq:fu0}
f_{u_0}(x)=f(u_0(x)\cdot x)
\end{equation}
for $x\in X$ is $(U,\chi)$-linear. More precisely, for a given
morphism $u_0$ the correspondence $({\bf X},{\bf p},f)\mapsto ({\bf
X},{\bf p},f_{u_0})$ is a covariant invertible functor from the
category of right $(U,\chi)$-linear bicrystals into itself.

\end{claim}


\begin{proposition}

\label{pr:right unipotent twisted} Let  $({\bf X},{\bf
p},f,\Theta)$ be a positive right $(U,\chi^{st})$-linear
bicrystal, let $u_0:X\to U$ be a morphism such $u_0(X^-)\subset
U^{w'}$ for some $w'\in W$, and the restriction $u_0|_{X^-}$ is a
$(\Theta,\Theta^{w'})$-positive morphism $X^-\to U^{w'}$.  Then:
\begin{enumerate} \item[(a)] The restriction of $f_{u_0}$ to $X^-$ is a
$\Theta$-positive function on $X^-$, that is, the quadruple $({\bf
X},{\bf p},f_{u_0},\Theta)$ is a positive right
$(U,\chi^{st})$-linear bicrystal. \item[(b)] The tropicalized
functions $\tilde f, \tilde f_{u_0}:\widetilde {X^-}\to \ZZ$
satisfy
\begin{equation}
\label{eq:inequality f_(u_0)}
\tilde f_{u_0}(\tilde b)\le \tilde f(\tilde b)
\end{equation}
for all $\tilde b\in \widetilde {X^-}=\Trop(X^-,\Theta)$ and,
therefore,  $\BB_{f_{u_0},\Theta}$ is an upper normal sub-crystal
of the upper normal crystal $\BB_{f_w,\Theta}$.

\end{enumerate}
\end{proposition}
\begin{proof}  We start with the following result.
\begin{lemma}
\label{le:positive right linear f} Let $({\bf X},{\bf
p},f,\Theta)$ be a positive right $(U,\chi^{st})$-linear
bicrystal. Then: \begin{enumerate}\item[(a)]  For any sequence
$\ii=(i_1,\ldots,i_\ell)\in I^\ell$, the regular function
$(\GG_m)^\ell\times X^-\to \AA^1$ given by
\begin{equation}
\label{eq:positive right linear fc}
(c_1,\ldots,c_\ell, x^-)\to f(x_{i_1}(c_1)\cdots x_{i_\ell}(c_\ell)\cdot x^-)
\end{equation}
is $\Theta_{(\GG_m)^\ell}\times \Theta$-positive. \item[(b)] For
any $w'\in W$,  the regular function $U^{w'}\times X^-\to \AA^1$
given by
\begin{equation}
\label{eq:positive right linear f}
(u, x^-)\to f(u\cdot x^-)
\end{equation}
is $\Theta^{w'} \times \Theta$-positive.
\end{enumerate}
\end{lemma}

\begin{proof} Prove (a). By the definition from Section \ref{subsect:Positive structures on geometric crystals and unipotent bicrystals},
the rational morphism $(\GG_m)^\ell\times X^-\to X^-$ given by
$$(c'_1,\ldots,c'_\ell, x^-)\to e_{i_1}^{c'_1}\cdots
e_{i_\ell}^{c'_\ell}(x^-)$$ is
$(\Theta_{(\GG_m)^\ell}\times\Theta,\Theta)$-positive. Denote
$x_k^-=e_{i_k}^{c'_k}\cdots e_{i_\ell}^{c'_\ell}(x^-)$ for
$k=1,\ldots,\ell$ and $x_{\ell+1}^-:=x^-$. Then substituting
$c'_k=c_k\varphi_{i_k}(x^-_{k+1})+1$, we see that the rational
morphism $F_\ii:(\GG_m)^\ell\times X^-\to X^-$ given by
$$F_\ii(c_1,\ldots,c_\ell; x^-)=
e_{i_1}^{c_k\varphi_{i_1}(x^-_{2})+1}\cdots
e_{i_\ell}^{c_k\varphi_{i_\ell}(x^-_{\ell+1})+1}(x^-)$$ is positive.

Taking into  account \eqref{eq:simple multiplicative generator ei}
and \eqref{eq:unipotent crystal}, we see that the latter positive
morphism is given by
 $$F_\ii(c_1,\ldots,c_\ell; x^-)=x_{i_1}(c_1)\cdots x_{i_\ell}(c_\ell)\cdot x^-\cdot x_{i_\ell}(-c''_\ell)\cdots x_{i_1}(-c''_1)\ , $$
 where $c''_k= {\frac{c_k}{\alpha_{i_k}(\gamma(x_{k+1}^-))(1+c_k\varphi_{i_k}(x_{k+1}^-))}}$ for $k=1,\ldots,\ell$.

Using the right $(U,\chi^{st})$-linearity of $f$, we obtain:
$$f(x_{i_1}(c_1)\cdots x_{i_\ell}(c_\ell)\cdot x^-)=f(F_\ii(c_1,\ldots,c_\ell; x^-))+\sum_{k=1}^\ell c''_k \ .$$
Therefore, the positivity of $F_\ii$,  $f|_{X^-}$, and of each
function $c_k'':(\GG_m)^\ell\times X^-\to \AA^1$ implies the
positivity of the function given by \eqref{eq:positive right linear
fc}. This proves (a).

Part (b) follows from (a) by taking $\ii=(i_1\ldots,i_\ell)$ to be a
reduced decomposition for $w'$. Lemma \ref{le:positive right linear
f} is proved.
\end{proof}


Since  $u_0|_{X^-}:X^-\to U^{w'}$ is a
$(\Theta,\Theta^{w'})$-positive morphism and \eqref{eq:positive
right linear f} defines a  $\Theta^{w'}\times \Theta$-positive
function by Lemma \ref{le:positive right linear f}(b),  it follows
that  the restriction of $f_{u_0}$ to $X^-$ is a $\Theta$-positive
function on $X^-$.  This proves part (a) of Proposition
\ref{pr:right unipotent twisted}.

Prove (b) now. Taking again $\ii=(i_1,\ldots,i_\ell)$ to be a
reduced decomposition of $w'$ and
$u=\theta_\ii^+(c_1,\ldots,c_\ell)=x_{i_1}(c_1)\cdots
x_{i_\ell}(c_\ell)$, we obtain for all $x\in X$ (based on Claim
\ref{cl:u-shift}):
$$f(x_{i_1}(c_1)\cdots x_{i_\ell}(c_\ell)\cdot x)=f(x)+\sum_{{\bf n}\in \ZZ_{\ge 0}^\ell\setminus \{\bf 0\}} c^{\bf n}\cdot f_{\bf n}(x) \ ,$$
where we abbreviated $c^{\bf n}=c_1^{n_1}\cdots c_\ell^{n_\ell}$
for ${\bf n}=(n_1,\ldots ,n_\ell)$, and each $f_{\bf n}$ is a
regular function on $X$. If we denote by $\tilde f, \tilde f_{\bf
n}$,  the respective $\Theta$-tropicalizations of $f|_{X^-}$ and
$f_{\bf n}|_{X^-}$, and by $\tilde F_\ii$ the tropicalization of
the function given by \eqref{eq:positive right linear fc}, then
Corollary~\ref{cor:trop of sum}(c) guarantees that
$$\tilde F_\ii({\bf m},\tilde x^-)=\min\bigg (\tilde f(\tilde x^-),
\min\limits_{{\bf n}\in \ZZ_{\ge 0}^\ell\setminus \{\bf 0\}}
\big({\bf m}\cdot {\bf n}+\tilde f_{\bf n}(\tilde x^-) \big)\bigg)
$$ for all ${\bf m}\in \ZZ^\ell$, $\tilde x^-\in \widetilde
{X^-}=\Trop(X^-,\Theta)$, where we abbreviated ${\bf m}\cdot {\bf
n}=m_1 n_1+\cdots +m_\ell n_\ell$ for ${\bf m}=(m_1,\ldots
,m_\ell)$, ${\bf n}=(n_1,\ldots ,n_\ell)$. This implies that
$$\tilde F_\ii({\bf m},\tilde x^-)\le \tilde f(\tilde x^-)$$
for all ${\bf m}\in \ZZ^\ell$, $\tilde x^-\in \widetilde {X^-}$.
More invariantly, if we denote by $\tilde F:\widetilde
{U^{w'}}\times \widetilde {X^-}\to \ZZ$ the $\Theta^{w'}\times
\Theta$-tropicalization of the positive function given by
\eqref{eq:positive right linear f}, we obtain
$$\tilde F(\tilde u,\tilde x^-)\le \tilde f(\tilde x^-)$$
for all $\tilde u\in \widetilde {U^{w'}}$, $\tilde x^-\in
\widetilde {X^-}$. Taking into   account that $\tilde
f_{u_0}(\tilde x^-)=\tilde F(\tilde u_0(\tilde x^-),\tilde
x^-)$, where $\tilde u_0:\widetilde {X^-}\to \widetilde {U^{w'}}$
is the tropicalization of the positive morphism $u_0$, we obtain
the inequality \eqref{eq:inequality f_(u_0)}. This proves (b).

Therefore, Proposition \ref{pr:right unipotent twisted} is proved.
\end{proof}


Our main example of right  $(U,\chi)$-linear functions is as follows.

\begin{lemma} The function $f_{w,\chi}$ on $(BwB,\id)$ defined in \eqref{eq:fwst} is right $(U,\chi)$-linear.
\end{lemma}

\begin{proof} Indeed, for $x\in BwB$, $u\in U$, we have
\begin{align*} f_{w,\chi}(xu)&=\chi(\pi^+(\overline
w^{\,-1}xu))=\chi(\pi^+(\overline w^{\,-1}x)u)\\
 &=\chi(\pi^+(\overline w^{\,-1} x))+\chi(u)=f_{w,\chi}(x)+\chi(u) \
 .\end{align*}

This proves the lemma.
\end{proof}

\begin{corollary}
 Let  $({\bf X},{\bf p})$ be any $U$-bicrystal such that ${\bf p}(X)\subset BwB$ for some  $w\in W$.
 Then for any morphism $u_0:X\to U$ and a character $\chi:U\to \AA^1$, the function $f_{w,\chi,u_0}:X\to \AA^1$
 given by $f_{w,\chi,u_0}(x)=f_{w,\chi}(u_0(x)\cdot {\bf p}(x))$ is right $(U,\chi)$-linear.
 That is, the correspondence $({\bf X},{\bf p})\mapsto ({\bf X},{\bf p},f_{w,\chi,u_0})$ is a functor
 from a full sub-category of the category of $U$-bicrystals into the category of right $(U,\chi)$-linear bicrystals.
\end{corollary}

\begin{remark}  The requirement ${\bf p}(X)\subset BwB$ is not restrictive at all because for each unipotent bicrystal
$({\bf X},{\bf p})$ of type $w$, there exists a dense $U\times
U$-invariant subvariety $X_0$ of $X$ such that ${\bf
p}(X_0)\subset BwB$.

\end{remark}

For each $w\in W$, denote simply by $f_w$ the restriction of the
function $f_{w,\chi^{st}}$ to $\overline X_w=U\overline w U$ (see Example
\ref{ex:standard bicrystal}). Then Claim \ref{cl:standard
bicrystal} implies the following result.

\begin{claim}
\label{cl:UwU}The quadruple $(\overline {\bf X}_w, \id,f_w,\Theta_w^-)$ is a
positive right $(U,\chi^{st})$-linear bicrystal. Therefore, in the
notation of  Proposition \ref{pr:upper normal Kashiwara}, the
Kashiwara crystal $\BB_{f_w,\Theta^-_w}$ is upper normal.

\end{claim}


Note that $\BB_{f_w,\Theta^-_w}^\lambda$ is empty unless $\lambda=0$, and $\BB_{f_w,\Theta^-_w}^0=\BB_{f_w,\Theta^-_w}$.

\begin{theorem}
\label{th:asociated crystal Schubert} For any standard parabolic
subgroup $P$ of $G$, the upper normal crystal
$\BB_{f_{w_P},\Theta^-_{w_P}}$ is isomorphic to the associated
crystal $\BB(\CC[U_P^\vee])$.
\end{theorem}

\begin{proof}   In   view of Corollary \ref{cor:directed family},
it will suffice to show that $\BB_{f_{w_P},\Theta^-_{w_P}}$ is also a
limit of the directed family $(\BB(V^\lambda),\tilde {\bf
f}_{\lambda,\mu})$ (where  $\BB(V^\lambda)$ is the crystal
associated to the $\bb^\vee$-module $V^\lambda$, which is the
restriction of simple $\gg^\vee$-module to $\bb^\vee$).

Define the   projection $pr_P:X_P=UZ(L_P)\overline {w_P}U\to \overline X_{w_0}=U\overline {w_0}U$  by
\begin{equation}
\label{eq:projection pr_0}
pr_P(tu\overline {w_P}u')=u\overline {w_P}u'
\end{equation}
for each $u,u'\in U$, $t\in Z(L_P).$ Clearly, $pr_P$ defines a
surjective morphism of unipotent bicrystals $({\bf X}_P,\id)\to
(\overline {\bf X}_{w_P},\id)$, which is
$(\Theta_P^-,\Theta_{w_P}^-)$-positive and, therefore, it defines
a projection  $\widetilde {pr_P}$ of marked sets $\widetilde
{Z(L_P)B_{w_0}^-}=X_\star(Z(L_P))\times \widetilde {B_{w_P}^-}\to
\widetilde {B_{w_P}^-}$ (see Sections \ref{subsect:Newton
polytopes and tropicalization} and \ref{subsect:from decorated
geometric crystals to combinatorial crystals}), where $\widetilde
{Z(L_P)B_{w_P}}=\Trop(Z(L_P)B_{w_P}^-,\Theta_P^-)$ and $\widetilde
{B_{w_P}}=\Trop(B_{w_P}^-,\Theta_{w_P}^-)$.

According to Theorem \ref{th:BlambdaL}, the normal crystal
$\BB_{f_P,\Theta^-_P}$ being considered as an upper normal crystal
is isomorphic to the union of all associated crystals $\BB(V^\lambda)$.

\begin{claim} For each standard parabolic subgroup $P\subset G$, one has:
\begin{enumerate} \item[(a)] The restriction of $\widetilde
{pr_P}$ to  the  normal crystal $\BB_{f_P,\Theta^-_P}$ is a
surjective homomorphism of upper normal Kashiwara crystals
$\BB_{f_P,\Theta^-_P}\to \BB_{f_{w_P},\Theta^-_{w_P}}$. \item[(b)]
For each $\lambda\in X_\star(Z(L_P))\cap X_\star(T)^+,$ the
restriction of $\widetilde {pr_P}$ to $\BB(V_\lambda)\subset
\BB_{f_P,\Theta^-_P}$ is an injective homomorphism $\tilde {\bf
j}'_\lambda:\BB(V^\lambda) \hookrightarrow
\BB_{f_{w_P},\Theta^-_{w_P}}$ of upper normal Kashiwara crystals.
\item[(c)] For each $\lambda,\mu\in X_\star(Z(L_P))\cap
X_\star(T)^+$, one has
$$\tilde {\bf j}'_\lambda(\BB(V^\lambda))\subset \tilde {\bf j}'_{\lambda+\mu}(\BB(V^{\lambda+\mu})$$
so that the induced homomorphism $\BB(V^\lambda)\hookrightarrow
\tilde \BB(V^{\lambda+\mu})$ equals $\tilde {\bf f}_{\lambda,\mu}$
from Corollary \ref{cor:directed family}.
\end{enumerate}
\end{claim}

Therefore, the   limit of  $(\BB(V^\lambda),\tilde {\bf
f}_{\lambda,\mu})$, $\lambda,\mu\in X_\star(Z(L_P))\cap
X_\star(T)^+$ is isomorphic to $\BB_{f_{w_P},\Theta^-_{w_P}}$. The
uniqueness of the limit and Corollary \ref{cor:directed family}
finish the proof of Theorem \ref{th:asociated crystal Schubert}.
\end{proof}

\begin{remark} An isomorphism $\BB(\CC[U_P])\osr \BB_{f_{w_P},\Theta^-_{w_P}}$ was constructed in \cite{bz-invent} by, first,
choosing a reduced decomposition $\ii$ of $w_P$ and, second, using an explicit Kashiwara parametrization of the dual canonical basis of $\CC[U_P]$.

\end{remark}



\section{Conjectures and open questions}

\subsection{Tropicalization and associated crystals}
We start the following conjecture which complements results of Section \ref{subsect:from positive unipotent bicrystals to combinatorial crystals}.

\begin{conjecture} \label{con:injective components} Let  $({\bf X},{\bf p},f,\Theta)$  be a positive $(U\times U,\chi^{st})$-linear bicrystal
of type $w$. Then the restriction of the structure map $\tilde {\bf
f}_w:\BB_\Theta\to \BB_{\Theta_{B^-}}$ (see \eqref{eq:homomorphism
free crystals}) to each connected component of $\BB_{f,\Theta}$ is
injective.

\end{conjecture}

\begin{remark} Informally speaking, Conjecture \ref{con:injective components}
means that $U\times U$-orbits in $X$ correspond to the components in
$\BB_\Theta$ or, more precisely, we expect a correspondence between
classes of orbits in $X^-$ under the rational $U$-action
\eqref{eq:new U action} and connected components of $\BB_\Theta$.
\end{remark}

The following is an immediate corollary from this conjecture and a refinement of Corollary \ref{cor:projection tilde m}.

\begin{conjecture}
\label{con:strongly positive crystal bases} Let $(X,{\bf
p},f,\Theta)$ be a strongly positive parabolic $(U\times
U,\chi^{st})$-linear bicrystal  of type $w_P$. Then each non-empty
the normal crystal $\BB_{f,\Theta}^\lambda$ is isomorphic to the
union of copies of the associated crystal $\BB(V_\lambda)$.

\end{conjecture}

\noindent{\bf Proof of the implication Conjecture
\ref{con:injective components}$\ \Rightarrow$ Conjecture
\ref{con:strongly positive crystal bases}}. By
Theorem~\ref{th:BlambdaL}, each non-empty crystal
$\BB_{f_P,\Theta_P^-}^\lambda$ is isomorphic to the crystal
$\BB(V_\lambda)$ associated to the finite-dimensional
$\gg^\vee$-module $V_\lambda$.

Denote by $\BB_0$  a connected component of
$\BB_{f,\Theta}^\lambda$. We have to prove that $\BB_0\cong
\BB_{f_P,\Theta_P^-}^\lambda$. According to Claim
\ref{cl:connected normal crystal}, $\BB_0$ is also a  normal
crystal. Since the homomorphism $\tilde {\bf f}$ from Proposition
\ref{pr:homomorphism normal crystals} commutes with the
tropicalization of the highest weight morphism $hw_X$, we obtain
$\tilde {\bf f}(\BB_0)\subset \BB_{f_P,\Theta_P^-}^\lambda$. Since
$\tilde {\bf f}(\BB_0)$ is normal and
$\BB_{f_P,\Theta_P^-}^\lambda$ is connected, we obtain by Claim
\ref{cl:connected normal crystal} that  $\tilde {\bf f}(\BB_0)=
\BB_{f_P,\Theta_P^-}^\lambda$. But according to Conjecture
\ref{con:injective components}, the restriction of $\tilde {\bf
f}$ to $\BB_0$ is an injective map $\BB_0 \hookrightarrow
\BB_{f_P,\Theta_P^-}^\lambda$. Therefore, the restriction of
$\tilde {\bf f}$ to $\BB_0$ is an isomorphism $\BB_0 \osr
\BB_{f_P,\Theta_P^-}^\lambda$. \hfill\quad \qed

\bigskip

For a complex projective variety $X^\vee$, we denote by $\widehat
{X^\vee}$ the affine cone over $X$. The following is an equivalent
reformulation of Conjecture \ref{con:strongly positive crystal
bases}.

\begin{conjecture}
\label{con:module} Let  $({\bf X},{\bf p},f,\Theta)$ of type $w_P$
be a strongly positive  parabolic $(U,\chi^{st})$-linear
bicrystal. Then there exist a based $\gg^\vee$-module
$(V_\Theta,{\bf B}_\Theta)$ and a $\gg^\vee$-linear map
$h:V_\Theta\to \CC[\widehat {G^\vee/P^\vee}]$ such that the
associated normal crystal $\BB(V_\Theta,{\bf B}_\Theta)$ is
isomorphic to   $\BB_{f,\Theta}$ and
 $h({\bf B}_\Theta)$ is a perfect basis for $\CC[\widehat {G^\vee/P^\vee}]$.
\end{conjecture}

In fact, we expect that the (infinite-dimensional) module $V_\Theta$
is the coordinate algebra $\CC[\widehat {X^\vee}]$ of some projective
$G^\vee$-variety $X^\vee$. Theorems \ref{th:asociated crystal Schubert}
guarantees that the variety $X^\vee=U_P^\vee$ (with the dual Levi
factor $L_P^\vee$ for $G^\vee$) is a suitable candidate.
However, we do not expect that any $G^\vee$-variety to $X^\vee$ is
related to a unipotent bicrystal.
 In fact, we will show in a separate paper that
for $G^\vee=GL_2(\CC)$, and the $4$-dimensional simple $G^\vee$-module
$Y^\vee=S^3(\CC^2)$ there is no unipotent bicrystal that would parametrize
the associated crystal $\BB(\CC[Y^\vee])$.

%
%
%
%

Below (Conjecture \ref{con:saturated valuation}), we will lay out
sufficient conditions which would allow to construct $(U\times U,\chi^{st})$-linear unipotent bicrystals for certain
$G^\vee$-varieties.

\begin{definition} Given a  commutative $\CC$-algebra ${\mathcal A}$ without zero divisors and a totally ordered free abelian group $\Gamma$, a map
$\nu:{\mathcal A}\setminus \{0\}\to \Gamma$ is said to be a {\it valuation} if
$$\nu(xy)=\nu(x)+\nu(y)$$
for all $x,y\in {\mathcal A}\setminus \{0\}$,
$$\nu(x+y)=\min(\nu(x),\nu(y))$$
for any $x,y\in {\mathcal A}\setminus \{0\}$ such that $\nu(x)\ne \nu(y)$ (we use the convention that
$\nu(0)=+\infty$, where $+\infty$ is greater than any element of
$\Gamma$).


We say that $\nu$ is {\it saturated} if for each $\lambda\in
\nu({\mathcal A}\setminus \{0\}),$ the entire ``half-line''
$\QQ_{\ge 0}\cdot \lambda \cap \Gamma$ also belongs to the
semi-group $\nu({\mathcal A}\setminus \{0\})$.)

\end{definition}

\begin{conjecture}
\label{con:saturated valuation} Let $Y^\vee$ be an affine
$G^\vee$-variety, $S$ be a split algebraic torus, and let $\prec$
be a total ordering on the co-character lattice $X_\star(S)$. Assume
that there exist a saturated valuation  $\nu:\CC[Y^\vee]\setminus
\{0\}\to (X_\star(S),\prec)$ and a perfect basis ${\bf B}$ for
$\CC[Y^\vee ]$ (see Section \ref{subsect:Perfect bases of
gg-modules}) such that the restriction of $\nu$ to ${\bf B}$ is an
injective map
\begin{equation}
\label{eq:crystal isomorphism}
{\bf B} \hookrightarrow   X_\star(S) \ .
\end{equation}
 Then there exists a strongly positive $(U,\chi^{st})$-linear bicrystal
  $({\bf X},{\bf p},f,\Theta)$ (where $\theta:S\osr X^-$) such that the image of \eqref{eq:crystal isomorphism} is $\BB_{f,\Theta}\subset  X_\star(S)$, i.e.,
   \eqref{eq:crystal isomorphism} defines
  an isomorphism of normal crystals $\BB(\CC[Y^\vee])\osr \BB_{f,\Theta}$.

\end{conjecture}

\subsection{Schubert cells and upper normal crystals} In the notation of Section
 \ref{subsect:Perfect bases of b-modules},
let  $B^\vee\backslash G^\vee$ be the (right) flag variety for
$G^\vee$. For $w\in W$ let $X_w^\vee:=B^\vee \backslash B^\vee
wB^\vee$ be the Schubert cell. Clearly, one has a
$B^\vee$-equivariant isomorphism $X_w^\vee\cong B^\vee(w)\backslash
B^\vee$, where $B^\vee(w)=B^\vee\cap w^{-1}B^\vee w$ (Note that if
$w={w_P}^{-1}$ is the inverse of the parabolic element $w_P\in W$
defined in \eqref{eq:parabolic element}, then
$X_{{w_P}^{-1}}^\vee=U_P^\vee$, in particular,
$X_{w_0}^\vee=U^\vee$).

Therefore, one has a $B^\vee$-equivariant surjective map
$\pi_w:U^\vee\twoheadrightarrow X_w^\vee$. In   turn, this defines
an embedding of coordinate algebras (and locally finite
$\bb^\vee$-modules)\linebreak
$\pi^*_w:\CC[X_w^\vee]\hookrightarrow \CC[U^\vee]$.  In
particular, the above embedding $\CC[U_P^\vee]\hookrightarrow
\CC[U^\vee]$ constructed in Section \ref{subsect:Perfect bases of
b-modules} is $\pi^*_{{w_P}^{-1}}$.

\begin{conjecture}
\label{conj:perfect basis Schubert} For each $w\in W$, there
exists a perfect basis of $\CC[X_w^\vee]$.

\end{conjecture}

It follows from  Proposition \ref{pr:perfect basis parabolic} that the conjecture is true for $w={w_P}^{-1}$.

Conjecture \ref{conj:perfect basis Schubert} would imply (based on
Definition \ref{def:general upper perfect basis}) the existence of
the associated crystal $\BB(\CC[X_w^\vee])$ for each $w\in W$
(where the coordinate algebra $\CC[X_w^\vee]$ of the Schubert cell
$X_w^\vee$ is regarded as a locally finite $\bb^\vee$-module). The
following is a (yet conjectural) refinement of Conjecture
\ref{conj:perfect basis Schubert}.

\begin{conjecture}
\label{con:asociated crystal Schubert} For any $w\in W$, the
crystal basis $\BB(\CC[X_{w^{-1}}^\vee])$ is isomorphic to the
upper normal crystal $\BB_{f_w,\Theta^-_w}$ (see Claim
\ref{cl:UwU}).
\end{conjecture}

We can provide the following partial justification of Conjecture \ref{con:asociated crystal Schubert}.

Note that the complex torus $T^\vee$ (see Section
\ref{subsect:Perfect bases of b-modules}) is dual to $T$, e.g.,
each co-character $\mu$ of $T$ is a character of $T^\vee$. Indeed,
for each $\mu\in X_\star(T)$, denote by $[\mu]:T^\vee\to
\CC^\times$ the corresponding character.

For each locally finite $B^\vee$-module $V$,  define  the
character $ch (V)\in \QQ[[T^\vee]]$ by the formula
$$ch(V)=\sum_{\mu\in X_\star(T)} (\dim_\CC V(\mu))\cdot [\mu] \ ,$$
where $V(\mu)$ is the $\mu$-th weight component of $V$ (see Section \ref{subsect:Perfect bases of b-modules}).

For each  Kashiwara crystal $\BB$ such that all  fibers of $\tilde \gamma:\tilde B\to X_\star(T)$ are finite, define the {\it
character} $ch(\BB)$ by the formula
$$ch(\BB)=\sum_{\tilde b\in \BB} [\tilde \gamma(\tilde b)] \ .$$

\begin{claim} If $\BB=\BB(V)$ is the crystal associated to $V$, then  $ch(V)=ch(\BB)$.
\end{claim}

\begin{lemma} For each $w\in W$, one has $ch(\CC[X_{w^{-1}}^\vee])=ch(\BB_{f_w,\Theta^-_w})$.
\end{lemma}


\begin{proof} We  need the following well-known facts.

\begin{claim}
\label{cl:character Schubert cell} For each $w\in W$, the
character $ch(\CC[X_{w^{-1}}^\vee])$ is given by the formula
$$ch(\CC[X_{w^{-1}}^\vee])=\prod_{\alpha^\vee\in R_+^\vee\cap
 w^{-1}(-R_+^\vee)} \frac{1}{1-[\alpha^\vee]}\ ,$$
where $R_+^\vee\subset X_\star(T)$ is the
set of positive coroots of $G$.

\end{claim}

\begin{claim} \cite[Formula (6.3)]{bz-invent}
\label{cl:corner minor} For each reduced decomposition
$\ii=(i_1,\ldots,i_\ell)$ of an element $w\in W$, one has
\begin{equation}
\label{eq:factorization plus}
\Delta_{\mu,w^{-1}\mu}(\theta_\ii^+(c_1,\ldots,c_\ell))=\prod_{k=1}^\ell c_k^{\langle \mu, \alpha_{(k)}^\vee \rangle}\ ,
\end{equation}
where $\alpha_{(k)}^\vee=s_{i_1}\cdots
s_{i_{k-1}}\alpha_{i_k}^\vee$, $k=1,\ldots,\ell$ is a normal
ordering of $R_+^\vee\cap w^{-1}(-R_+^\vee)$.

\end{claim}

Denote by $\widetilde {\eta^w}:\widetilde {U^w}\to \widetilde
{B^-_w}$ the $(\Theta^w,\Theta_w^-)$-tropicalization of the
positive isomorphism $\eta^w:U^w\to B^-_w$ (See \eqref{eq:etaw}
and Claim \ref{cl:positive eta}). By Claim \ref{cl:positive eta},
$\widetilde {\eta^w}$ is a bijection $\widetilde {U^w}\osr
\widetilde {B^-_w}$.

Since $\Delta_{\mu,\mu}(b)=\Delta_{\mu,w\mu}(\eta^{w^{-1}}(b))$
for each $b\in B^-_w$, $\mu\in X^\star(T)^+$, applying the
tropicalization to \eqref{eq:factorization plus} with respect to
the toric chart $\theta_\ii^+:(\GG_m)^\ell\osr U^w$(and,
therefore, identifying  $\widetilde {U^w}$ with $\ZZ^\ell$), we
obtain
$$\tilde \gamma(\tilde b)=\sum_{k=1}^\ell \tilde u_k\alpha^\vee_{(k)}\ ,$$
where $$(\tilde u_1,\ldots,\tilde u_\ell)=(\widetilde
{\eta^w})^{-1}(\tilde b)$$ and we identified $\widetilde {U^w}$
with $\ZZ^\ell$ using the toric chart $\theta_\ii^+$. Note also
that $\tilde f_w(\tilde b)\ge 0$ if and only if $\tilde u_k\ge 0$
for all $k$. Therefore,
\begin{align*} ch(\BB_{f_w,\Theta^-_w})&=\sum_{\tilde u_1,\ldots,\tilde
u_\ell\in (\ZZ_{\ge 0})^\ell} \left [\sum_{k=1}^\ell \tilde
u_k\alpha^\vee_{(k)}\right ] =\sum_{\tilde u_1,\ldots,\tilde
u_\ell\in (\ZZ_{\ge 0})^\ell}\prod_{k=1}^\ell \left
[\alpha^\vee_{(k)}\right ]^{\tilde u_k}\\ &=\prod_{k=1}^\ell
\frac{1}{1-[\alpha_{(k)}^\vee]} \ .\end{align*} Since
$\alpha_{(1)}^\vee,\ldots, \alpha_{(\ell)}^\vee$ is a
linear ordering of the set $R_+^\vee\cap w^{-1}(-R_+^\vee)$, we
obtain in conjunction with Claim \ref{cl:character Schubert cell}
the identity $ch(\CC[X_{w^{-1}}^\vee])=ch(\BB_{f_w,\Theta^-_w})$.

The lemma is proved. \end{proof}






\subsection{Tensor product multiplicities and  combinatorial central charge}

We start with a conjecture which is a complement of Claim \ref{cl:geometric tensor product}.

\begin{conjecture}
\label{con:geometric tensor product}  The inverse of \eqref{eq:positive F}
is an isomorphism of positive decorated geometric crystals
$$({\mathcal Z}_{w_0},\Theta_Z) \osr (\XX_B,\Theta_B)\times (\XX_B,\Theta_B) \ .$$

\end{conjecture}

Denote by ${\mathcal B}_{f_Z,\Theta_Z}$ the normal Kashiwara crystal
obtained  by the tropicalization of \eqref{eq:Zwnot} with respect to
the strongly positive structure $\Theta_Z$ from Claim
\ref{cl:strongly positive Zwnot}.

By definition  \eqref{eq:Zw}, one has an invariant projection
${\mathcal Z}_{w_0}\to T\times T$, which is obviously
$(\Theta_Z,\Theta_T\times \Theta_T)$-positive; and denote by
$\tilde \pi':{\mathcal B}_{\Theta_Z}\to X_\star(T)\times
X_\star(T)$ its tropicalization. In fact, the restriction of
$\tilde \pi$ to the normal sub-crystal ${\mathcal
B}_{f_Z,\Theta_Z}$ is an invariant projection ${\mathcal
B}_{f_Z,\Theta_Z}\to X_\star(T)^+\times X_\star(T)^+$. For each
$\lambda,\nu\in X_\star(T)^+$, we denote by ${\mathcal
B}_{f_Z,\Theta_Z;\lambda,\nu}$ the fiber of the latter projection.

By definition, we have a decomposition into the normal Kashiwara
sub-crystals
$${\mathcal B}_{f_Z,\Theta_Z}=\bigsqcup_{\lambda,\nu\in X_\star(T)^+} {\mathcal B}_{f_Z,\Theta_Z;\lambda,\nu} \ .$$

\begin{lemma}
\label{le:tensor multiplicity} For each $\lambda,\nu\in
X_\star(T)^+$, the component ${\mathcal
B}_{f_Z,\Theta_Z;\lambda,\nu}$ is isomorphic to the associated
crystal $\BB(V_\lambda\otimes V_\nu)$.

\end{lemma}

\begin{proof} One can easily see that each choice of toric chart ${\bf j}_{\ii,\ii'}\in \Theta_Z$
defines an isomorphism of normal Kashiwara crystals
$${\mathcal B}_{f_Z,\Theta_Z;\lambda,\nu}\cong \bigsqcup_{\mu\in X_\star(T)^+} {\mathcal C}_{\lambda,\nu}^\mu
\times \BB_{f_B,\Theta^-_B}^\mu \ ,$$ where each ${\mathcal
C}_{\lambda,\nu}^{\mu}={\mathcal
C}_{\lambda,\nu}^{\mu}(\ii)$ is a finite set considered
as a trivial Kashiwara crystal.  It is easy to see that ${\mathcal
C}_{\lambda,\nu}^{\mu}$ is precisely the set defined by
conditions (1)-(4) in\linebreak \cite[Theorem 2.3]{bz-invent}. On
the other hand, the latter result asserts that the cardinality of
${\mathcal C}_{\lambda,\nu}^{\mu}$ is equal to the
multiplicity of $V_\mu$ in the tensor product $V_\lambda\otimes
V_\nu$. This and the isomorphism $\BB_{f_B,\Theta^-_B}^\mu\cong
\BB(V_\mu)$ finish the proof of the lemma.
\end{proof}

Then  Conjecture \ref{con:geometric tensor product} implies the following refinement of Lemma \ref{le:tensor multiplicity}.

\begin{conjecture}
\label{con:combinatorial tensor product} The tropicalization of
\eqref{eq:positive F} is an isomorphism of torsion-free Kashiwara
crystals  $\BB_{\Theta^-_B}\times \BB_{\Theta^-_B}\to {\mathcal
B}_{f_Z,\Theta_Z}$ and its restriction to
$\BB_{f_B,\Theta^-_B}\times \BB_{f_B,\Theta^-_B}$ is an isomorphism
of normal Kashiwara crystals $\tilde {\bf
F}_{w_0}:\BB_{f_B,\Theta^-_B}\times \BB_{f_B,\Theta^-_B}\to
{\mathcal B}_{f_Z,\Theta_Z}$. The restriction $\tilde {\bf F}_{w_0}$
to each $\BB_{f_B,\Theta^-_B}^\lambda \times
\BB_{f_B,\Theta^-_B}^\nu$, $\lambda,\nu\in X_\star(T)^+$  is an
isomorphism of Kashiwara crystals
\begin{equation}
\label{eq:F lambda mu} \BB (V_\lambda)\times \BB(V_\nu) \osr
{\mathcal B}_{f_Z,\Theta_Z;\lambda,\nu} .
\end{equation}

\end{conjecture}

\begin{remark} Note that the positivity of the birational isomorphism \eqref{eq:positive F} implies only
surjectivity of each \eqref{eq:F lambda mu}. Injectivity of
\eqref{eq:F lambda mu} follows from \cite[Theorem~2.3]{bz-invent}
or from Conjecture \ref{con:geometric tensor product} above.
Conversely, Claim \ref{cl:geometric tensor product} and Conjecture
\ref{con:geometric tensor product} taken together imply Theorems
2.3 and 2.4 of \cite{bz-invent}.
\end{remark}

\label{subsect:Crystal multiplication and the central charge} Let
$({\bf X},{\bf p},f,\Theta)$ be a strongly positive parabolic
$(U\times U,\chi^{st})$-linear bicrystal of type $w_P$ and let
$\kappa:X\to T'$ be a $U\times U$-invariant morphism (where $T'$
is an algebraic torus) such that the restriction of $\kappa$ to
$X^-$ is  $(\Theta,\Theta_{T'})$-positive. This data defines a
normal Kashiwara crystal $\BB_{f,\tilde \Theta}$ along with
$\tilde e_i^{\,n}$-invariant maps $\tilde \kappa:\BB_{f,\Theta}\to
X_\star(S)$ and $\tilde \Delta_X:\BB_{f,\Theta}\to \ZZ$.
For each $\mu\in X_\star(S),$ let  $\BB_{\lambda'}=\{\tilde b\in
\BB:\tilde \kappa(\tilde b)=\lambda'\}$, i.e., $\BB_{\lambda'}$ is the
$\lambda'$-th fiber of $\kappa$. By definition, $\BB_{\lambda'}$
is a normal sub-crystal of $\BB$. For each $\mu\in
X_\star(Z(L_P))\cap X_\star(T)^+$, denote by ${\mathcal
C} _{\lambda'}^{(\mu)}$ the set of connected components in $\BB _{\lambda'}$
of type (i.e., of the highest weight) $\mu$. In particular, if
$\BB _{\lambda'}$ is an associated crystal of some based $\gg^\vee$-module, then ${\mathcal
C}_{\lambda'}^{(\mu)}$ is the multiplicity of $\BB(V_\mu)$ in
$\BB_{\lambda'}$. Note that the restriction of $\tilde \Delta_X$
to each ${\mathcal C}^{(\mu)}_{\lambda'}$ is a well-defined
function ${\mathcal C}^{(\mu)}_{\lambda'}\to \ZZ$.

For each $\lambda'\in X_\star(S)$ and $\mu\in X_\star(Z(L_P))\cap
X_\star(T)^+$ such that $|{\mathcal C}^{(\mu)}_{\lambda'}|<\infty$, define the
$q$-{\it multiplicity} function $[{\mathcal C}^{(\mu)}_{\lambda'}]_q$
by
$$[{\mathcal
C} _{\lambda'}^{(\mu)}]_q=\sum_{\tilde c\in {\mathcal C}^{(\mu)}_{\lambda'}}
 q^{\tilde \Delta_X(\tilde c)}$$

We apply this construction in the case when
$$({\bf X},{\bf P},f,\Theta_X)=(Bw_0B,\id,f,\Theta_B^-)*(Bw_0B,
\id,f,\Theta_B^-)*\cdots*(Bw_0B,\id,f,\Theta_B^-)$$ is the $k$-th
power of the standard strongly positive $(U\times
U,\chi^{st})$-linear bicrystal $(Bw_0B,\id,f,\Theta_B^-)$, and
$\kappa: {\bf X}\to T^k$ is given by
$$\kappa(x_1*x_2\cdots *x_k)=(hw(x_1),hw(x_2), \ldots, hw(x_k)) \ ,$$
where $hw:Bw_0B\to T$ is the highest weight morphism (see
\eqref{eq:hw_X} above) given by $hw(ut\overline {w_0}u')=t$ for all
$u,u'\in U$. Tropicalizing with respect to
$\Theta_X:=\Theta_B^-*\Theta_B^-*\cdots * \Theta_B^-$, one
obtains for each $\lambda'=(\lambda_1,\lambda_2,\ldots,\lambda_k)\in
(X_\star(T)^+)^k=X_\star(T^k)^+$:
$$\BB_{(\lambda_1,\lambda_2,\ldots,\lambda_k)}=\BB(V_{\lambda_1}\otimes V_{\lambda_2}\otimes
\cdots \otimes V_{\lambda_k})=\BB(V_{\lambda_1})\otimes \BB(V_{\lambda_2})\otimes \cdots \otimes \BB(V_{\lambda_k}) \ .$$

Therefore, the polynomial $[{\mathcal
C}^{(\mu)}_{(\lambda_1,\lambda_2,\ldots,\lambda_k)}]_q$ is a new
$q$-deformation of the tensor product multiplicity
$[V_\mu:V_{\lambda_1}\otimes V_{\lambda_2}\otimes \cdots \otimes
V_{\lambda_k}]$.

We expect that the polynomials $[{\mathcal
C}^{(\mu)}_{(\lambda_1,\lambda_2,\ldots,\lambda_k)}]_q$ are related
to parabolic Kazhdan-Lusztig polynomials (see \cite{lt}).


\end{document}